\documentclass[reqno]{amsart}
\usepackage[utf8]{inputenc}

\usepackage[margin=3.5cm]{geometry}
\usepackage[english]{babel}
\usepackage[pdftex]{graphicx}
\usepackage{amsfonts}
\usepackage{amsmath, bm}
\usepackage{amsbsy}
\usepackage{amsthm}
\usepackage{subcaption}
\usepackage{mathtools}
\usepackage[dvipsnames]{xcolor}
\usepackage{gensymb}
\usepackage{textcomp}
\usepackage{cleveref}
\usepackage[utf8]{inputenc}
\usepackage{makecell}
\usepackage[square,sort,comma,numbers]{natbib}

\DeclareFontFamily{U}{skulls}{}
\DeclareFontShape{U}{skulls}{m}{n}{ <-> skull }{}

\mathtoolsset{centercolon}

\newtheorem{theorem}{Theorem}

\newtheorem{proposition}[theorem]{Proposition}
\newtheorem{remark}[theorem]{Remark}

\makeatletter
\renewcommand*\env@matrix[1][*\c@MaxMatrixCols c]{%
  \hskip -\arraycolsep
  \let\@ifnextchar\new@ifnextchar
  \array{#1}}
\makeatother

\long\def\MSC#1\EndMSC{\def\arg{#1}\ifx\arg\empty\relax\else
{\narrower\noindent%
{2010 Mathematics Subject Classification}: #1\\} \fi}
\long\def\PACS#1\EndPACS{\def\arg{#1}\ifx\arg\empty\relax\else
{\narrower\noindent%
{PACS numbers}: #1}\fi}
\long\def\KEY#1\EndKEY{\def\arg{#1}\ifx\arg\empty\relax\else
{\narrower\noindent%
Keywords: #1\\}\fi}

\renewcommand{\v}[1]{ \textbf{#1} }

\newcommand{\trace}{ {\operatorname{trace}} }

\providecommand{\trace}{{\operatorname{trace}}}

\providecommand{\R}{{\mathbb{R}}}

\providecommand{\N}{{\mathbb{N}}}

\providecommand{\vA}{{\bm{A}}}
\providecommand{\vB}{{\bm{B}}}
\providecommand{\vC}{\bm{C}}

\providecommand{\vF}{\bm{F}}

\providecommand{\vI}{\bm{I}}

\providecommand{\vH}{\bm{H}}
\providecommand{\vK}{\bm{K}}
\providecommand{\vL}{\bm{L}}
\providecommand{\vM}{\bm{M}}
\providecommand{\vN}{\bm{N}}
\providecommand{\vQ}{\bm{Q}}

\providecommand{\vU}{\bm{U}}
\providecommand{\vV}{\bm{V}}

\providecommand{\vX}{{\bm{X}}}
\providecommand{\vY}{{\bm{Y}}}

\providecommand{\vf}{{\bm{f}}}

\providecommand{\vp}{{\bm{p}}}

\providecommand{\vu}{{\bm{u}}}

\providecommand{\vx}{{\bm{x}}}
\providecommand{\vy}{{\bm{y}}}

\providecommand{\vGamma}{{\bm {\mathit{\Gamma}}}}

\providecommand{\vPsi}{{\bm {\mathit{\Psi}}}}
\providecommand{\vSigma}{{\bm {\mathit{\Sigma}}}}

\providecommand{\noise} {{ \mathrm{noise} }}

\providecommand{\bdry} {{ \mathrm{bdry} }}

\providecommand{\Ltwo} {{L^2(\Omega)}}

\usepackage{natbib}

\begin{document}

\title[Inverse source problem for iron loss]{Inverse Heat Source Problem and Experimental Design for Determining Iron Loss Distribution}

\author[A.~Hannukainen]{Antti Hannukainen}
\address[A.~Hannukainen]{Department of Mathematics and Systems Analysis, Aalto University, P.O. Box~11100, 00076 Helsinki, Finland.}
\email{antti.hannukainen@aalto.fi}

\author[N.~Hyv\"onen]{Nuutti Hyv\"onen}
\address[N.~Hyv\"onen]{Department of Mathematics and Systems Analysis, Aalto University, P.O. Box~11100, 00076 Helsinki, Finland.}
\email{nuutti.hyvonen@aalto.fi}

\author[L.~Perkki\"o]{Lauri Perkki\"o}
\address[L.~Perkki\"o]{Department of Mathematics and Systems Analysis, Aalto University, P.O. Box~11100, 00076 Helsinki, Finland.}
\email{lauri.perkkio@aalto.fi}

\date{\today}

\begin{abstract}
Iron loss determination in the magnetic core of an electrical machine, such as a motor or a transformer, is formulated as an inverse heat source problem. The sensor positions inside the object are optimized in order to minimize the uncertainty in the reconstruction in the sense of the A-optimality of Bayesian experimental design. This paper focuses on the problem formulation and an efficient numerical solution of the discretized sensor optimization and source reconstruction problems. A semirealistic linear model is discretized by finite elements and studied numerically.
\end{abstract}

\maketitle

\KEY
Bayesian inversion,
electric machine,
inverse source problem,
iron loss,
optimal experimental design
\EndKEY

\MSC
65N21, 62K05, 35K20
\EndMSC

\section{Introduction}
A dynamic electromagnetic field induces heat generation in the core materials of electric machines, such as transformers or electric motors. These unwanted phenomena are called iron losses, and they constitute a major portion of the total power loss in an electrical machine. The iron loss depends on the electromagnetic field via complicated mechanisms, and there exist several different models that attempt to estimate such losses~\cite{Kowal2015}. However, the validity of these models cannot be verified directly; instead, it has to be examined indirectly via temperature and calorimetric measurements.

In the approach chosen in this paper, the iron loss acts as an unknown (volume) source field in the heat equation, and this source is reconstructed by measuring the temperature on an easily reachable surface of the machine as well as at a limited number of sensors inside the machine. In other words, we consider an {\em inverse heat source problem} with a time-independent source term, which has been studied both in theory and in practice for some decades~\cite{Cannon1968, Engl1994, Hao2017}. To the best of our knowledge, the inverse source problem approach to the iron loss determination has been studied only recently (e.g.,~\cite{Krings2012, Nair2017}). 

The first aim is to investigate the overall feasibility of obtaining a good reconstruction of the unknown source. As the studied inverse problem is severely ill-posed, the reconstruction is extremely sensitive to measurement noise and model errors, and so the problem has to be regularized or treated statistically. Based on our numerical tests with simulated data, a boundary measurement (thermal camera) has to be augmented by sensors inside the object to obtain information on the source in the most crucial area close to the the windings of the electric machine. The number of these sensors is limited for practical reasons, so a proper sensor placement is studied in the Bayesian framework.  { We consider the optimal experimental design (OED) applied to the iron loss determination problem, and especially some of the introduced computational tools, to be the main novelty of this paper.}

\begin{figure}
	\includegraphics[width=0.54\columnwidth]{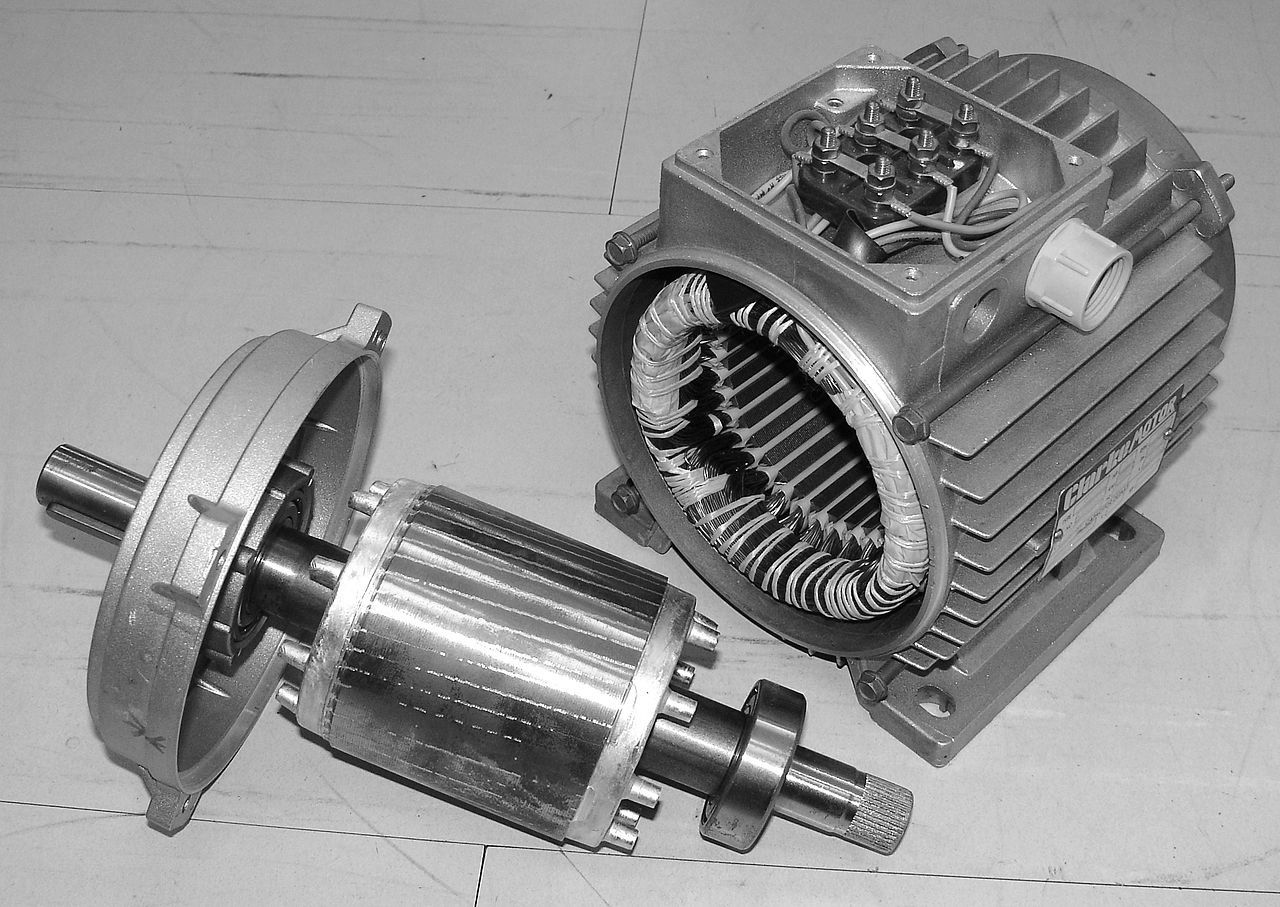}
	\includegraphics[width=0.40\columnwidth, clip=true, trim = 0 0.8cm 0 0]{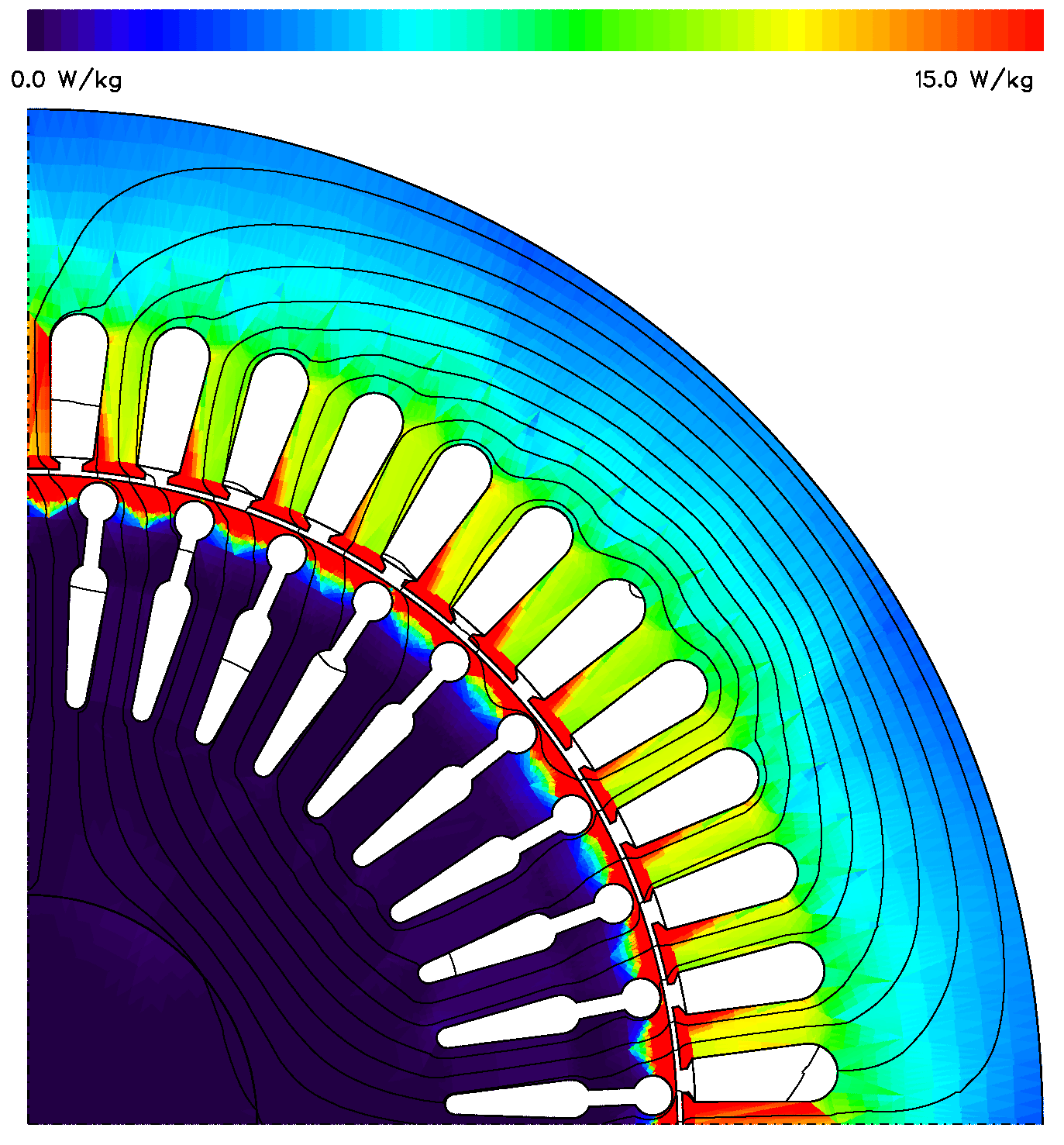}
	\caption{\small {\sc Left}: A rotating electrical machine. {\sc Right}: Iron loss (i.e., heat source) in the machine simulated by electromagnetic FE analysis~\cite{Rasilo2012}.}
	\label{fig:arkkio}
\end{figure}

Figure~\ref{fig:arkkio} shows an estimated heat loss distribution computed by electromagnetic {\em finite element} (FE) analysis, using an existing heat loss model~\cite{Rasilo2012}. The loss is expected to be a smooth function that takes large values and varies quickly close to the windings and decays towards the outer boundary of the machine. This general information could in principle be included in the reconstruction process, but in order not to bias our results by assuming too much prior information on the source, we exclude such considerations in this preliminary study. However, in some numerical tests we assume the variations in the source are correlated with the anisotropy in the structure of the examined machine.

A detailed heat model for a rotating machine is nontrivial in general, as it involves an air (or other coolant) flow in a complicated geometry~\cite{Yoheswaran2014}. In addition, the end windings transfer a significant amount of heat out of the machine, and this phenomenon is not easy to model~\cite{Boglietti2007}. Thus, an electric transformer, having a considerably simpler heat conduction model, is studied as a test problem in this work.

The forward problem, namely a linear parabolic initial/boundary value problem, is discretized spatially by FEs, and temporally by a suitable implicit difference method. The resulting system involves a large number of degrees of freedom, if a realistic three-dimensional geometry is considered. As the sensor location optimization requires repetitive forward solutions as well as evaluating traces of related posterior covariance matrices with dimensions equaling the number of degrees of freedom in the parametrization for the heat source, a main focus of this paper is on efficient computational tools.

Considering an electric machine as in Figure~\ref{fig:arkkio}, the inverse heat source problem has previously been treated by reduced models, such as thermal networks~\cite{Nair2017}. In such models, the source field is not an `arbitrary' function, but can be understood as a piecewise constant function in different subregions (e.g., in windings, teeth, inner core, outer core), leading to a low number of parameters to be solved in the inversion. In some sense, this `model order reduction' gives a suitable regularization `by discretization' for the inverse heat source problem \cite{Kirsch11}. In contrast, the source is reconstructed as a (FE-discretized) continuous function in this paper. 

This text is organized as follows. Section~\ref{sec:formulation} introduces the continuum forward and inverse problems, and briefly studies their unique solvability. Section~\ref{sec:discretization} describes the discretization, the basic principles of Bayesian inversion and experimental design, as well as the employed model order reduction strategy. The algorithms and computational techniques for optimizing the internal sensor positions are then introduced in Section~\ref{sec:sensoropt}. Finally, Section~\ref{sec:numresults} illustrates numerical results for simulated data, with a semirealistic geometry and parameters for the considered parabolic {\em partial differential equation} (PDE). An appendix describes the adjoint of the forward problem and a functional derivative that is used in the sensor location optimization.

\section{Setting}
\label{sec:formulation}

In this section, we first describe the idealized parabolic model with continuum boundary measurements for the heat loss in an electric machine and then consider the unique solvability of the associated inverse source problem. The section is completed by introducing a more realistic model for the boundary measurements by a heat camera as well as for the internal heat sensors.

\subsection{Forward model}
\label{sec:fwd_model}
This paper focuses on a linear parabolic PDE, supposedly capturing the essential properties of a practical iron loss determination problem. The bounded physical domain $\Omega \subset \R^d$, $d = 2$ or $3$, is assumed to have a Lipschitz boundary $\partial \Omega = \overline{\Gamma}_{\mathrm R} \cup \overline{\Gamma}_{\mathrm N}$, where the relatively open subsets $\Gamma_{\mathrm R}, \Gamma_{\mathrm N} \subset \partial \Omega$ are such that $\Gamma_{\mathrm R} \cap \Gamma_{\mathrm N} = \emptyset \not= \Gamma_{\mathrm R}$. Here, $\Gamma_{\mathrm R}$ is the conducting boundary and $\Gamma_{\mathrm N}$ is the insulated (or symmetry) boundary. The measurement time interval $[0, T]$, with $T > 0$, can be chosen to be as long as needed.

The time-dependent temperature $u : \Omega \times (0, T) \to \R$ is governed by a linear parabolic initial/boundary value problem with Robin/Neumann boundary conditions:
\begin{subequations} \label{eq:fwdPDE}
\begin{align}
	\rho \partial_t u - \nabla \cdot (\kappa \nabla u) &= f 	& \text{ in } \Omega \times (0, T),  \label{eq:heateq} \\
	\nu \cdot \kappa \nabla u &= h( u_{\mathrm{out}} - u) 	& \text{ on } \Gamma_{\mathrm R} \times (0, T),  \label{eq:robin}\\
	 \nu \cdot \kappa \nabla u &= 0 	& \text{ on } \Gamma_{\mathrm N} \times (0, T),  \label{eq:neumann}\\
	u &= u_{\mathrm{init}} 	& \text{ on } \Omega \times \{ t = 0 \}, \label{eq:initcond}
\end{align}
\end{subequations}
where $\nu \in L^\infty(\partial \Omega, \R^d)$ is the exterior unit normal of $\partial \Omega$. Moreover, $\rho \in L^\infty_+(\Omega)$ is the product of the material density and the heat capacity, $\kappa \in [L^\infty_+(\Omega)]^{d \times d}$ is the heat conductivity, $f \in L^2(\Omega)$ is the time-independent heat source, $h \in L_{+}^\infty(\Gamma_{\mathrm R})$ is the boundary heat transfer coefficient, $u_{\mathrm{out}} \in L^2(\Gamma_{\mathrm R})$ is the ambient temperature, and $u_{\mathrm{init}} \in \Ltwo$ is the initial temperature. Here,
$$
L^\infty_+(D) := \{ v \in L^\infty(D) \ | \ {\rm ess} \inf v > 0 \} \qquad {\rm with} \ D = \Omega \ {\rm or} \ \Gamma_{\mathrm R},
$$
and analogously, the elements of $[L^{\infty}_+(\Omega)]^{d \times d}$ are symmetric matrices with coefficients in $L^{\infty}(\Omega)$ and with a positive essential infimum for the smallest eigenvalue. The multiplier field for all considered function spaces is~$\R$. Without too severe loss of generality, we assume that $u_{\mathrm{init}} \equiv 0 \equiv u_{\mathrm{out}}$.

 For our purposes, it is convenient to introduce a weighted inner product for $L^2(\Omega)$ through
\begin{equation}
  \label{eq:rho_inner}
(w, v)_\rho := \int_{\Omega} \rho w v \, {\rm d} x, \qquad w,v \in L^2(\Omega).
\end{equation}
As $\rho \in L^\infty_+(\Omega)$, this new inner product does not alter the topology of $L^2(\Omega)$. We denote $L^2(\Omega)$ equipped with this new inner product by $L^2_\rho(\Omega)$ in order to remind the reader about the appropriate interpretation of orthogonality. Moreover, we define the scaled heat source by $f_\rho := f/\rho \in L^2_\rho(\Omega)$.

The solution operator for the elliptic (steady-state) part of \eqref{eq:fwdPDE},
\begin{equation*}
\label{eq:Koper0}
  K :
  \left\{
  \begin{array}{l}
    f_\rho \mapsto w, \\[1mm]
  L^2_\rho(\Omega) \to H^1(\Omega),
\end{array}
  \right.
    \end{equation*}
is defined by the following problem: Given $f_\rho \in L^2_\rho(\Omega)$, find $w = K f_\rho $ such that
\begin{equation}
	a(w, v) = (f_\rho,  v )_\rho \qquad \text{ for all } v \in H^1(\Omega),
	\label{eq:Koper}
\end{equation}
where
\begin{align}
	 a(w, v) &:= \int_{\Omega} \kappa  \nabla w \cdot \nabla v \, {\rm d} x + \int_{\Gamma_{\mathrm R}} h  w  v \, {\rm d} S, \qquad w, v \in H^1(\Omega),
    \label{eq:bilinform}
\end{align}
and the boundary integral is understood in the sense of traces. Since the bilinear form $a(\cdot, \cdot): H^1(\Omega) \times H^1(\Omega) \to \R$ is symmetric, bounded and coercive (see,~e.g.,~\cite{Nazarov15}), the unique solvability of \eqref{eq:Koper} as well as the boundedness of $K: L^2_\rho(\Omega) \to H^1(\Omega)$ follows immediately from the Lax--Milgram theorem. The reason for defining~\eqref{eq:Koper} and $K$ with the help of the inner product of $L^2_\rho(\Omega)$ will become more apparent when the unique solvability of our (idealized) inverse problem is tackled in Section~\ref{sec:measmodel}.

The solution operator for the complete time-dependent problem~\eqref{eq:fwdPDE}, mapping the time-independent source $f$ to the time-dependent temperature $u$, is denoted by
\begin{equation*}
  \mathit{\Psi}:
\left\{
  \begin{array}{l}
    f \mapsto u, \\[1mm]
    L^2(\Omega) \to \mathcal{H}^1\big((0, T); \Omega \big),
  \end{array}
  \right.
\end{equation*}
where
\begin{equation}
  \mathcal{H}^s\big((0, T); \Omega\big) := H^{s-1}\big((0, T) ; H^1(\Omega) \big) \, \cap
  \, H^s\big((0, T) ; H^1(\Omega)^* \big) , \qquad s \in \R.
    \label{eq:Uspace}
\end{equation}
In particular, $\mathcal{H}^1\big((0, T); \Omega\big) \subset \mathcal{C} ( [0, T] ; L^2(\Omega) )$ \cite[Chapter~10]{Renardy93}. To be more precise, \eqref{eq:fwdPDE} is interpreted in a weak sense: the temperature $u \in \mathcal{H}^1((0, T); \Omega)$ satisfies the variational equation
\begin{align}
	 \langle \partial_t u, v \rangle_\rho +  a(u, v) &=   ( f_\rho,  v )_\rho \qquad \text{ for all } v \in H^1(\Omega) \text{ and almost all } t \in (0, T),
	\label{eq:weakform}
\end{align}
together with the initial condition of \eqref{eq:fwdPDE}.
Here $\langle \, \cdot \, , \, \cdot \, \rangle_\rho: H^1(\Omega)^* \times H^1(\Omega) \to \R$ is the dual pairing between $H^1(\Omega)$ and its dual $H^1(\Omega)^*$ with $L^2_\rho(\Omega)$ as the pivot space, that is, $\langle w, v \rangle_\rho = ( w,  v )_\rho$ if $w, v \in L^2(\Omega)$.  It follows from the standard theory for parabolic PDEs that \eqref{eq:weakform} has a unique solution in $\mathcal{H}^1\big((0, T); \Omega\big)$~\cite[Chapter~10]{Renardy93}.

As $\kappa$ and $\rho$ are independent of time and $u_{\rm init} \equiv 0$, the unique solution of \eqref{eq:weakform} in fact carries more time-regularity and belongs to $\mathcal{H}^2((0, T); \Omega)$ \cite[Chapter~10]{Renardy93}. In particular, the fixed-time operator
\begin{equation}
  \label{eq:Psit}
  \mathit{\Psi}(t):
\left\{
  \begin{array}{l}
    f \mapsto u(\,\cdot \, , t), \\[1mm]
    L^2(\Omega) \to H^1(\Omega),
  \end{array}
  \right.
\end{equation}
is well defined for all $t \in [0,T]$. Observe that the extra time-regularity also means $\partial_t u \in \mathcal{H}^1((0, T); \Omega) \subset \mathcal{C} ( [0, T] ; L^2(\Omega) )$, and thus the dual bracket in \eqref{eq:weakform} can actually be interpreted as the inner product of $L^2_\rho(\Omega)$.

\begin{remark}
  In the considered inverse problem, the heat source $f$ is the unknown, and the aim is to reconstruct it from limited temperature measurements. The assumption $f \in L^2(\Omega)$ in~\eqref{eq:fwdPDE} provides (more than) enough regularity for the forward problem to be well defined, but in the inversion one should assume more prior information on $f$, as discussed in Section~\ref{sec:BayOpt} below. The other parameters in~\eqref{eq:fwdPDE}, $ \rho$, $\kappa$, $h$, as well as $u_{\rm init}\equiv 0 \equiv  u_{\rm out}$, are assumed to be precisely known in our considerations.

  A more complete model for the inverse problem would  consider the boundary heat transfer coefficient $h$ as a second unknown, as it cannot typically be measured reliably in practice~\cite{Staton2008}. In addition, the parameters and the heat source $f$ are in reality temperature dependent, which means that a more accurate model for the evolution of the temperature inside $\Omega$ would involve a nonlinear parabolic PDE. However, taking these observations into account would lead to a more complicated nonlinear inverse source problem, which we wish to avoid in this initial study.
\end{remark}

\subsection{Idealized inverse problem and its unique solvability}
\label{sec:measmodel}
Let us start by considering the inverse problem with a continuum of boundary measurements without any interior sensors. To be more precise, we assume the temperature $u$ can be measured on a certain nonempty relatively open part $\mathcal M_\bdry \subset \partial \Omega$ of the object boundary over the time interval $(0,T)$. The idealized forward operator, sending the unknown heat source to the boundary measurement, is thus defined via
\begin{equation}
  \label{eq:forwardop}
\mathcal{F}: \left\{
  \begin{array}{l}
    f \mapsto u|_{\mathcal{M}_\bdry}, \\[2mm]
  L^2(\Omega) \to L^2\big((0,T); L^2(\mathcal{M}_\bdry) \big) \cong  L^2\big( \mathcal{M}_\bdry \times (0,T)\big),
\end{array}
  \right.
\end{equation}
  where $u \in \mathcal{H}^1((0,T); \Omega)$ is the solution to \eqref{eq:weakform}. Observe that for all $t \in (0,T)$,
  $$
  \| u|_{\mathcal{M}_\bdry}(\,\cdot\, , t) \|_{L^2(\mathcal{M}_\bdry)} \leq C \| u(\,\cdot\, , t )\|_{H^1(\Omega)}
  $$
  by the trace theorem. Hence,
\begin{align*}
  \| \mathcal{F}f\|_{L^2(\mathcal{M}_\bdry \times (0,T))} &=
  \| u|_{\mathcal{M}_\bdry} \|_{L^2((0,T); L^2(\mathcal{M}_\bdry)) } \\[1mm]
   & \leq C \| u \|_{L^2((0,T); H^1(\Omega))} \leq C \| f \|_{L^2(\Omega)},
\end{align*}
where the last step follows from the continuity of the solution operator $\mathit{\Psi}$ for \eqref{eq:weakform} as a mapping from $L^2(\Omega)$ to $L^2((0,T); H^1(\Omega))$~\cite[Chapter~10]{Renardy93}. In other words, $\mathcal{F}$ is well defined and bounded.

Our idealized inverse problem is defined as follows:
\begin{align}
  \textrm{Given } g \in L^2( \mathcal{M}_\bdry \times (0,T)),  \textrm{ find } f \in L^2(\Omega) \, \text { such that } \, \mathcal{F} f = g.
    \label{eq:ideal}
\end{align}
It is obvious that \eqref{eq:ideal} does not have a solution for all --- or actually most --- $g \in L^2( \mathcal{M}_\bdry \times (0,T))$. However, if there exists a solution, it is unique under only mild regularity assumptions on $\kappa$. More details can be found in~\cite{Engl1994}, where similar uniqueness results are proven for a more general $f$, but a smoother $\kappa$ in~\eqref{eq:heateq}.
\begin{proposition}
Assume $\kappa \in [L^\infty_+(\Omega)]^{d \times d}$ is regular enough to allow unique continuation of Cauchy data from $\mathcal{M}_\bdry$ to $\Omega$ for the elliptic steady-state equation \eqref{eq:Koper}. Then,  $\mathcal{F}: L^2(\Omega) \to L^2( \mathcal{M}_\bdry \times (0,T))$ is injective.
    \label{lemma:ideal}
\end{proposition}
\begin{proof}
It is easy to check that the steady-state operator $K$, defined by~\eqref{eq:Koper}, is self-adjoint and positive definite when interpreted as an operator from $L^2_\rho(\Omega)$ to itself. It is also compact due to the compactness of the embedding $H^1(\Omega) \hookrightarrow L^2(\Omega)$ and the equivalence of the topologies of $L^2(\Omega)$ and $L^2_\rho(\Omega)$. Hence, it follows from fundamental spectral theory that $K$ has an orthonormal eigenbasis $\{ v_i \}_{i=1}^\infty$ for $L^2_\rho(\Omega)$ with the corresponding eigenvalues $\{\lambda_i\}_{i=1}^\infty \subset \R_+$ repeated according to their multiplicity and satisfying $\lambda_i \to 0$ as $i \to \infty$. Since $v_i = \lambda_i^{-1} K v_i$, it follows that $v_i \in H^1(\Omega)$ for all $i \in \N$. \cite{Yosida95}

If the scaled source $f_\rho = f/\rho \in L^2_\rho(\Omega)$ is expanded in the aforementioned eigenbasis as $f_\rho = \sum_{i=1}^\infty  (f_\rho,v_i)_\rho v_i$, then the solution to~\eqref{eq:weakform} reads
\[
    u(\,\cdot \, , t) = \sum_{i=1}^\infty \lambda_i (f_\rho,v_i)_\rho \big(1 - {\rm e}^{-t/\lambda_i}\big) v_i, \qquad t \in [0,T],
\]
    as can  be verified via a straightforward calculation. In particular, this leads to the representation
    \begin{equation}
      \label{eq:spectral}
    \mathcal{F}: f \mapsto \sum_{i=1}^\infty \lambda_i (f_\rho,v_i)_\rho \big(1 - {\rm e}^{-t/\lambda_i}\big) v_i|_{\mathcal{M}_\bdry}, \qquad t \in [0,T],
    \end{equation}
    for any $f \in L^2(\Omega)$.

    Let $\{ \widehat{\lambda}_i \}_{i=1}^\infty \subset \R_+$ denote the {\em distinct} eigenvalues of $K$, let $\{ \mu_i \}_{i=1}^\infty \subset \N$ be their respective multiplicities, and let $\widehat{v}_{i,j}$, $j=1, \dots, \mu_i$, form an $L^2_\rho(\Omega)$-orthonormal basis for the eigenspace corresponding to $\widehat{\lambda}_i$. It immediately follows from \eqref{eq:spectral} that for almost all $x \in \mathcal{M}_\bdry$ and all $t \in [0,T]$, the output $\mathcal{F}f$ can be given as a (generalized) {\em Dirichlet series} of the form (see,~e.g.,~\cite[Chapter~VI, Section~2]{Serre1973})
\begin{equation}
    \label{eq:dirichlet_series}
    (\mathcal{F}f)(x,t) = a_0(x) + \sum_{i=1}^\infty a_i(x) \, {\rm e}^{-t/\widehat{\lambda}_i},
    \end{equation}
    where
    $$
    a_0 = \sum_{i=1}^\infty \lambda_i (f_\rho,v_i)_\rho v_i|_{\mathcal{M}_\bdry} \qquad {\rm and} \qquad a_i = - \widehat{\lambda}_i \sum_{j=1}^{\mu_i}  (f_\rho,\widehat{v}_{i,j})_\rho  \widehat{v}_{i,j}|_{\mathcal{M}_\bdry}.
    $$
Observe that \eqref{eq:dirichlet_series} really is a Dirichlet series since $\widehat{\lambda}_{i}^{-1} \to \infty$ as $i \to \infty$, and hence it defines a holomorphic function in $t$ in the whole open right half of the complex plane for almost every fixed $x \in \mathcal{M}_\bdry$~\cite[Chapter~VI, Corollary~1]{Serre1973}.

Suppose now $\mathcal{F}f = 0$ for some $f \in L^2(\Omega)$, which via analytic continuation means that $(\mathcal{F}f)(x,\, \cdot \,)$ given by \eqref{eq:dirichlet_series} equals zero in the right half of the complex plane for almost all $x \in \mathcal{M}_\bdry$. Due to the uniqueness of the coefficients in a Dirichlet series \cite[Chapter~VI, Corollary~4]{Serre1973}, it must thus hold $a_i = 0$, $i \in \N_0$, almost everywhere on $\mathcal{M}_\bdry$. Since $\widehat{v}_{i,j}$ are eigenfunctions of the corresponding steady-state operator $K$ and, in particular, satisfy the boundary conditions of \eqref{eq:fwdPDE} in the appropriate sense of traces,
\[
w_i :=  \sum_{j=1}^{\mu_i}  (f_\rho,\widehat{v}_{i,j})_\rho  \widehat{v}_{i,j} \in H^1(\Omega), \qquad i \in \N,
\]
weakly satisfies
\[
  - \nabla \cdot (\kappa \nabla w_i) = \widehat{\lambda}_i^{-1} w_i  \quad \text{in} \ \Omega, \qquad
  -a_i/\widehat{\lambda}_i =  w_i = \nu \cdot \kappa \nabla w_i = 0  \quad \text{on} \ \mathcal{M}_\bdry.
\]
The principle of unique continuation thus yields $w_i = 0$ in $\Omega$ for all $i \in \N$. As $\{ v_i\} = \{ \widehat{v}_{i,j}\}$ is an orthonormal basis for $L^2_\rho(\Omega)$, it must thus, in fact, hold $f/\rho = f_\rho = 0$. Since $\rho \in L^\infty_+(\Omega)$, the proof is complete.
\end{proof}

\begin{remark}
  The exact smoothness requirement on $\kappa$ for \eqref{eq:Koper} to allow unique continuation of Cauchy data is dimension dependent. However, in all spatial dimensions (suitably defined) piecewise `smooth enough' regularity of $\kappa$ is both a sufficient theoretical condition as well as a reasonable assumption from the practical standpoint (cf.,~e.g.,~\cite{Druskin1998}).
  \end{remark}

\begin{remark}
According to Proposition~\ref{lemma:ideal}, the yet-to-be-introduced inner sensors are not required for the unique solvability of the idealized inverse problem~\eqref{eq:ideal}. However, they {\em significantly} improve numerical reconstructions, as demonstrated by the numerical examples in Section~\ref{sec:numresults}.
\end{remark}

\subsection{Discrete measurements and interior sensors}
\label{sec:disc_sensors}
In practice, the temperature measurements on the accessible boundary $\mathcal{M}_\bdry$ are obtained by a (digital) thermal camera. To augment the reconstruction process, a small number of finite-size temperature sensors can be placed inside the domain $\Omega$ by installing them onto a circuit board, and subsequently inserting it between the iron lamination sheets of the examined transformer core.

According to a more realistic model, a single measurement corresponds to the mean temperature at a given time over $S_i$, where $S_i \subset \mathcal{M}_\bdry$ (or $S_i \subset \Omega$) is identified with the considered heat camera pixel (or interior sensor). To be more precise, the whole spatial measurement is modeled by the finite-dimensional mapping
\[
B = \left[ \begin{array}{c}
    \! B_1 \! \\
    \vdots \\
    \! B_{m_s} \!
  \end{array}
  \right]: \ H^1(\Omega) \to \R^{m_s},
\]
with
\begin{equation}
	B_{i} w  = \dfrac{1}{ | S_i | }  \int_{S_i} w \, {\rm d} S, \qquad i=1, \dots, m_s.
	\label{eq:pointmeas}
\end{equation}
Here $| S_i |$ denotes the area (or volume) of the $i$th sensor and $m_s$ is the total number of pixels and sensors. Note that the integral on the right-hand side of \eqref{eq:pointmeas} is well defined for all $i = 1, \dots, m_s$ due to the trace theorem.

Let us denote by $t_j$, $j = 1, \dots, m_s$, the discrete measurement times and set $m=m_s m_t$. The {\em realistic} forward map
\begin{equation}
  \label{eq:Fdisc}
F = \left[ \begin{array}{c}
    \! B \mathit{\Psi} (t_1) \! \\
    \vdots \\
    \! B \mathit{\Psi}(t_{m_t}) \!
  \end{array}
  \right]: \ L^2(\Omega) \to \R^{m},
\end{equation}
is well defined, linear, bounded and compact; see \eqref{eq:Psit}. To sum up, the realistic inverse problem that we aim to numerically tackle reads:
\begin{align}
  \textrm{Given } y \in \R^{m},   \textrm{ find } f \in L^2(\Omega) \, \text { such that } \, F f = y.
    \label{eq:inverse2}
\end{align}
In particular, due to the finite-dimensionality of $F$, \eqref{eq:inverse2} is obviously not as such uniquely solvable for any $y \in\R^{m}$.

\begin{remark}
  \label{remark:weight}
  It would be physically more realistic to model the measurements as
  \begin{equation}
    \label{eq:weighted_meas}
  \dfrac{1}{|S_i|} \int_0^T \int_{S_i} u(x,t) \omega(t-t_j) \, {\rm d} S \, {\rm d} t, \qquad i=1, \dots, m_s, \ j = 1, \dots, m_t,
  \end{equation}
  where $u$ is the weak solution to \eqref{eq:fwdPDE} and $\omega \geq 0$ is a `device function' that is concentrated around the origin and integrates to one. This would account for the fact that no temperature measurement can be instantaneous, but it actually lasts in reality over a finite time interval. Moreover, from the mathematical standpoint, computing \eqref{eq:weighted_meas} is more stable than applying $F$ of \eqref{eq:Fdisc} to a given source $f \in L^2(\Omega)$. However, for the sake of notational and conceptual simplicity, we stick with the pointwise measurements in time.
  \end{remark}

\section{Discretization and the Bayesian setting}
\label{sec:discretization}

Although it would be possible to formulate the principles of Bayesian optimal experimental design for an infinite-dimensional unknown (cf.~\cite{Stuart2010}), we discretize the forward operator \eqref{eq:Fdisc} before presenting the Bayesian formulation for the inverse problem \eqref{eq:inverse2}. The reason for this choice is that our main objective is to introduce an efficient {\em computational} framework for choosing optimal locations for the interior sensors; in particular, we want to separate this task from any extra complications caused by an infinite-dimensional setup.

\subsection{Discretization of the forward operator}
The temperature $u$ in \eqref{eq:fwdPDE} is discretized spatially by standard $H^1$-FEs and time-integrated by the implicit midpoint rule (unless stated otherwise). The discretization is assumed to be ``good enough'', so that the discretization error is negligible compared to other sources of error related to, e.g., mismodeling and measurement noise. Throughout this text, discretized objects are written in bold.

A Lagrangian FE basis corresponding to a discretization of $\Omega$ is denoted by $\{ \phi_i \}_{i=1}^n \subset H^1(\Omega)$, and the nodal values of the discretized temperature field at the $i$th time step of the implicit midpoint rule is denoted by $\vu_i \in \R^n$. In our numerical examples, the unknown source $f$ is assumed to be an element of $ H^1(\Omega)$, and so it is reasonable to express the discretized source field in the same FE basis as the temperature field, although one could in principle use a coarser discretization or some other lower-dimensional basis for $f$ to reduce the computational cost.

As mentioned in Section~\ref{sec:disc_sensors}, the temperature is measured at $m_s$ locations and $m_t$ observation times, so the measurement can be interpreted as a vector $\vy \in \R^{m}$, with $m = m_s m_t$. The discrete forward operator $\vF \in \R^{m \times n}$ approximates $F$ defined in \eqref{eq:Fdisc}: Given the nodal values of a (discretized) heat source $\vx \in \R^n$, the corresponding simulated measurement is given by
\[
\vy =
\left[ \begin{array}{c}
    \! \vy_{1} \! \\
    \vdots \\
    \! \vy_{m_t} \!
  \end{array}
  \right] =
\left[ \begin{array}{c}
    \! \vF_1 \vx  \! \\
    \vdots \\
    \! \vF_{m_t} \vx \!
  \end{array}
  \right]
    =:  \vF \vx,
\]
where $\vy_i \in \R^{m_s}$ carries the measured temperatures at the $i$th observation time. Moreover,
\begin{align}
  \vF_i := \vB \vPsi_{\! i} \vM, \qquad i=1, \dots, m_t,
  \label{eq:fwd}
\end{align}
where the FE mass matrix $\vM_{i,j} := (\phi_i, \phi_j)_{L^2(\Omega)}$, $i,j = 1, \dots, n$, maps a nodal vector $\vx$ to a load vector $\vf$, $\vPsi_{\! i} \in \R^{n \times n}$ discretizes $\Psi(t_i)$ defined by \eqref{eq:Psit} in the FE basis, and $\vB\in\R^{m_s \times n}$ discretizes the spatial measurement operator~\eqref{eq:pointmeas}.

More specifically, $\vPsi_{\! i} \vf$ gives a numerical solution to the system of ordinary differential equations
\begin{align}
  \label{eq:ODE}
    \vM_\rho  \vu'(t) + \vK \vu(t) &= \vf, \\
    \vu(0) &= {\bf 0}, \nonumber
\end{align}
evaluated at the time $t_i$. The matrices in \eqref{eq:ODE} originate from~\eqref{eq:bilinform} and~\eqref{eq:weakform}, that is,
\begin{align*}
  \left[ \vM_\rho \right]_{i,j} &= \int_{\Omega} \rho \phi_i \phi_j \, {\rm d}x, \\
  \vK_{i, j} &= \int_{\Omega} \kappa \nabla \phi_i \cdot \nabla \phi_j \, {\rm d} x + \int_{\Gamma_{\rm R}} h \phi_i \phi_j \, {\rm d} S,  \qquad i,j = 1, \dots, n. \\
  \vf_{i} &= \int_{\Omega} f \phi_i \, {\rm d} x,
  \end{align*}
The  solution to \eqref{eq:ODE} can either be approximated by a suitable time-integration scheme (in our case the implicit midpoint rule), or given by the explicit formula
$$
\vu(t) = \big(\vI + \exp(- \vM_\rho^{-1} \vK t)\big) \vK^{-1} \vf =: \vPsi_t^{\mathrm{exct}} \vf,
$$
if the problem is sufficiently small so that the matrix exponential can be numerically evaluated.

Initially, one would expect that the explicit construction of $\vF$ requires one numerical solution of~\eqref{eq:ODE} for each FE degree of freedom. However, by noticing that $\vPsi_t^{\mathrm{exct}}$ is a symmetric matrix and the same holds when $\vPsi_t$ corresponds to, e.g., the implicit midpoint rule, $\vF^T$ can be computed cheaply by transposing~\eqref{eq:fwd}:
\begin{align}
	\vF^T = \vM \left[ \vPsi_{\! 1} \vB^T \; \vPsi_{\! 2} \vB^T \; \cdots \; \vPsi_{\! m_t} \vB^T  \right].
	\label{eq:Ftransp}
\end{align}
In other words, $\vF$ can be formed by solving one parabolic forward problem for each sensor location,~i.e., for each column of $\vB^T$; this observation is interpreted in the non-discretized setting in Appendix~\ref{sec:adjoints}.  However, it is more reasonable to handle the boundary measurement corresponding to a large number of sensors/pixels with the help of a low-dimensional approximation, as explained in Section~\ref{sec:lowrank} below.

\subsection{Bayesian inversion and A-optimal design}
\label{sec:BayOpt}
In Bayesian inversion all parameters carrying uncertainty are treated as random variables. The prior probability distributions for these parameters reflect the available information before the measurements are carried out. The measurement is modeled as a realization of a random variable depending on both the noise process and the random parameters in the forward model, as well as on the so-called design parameters that define the measurement setup. The Bayes’ formula is then employed to form the posterior probability density that updates the prior based on the information in the measurement. In our setting, the experimental design variables are the positions of the internal sensors, and our ultimate aim is to choose them so that the posterior density of the heat source is as `localized as possible' in the sense of the A-optimality criterion of Bayesian optimal experimental design. \cite{Chaloner95, Kaipio2004statistical}

Let $\vy \in \R^m$ carry the (noisy) temperature measurements, $\vp \in \R^N$ be a vector parametrizing the positions of the internal sensors, and suppose our prior information on the (discretized) heat source is encoded in a probability density $\pi_{\rm pr}: \R^n \to \R_+$. By the Bayes' formula, the posterior density for the (randomized) nodal source $\vX$ reads~\cite{Kaipio2004statistical}
\begin{equation}
  \label{eq:Bayes}
\pi(\vx \, | \, \vy;  \vp) = \frac{\pi(\vy \, | \, \vx;  \vp) \pi_{\rm pr}(\vx)}{\pi(\vy; \vp)},  \qquad \vx \in \R^n,
\end{equation}
where $\pi(\vy \, | \, \cdot \, ;  \vp): \R^n \to \R_+$ is the so-called likelihood function and the normalizing term in the denominator is the marginal density of the random measurement $\vY$ evaluated at the data $\vy$.

In this work we assume that the prior is Gaussian, i.e.~$\vX \sim \mathcal{N}(\vx_{\mathrm{pr}}, \vGamma_{\!\rm pr})$, and the measurement can be modeled as a realization of the random variable
\[
\vY = \vF(\vp) \vX + \vN,
\]
where $\vN \sim \mathcal{N}(\v0, \vGamma_{\!\rm noise})$ is independent of $\vX$. Here, $\vGamma_{\!\rm pr} \in \R^{n \times n}$ and $\vGamma_{\!\noise} \in \R^{m \times m}$ are symmetric and positive definite covariance matrices, ${\vx}_{\mathrm{pr}} \in \R^n$ is the prior mean for $\vX$, and we have explicitly indicated the nonlinear dependence of the discrete forward operator $\vF(\vp)$ on the positions of the internal sensors{; see Section~\ref{sec:lowrank} for its computational implementation}. Under these simplifying assumptions, the posterior in \eqref{eq:Bayes} is also Gaussian with the covariance matrix and mean~\cite{Kaipio2004statistical}
\begin{subequations}
   \label{eq:posterior_cov}
  \begin{align}
  \vGamma_{\!\rm post}(\vp) &= \big(\vGamma_{\!\rm pr}^{-1} + \vF(\vp)^T \vGamma_{\!\rm noise}^{-1} \vF(\vp) \big)^{-1},\\[1mm]
  \widehat{\vx}(\vp) &= \vGamma_{\!\rm post}(\vp) \big( \vGamma_{\!\rm pr}^{-1}{\vx}_{\mathrm{pr}} + \vF(\vp)^T \vGamma_{\!\rm noise}^{-1} \vy \big),
  \label{eq:postmean}
\end{align}
\end{subequations}
respectively, as can be deduced by a straightforward completion of squares in \eqref{eq:Bayes}. Using the Woodbury matrix identity, these equations can alternatively be represented as (cf.~\cite{Kaipio2004statistical})
\begin{subequations}
   \label{eq:posterior_cov2}
  \begin{align}
    \vGamma_{\!\rm post}(\vp) &= \vGamma_{\!\rm pr} - \vGamma_{\!\rm pr} \vF(\vp)^T \big(\vF(\vp) \vGamma_{\!\rm pr} \vF(\vp)^T + \vGamma_{\!\rm noise}  \big)^{-1} \vF(\vp) \vGamma_{\!\rm pr} ,\\[1mm]
    \widehat{\vx}(\vp) &= {\vx}_{\mathrm{pr}} +  \vGamma_{\!\rm pr} \vF(\vp)^T \big(\vF(\vp) \vGamma_{\!\rm pr} \vF(\vp)^T + \vGamma_{\!\rm noise}  \big)^{-1} (\vy -  \vF(\vp) {\vx}_{\mathrm{pr}}).
      \label{eq:postmean2}
\end{align}
\end{subequations}

In Bayesian optimal experimental design, one often considers minimizing the expected squared distance of the unknown in a given (semi)norm  around some chosen point estimate, which corresponds to the so-called A-optimal design. Assuming the point estimate of interest is the posterior mean and the employed seminorm is induced by the positive semidefinite matrix $\vA^T \!\vA$ for a given $\vA \in \R^{l \times n}$, in our simple, i.e.~Gaussian, linear and finite-dimensional, setting, A-optimality corresponds to choosing a design parameter $\vp_* \in \R^N$ satisfying \cite{Chaloner95,Hyvonen2014}
\begin{equation}
\label{eq:Aoptimal}
\vp_* = {\rm arg} \min_{\vp} \, {\rm tr}  \big(\vA \vGamma_{\!\rm post}(\vp) \vA^T\big) = {\rm arg} \min_{\vp} \, {\rm tr}  \big(\vGamma_{\!\rm post}(\vp) \vA^T \! \vA \big),
\end{equation}
where the second equality is a consequence of the matrix trace being invariant under cyclic permutations. One natural choice for measuring the deviation from the posterior mean is arguably the $L^2(\Omega)$-norm. Since the heat source is represented in the FE basis, one could thus choose $\vA^T \!\vA$ to be the mass matrix $\vM$ associated to the FE discretization. Another possible choice is $\vA^T \!\vA = \vGamma_{\!\rm pr}^{-1}$, which renders the metric for A-optimality to be the same as in the penalty term of the Tikhonov functional corresponding to the `regularized solution' given by the second equation of \eqref{eq:posterior_cov}. As discussed in Section~\ref{sec:eval_object} below, this latter choice makes the formula for (approximately) evaluating the trace needed in \eqref{eq:Aoptimal} particularly simple.

Before investigating the optimal positioning of the internal sensors following the above guidelines, we still need to tackle some computational issues: For a realistic three-dimensional forward problem, the number of degrees of freedom in the parametrization for the unknown source can easily be of the order $n \sim 10^5$, making the repetitive formation of $\vGamma_{\!\rm post}(\vp)$ for different $\vp$ based on \eqref{eq:posterior_cov} impractical (cf.~\eqref{eq:Aoptimal}). If the total number of measurements $m$ is low(ish), this problem can be circumvented by resorting to the alternative formulation \eqref{eq:posterior_cov2}. However, as the number of sensors in our setting is (slightly) higher than the number of pixels in the employed thermal camera, one cannot initially assume that $m$ is of moderate size. Hence, we combine \eqref{eq:posterior_cov2} with a low-rank approximation for $\vF$, as explained in Section~\ref{sec:lowrank} below.

\begin{remark}
  \label{remark:regularization}
  Although most of our analysis only requires $\vGamma_{\!\rm pr} \in \R^{n \times n}$ to simply be positive definite, some choices become more transparent if it is noted that in our numerical experiments, we choose
\begin{equation}
  \label{eq:prior}
	\left[ \vGamma_{\!\rm pr}^{-1} \right]_{i, j} = \int_{\Omega} \big( \beta \phi_i \phi_j + \alpha \nabla \phi_i \cdot \nabla \phi_j \big) {\rm d}x, \qquad i,j = 1, \dots, n,
\end{equation}
where $\beta \in L^\infty_+(\Omega)$ and $\alpha \in [L^\infty_+(\Omega)]^{d \times d}$ are positive (definite) weight functions. In other words, the {\em inverse} covariance matrix corresponds to the FE discretization of an elliptic PDE, and so multiplying with $\vGamma_{\!\rm pr}^{-1}$ is extremely cheap and multiplication by $\vGamma_{\!\rm pr}$ itself is also computationally tractable as it corresponds to solving one elliptic boundary value problem. Note that~\eqref{eq:postmean} and~\eqref{eq:postmean2} essentially correspond to Tikhonov regularization with a mixed $L^2(\Omega)$--$H^1(\Omega)$ penalty term in our numerical studies. { Moreover, there exists no corresponding infinite-dimensional Gaussian random field for a finite-dimensional covariance matrix of the form~\eqref{eq:prior}~\cite{Stuart2010}. See Remark~\ref{remark:trace} for further information on how this is expected to affect the numerical implementation.}
\end{remark}

\subsection{Dimension reduction for $\vF(\vp)$}
\label{sec:lowrank}

To start with, we divide the discretized spatial measurement operator $\vB \in \R^{m_s \times n}$ into two parts as
\[
  \vB(\vp) = \begin{bmatrix} \vB_{\rm int}(\vp) \\[1mm] \vB_{\rm bdry} \end{bmatrix},
\]
where $\vB_{\rm int}(\vp) \in \R^{m_{\rm int} \times n}$ and $\vB_{\rm bdry} \in \R^{m_{\rm bdry} \times n}$ correspond to the internal and boundary sensors, respectively, with $m_{\rm int} + m_{\rm bdry} = m_s$. Typically, $m_{\rm bdry} \gg m_{\rm int}$ since there are far more pixels in a thermal camera image than there are internal sensors. In particular, only the upper part of $\vB(\vp)$ depends on the parameters $\vp \in \R^N$ defining the positions of the internal measurements. { To be more precise, $\vB(\vp)$ is a FE-based discretization of \eqref{eq:pointmeas} for the employed set of internal sensors with predefined geometric specifications and their positions parametrized by $\vp$.}

Let us then abuse the notation by redefining the forward operator $\vF$ as
\[
  \vF(\vp) = \begin{bmatrix} \vF_{\rm int}(\vp) \\ \vF_{\rm bdry} \end{bmatrix},
\]
where
\begin{equation}
\label{eq:FintFbdry}
\vF_{\rm int}(\vp) =  \begin{bmatrix} \vB_{\rm int}(\vp) \vPsi_1 \\ \vdots \\ \vB_{\rm int}(\vp) \vPsi_{m_t}  \end{bmatrix} \! \vM
\qquad \text{and} \qquad
\vF_{\rm bdry} =  \begin{bmatrix} \vB_{\rm bdry} \vPsi_1 \\ \vdots \\ \vB_{\rm bdry} \vPsi_{m_t}  \end{bmatrix} \! \vM.
\end{equation}
In other words, we simply reorder the rows of $\vF$ so that all measurements corresponding to the boundary sensors are at the bottom --- and only the top part of $\vF$ depends on $\vp$.  In what follows, we implicitly assume that all measurement vectors are also ordered in an analogous manner. Moreover, the noise processes corrupting the internal and boundary measurements are assumed to be mutually independent, that is,
\begin{equation}
  \label{eq:independence}
\vGamma_{\!\rm noise} =
\begin{bmatrix}
  \vGamma_{\!\rm int} & 0 \\[1mm]
  0 & \vGamma_{\!\rm bdry}
  \end{bmatrix},
\end{equation}
where $\vGamma_{\!\rm int} \in \R^{m_{\rm int}m_t \times m_{\rm int}m_t }$ and $\vGamma_{\!\rm bdry} \in \R^{m_{\rm bdry}m_t \times m_{\rm bdry}m_t }$ are symmetric and positive definite.

With a realistic measurement noise level, a thermal camera video, with thousands of sensors (i.e.,~pixels) and a high number of observation times, gives a huge number of data points but only a small amount of computationally extractable information on a source (deep) inside $\Omega$. In consequence, $\vF_\bdry$ has a very low (numerical) rank compared to its size. Motivated by this observation, we compute a low-rank approximation for a variant of $\vF_\bdry$. As $\vF_\bdry$ is independent of the inner sensor locations, this precomputed approximation can then be reused in each iteration of the employed algorithm for optimizing the positions of the internal sensors. On the other hand, because well-placed sensors inside the domain are expected to be `more informative' than those on the boundary, $\vF_{\rm int}(\vp)$ is computed exactly by replacing $\vB$ with $\vB_{\rm int}(\vp)$ in~\eqref{eq:Ftransp}.

Low-rank approximations for related Bayesian or Tikhonov-regularized inverse problems are discussed,~e.g.,~in~\cite{Alexanderian2014,Flath2011}. Here we resort to related ideas. We begin by introducing a decomposition
\begin{equation}
	\vGamma_{\!\rm pr}^{-1} = \vL^T \vL,
	\label{eq:chol}
\end{equation}
which is a fill-in reducing sparse Cholesky factorization if a prior such as~\eqref{eq:prior} is used. Since $\vF(\vp)$ is almost always applied in composition with $\vL^{-1}$ in our numerical considerations (cf.~\eqref{eq:posterior_cov4} below), we introduce a `prior-conditioned' forward operator (cf.~\cite{Calvetti05})
\begin{equation}
  \vF^{\rm pr}(\vp) =
  \begin{bmatrix} \vF^{\rm pr}_{\rm int}(\vp) \\[1mm] \vF^{\rm pr}_{\rm bdry} \end{bmatrix} : = \begin{bmatrix} \vF_{\rm int}(\vp)\vL^{-1} \\[1mm] \vF_{\rm bdry}\vL^{-1} \end{bmatrix}
    = \vF(\vp)\vL^{-1}.
	\label{eq:Fpr}
\end{equation}
As noted in Remark~\ref{remark:regularization}, $\vGamma_{\!\rm pr}^{-1}$ originates from a FE discretization of an elliptic PDE in our numerical experiments, and so it makes sense to compute a truncated singular value decomposition for the prior-conditioned matrix $\vF^{\rm pr}_{\rm bdry}$ instead of mere $\vF_{\rm bdry}$. Indeed, because multiplication with (the continuum version of) $\vL^{-1}$ is a smoothening operation, it causes singular values to decay faster and thus enables the use of lower rank approximations for $\vF^{\rm pr}_{\rm bdry}$ than for $\vF_\bdry$ itself.

To be more precise, we introduce an approximation
\begin{equation}
	 \vF^{\rm pr}_{\rm bdry} \approx  \vU_r \vSigma_r \vV^T_r \in \R^{m_{\rm bdry} m_t \times n},
	\label{eq:Fapprox}
\end{equation}
where the columns of $\vV_r \in \R^{n \times r}$ and $\vU_r \in \R^{m_\bdry m_t \times r}$ are orthonormal, $\vSigma_r \in \R^{r \times r}$ is diagonal with positive entries, and $r \in \N$ is small compared to both $m_{\rm bdry} m_t$ and $n$. As forming an (exact) truncated singular value decomposition for $\vF^{\rm pr}_{\rm bdry}$ is expensive, we instead utilize an algorithm for computing an approximate version: A matrix $\vQ \in \R^{m_{\rm bdry} m_t \times r}$ with orthonormal columns effectively spanning the range of $\vF^{\rm pr}_{\rm bdry}$ is computed iteratively. This is accomplished by repeatedly applying $\vF^{\rm pr}_{\rm bdry}$, $(\vF^{\rm pr}_{\rm bdry})^T$ and orthonormalization to a randomly chosen $\vQ_0 \in \R^{ m_{\rm bdry} m_t \times r}$. After an appropriate $\vQ$ is found, a reduced singular value decomposition $\widetilde \vU_r \vSigma_r \vV^T_r$ is computed for the small(ish) matrix $\vQ^T\vF^{\rm pr}_{\rm bdry} \in \R^{r \times n}$  by some conventional method. Finally, the dimension-reduced representation is~\eqref{eq:Fapprox} with $\vU_r = \vQ \widetilde \vU_r$. The reduced dimension $r$ should be chosen to be slightly larger than the (numerical) rank of $\widetilde \vF^{\rm pr}_{\rm bdry}$ to obtain a good convergence rate when forming $\vQ$. Consult~\cite{Halko2011finding} for further details and analysis on the reliability of similar dimension reduction methods.

With the approximation \eqref{eq:Fapprox} in hand, we define
\begin{equation}
  \label{eq:tFpr}
	\widetilde{\vF}^{\rm pr}(\vp) = \begin{bmatrix} \widetilde{\vF}_{\rm int}^{\rm pr}(\vp) \\[1mm]  \widetilde{\vF}^{\rm pr}_{\rm bdry} \end{bmatrix} := \begin{bmatrix} \vF_{\rm int}^{\rm pr}(\vp) \\[1mm] \vSigma_r \vV_r^T \end{bmatrix}
	\in \R^{(m_{\rm int} m_{\rm t} + r) \times n}
\end{equation}
and subsequently utilize \eqref{eq:independence} and \eqref{eq:Fapprox} to write the first equation of \eqref{eq:posterior_cov} approximately as
\begin{equation}
   \label{eq:posterior_cov3}
  \vGamma_{\!\rm post}(\vp) \approx \big(\vGamma_{\!\rm pr}^{-1} +  \vL^T \widetilde{\vF}^{\rm pr}(\vp)^T \tilde{\vGamma}_{\rm noise}^{-1} \widetilde{\vF}^{\rm pr}(\vp) \vL \big)^{-1},
\end{equation}
where
\begin{equation*}
	\widetilde{\vGamma}_{\rm noise}^{-1} = \begin{bmatrix} \vGamma_{\!\rm int}^{-1}  & \bm{0} \\[1mm]
	 \bm{0} & \vU_r^T \vGamma_{\!\rm bdry}^{-1} \vU_r \end{bmatrix} \in \R^{(m_{\rm int} m_{\rm t} + r) \times (m_{\rm int} m_{\rm t} + r)}.
\end{equation*}
By applying the Woodbury matrix identity and \eqref{eq:chol} to \eqref{eq:posterior_cov3}, one finally arrives at
\begin{equation}
  \label{eq:posterior_cov4}
   \vGamma_{\!\rm post}(\vp) \approx \widetilde{\vGamma}_{\!\rm post}(\vp) := \vGamma_{\!\rm pr} -  \vL^{-1}\widetilde\vF^{\rm pr}(\vp)^T \big(\widetilde\vF^{\rm pr}(\vp) \widetilde\vF^{\rm pr}(\vp)^T + \widetilde \vGamma_{\!\rm noise} \big)^{-1} \widetilde\vF^{\rm pr}(\vp) \vL^{-T},
\end{equation}
which is our low-rank approximation for the posterior covariance.

Using \eqref{eq:posterior_cov4}, it is now possible to repetitively and efficiently approximate $\vGamma_{\!\rm post}(\vp)$ while iteratively optimizing the positions of the internal sensors (cf.~\eqref{eq:Aoptimal}), assuming $m_{\rm int} m_{\rm t}$ and $r$ are of moderate size. Indeed, $\vL^{-1}\widetilde\vF^{\rm pr}(\vp)^T$ can be formed explicitly since one essentially only needs to apply $\vGamma_{\!\rm pr} = \vL^{-1}\vL^{-T}$ to the $m_{\rm int} m_t$ rows of $\vF_{\rm int}(\vp)$ and $\vL^{-1}$ to the $r$ columns of $\vV_r$. Moreover, although $\widetilde\vF^{\rm pr}(\vp) \widetilde\vF^{\rm pr}(\vp)^T + \widetilde{\vGamma}_{\!\rm noise} \in \R^{r \times r}$ is full, it is small(ish) and symmetric, and so the associated linear systems can be solved directly.

\begin{remark}
  Let us reiterate that the approximation \eqref{eq:Fapprox} is not affected by the locations of the internal sensors, and thus $\widetilde{\vGamma}_{\! \rm noise}$ and $\widetilde{\vF}^{\rm pr}_{\rm bdry}$ can be precomputed and used, say, in connection with different number or type of internal sensors when optimizing their positions (cf.~\eqref{eq:Aoptimal}).
  \end{remark}

\section{Sensor location optimization}
\label{sec:sensoropt}
Since our ultimate aim is to find an A-optimal measurement design as defined in~\eqref{eq:Aoptimal}, the to-be-minimized target functional is\footnote{In Section~\ref{sec:numerics}, we frequently abuse the notation by denoting the A-optimality target functional as $\mathit{\Phi}_{\!\vA}$ even if the exact posterior $\vGamma_{\!\rm post}$ is used in \eqref{eq:traceminim} in place of its low-rank approximation $\widetilde{\vGamma}_{\!\rm post}$. The choice between the two options should be clear from the context.}
\begin{equation}
	\mathit{\Phi}_{\vA}(\vp) := {\rm tr} \big(\widetilde{\vGamma}_{\!\rm post}(\vp) \vA^T \! \vA\big) = {\rm tr} \big(\vA \widetilde{\vGamma}_{\!\rm post}(\vp) \vA^T\big),
	\label{eq:traceminim}
\end{equation}
where we have replaced the exact posterior covariance by its approximation from \eqref{eq:posterior_cov4} and $\vA^T \! \vA$ defines the seminorm in which we aim to minimize the expected squared distance of the posterior mean from the unknown. In addition, we included an additive term penalizing for overlapping sensors in the numerical tests, but it is omitted from the following text for brevity.

The approximate posterior covariance $\widetilde{\vGamma}_{\!\rm post}(\vp)$ of \eqref{eq:posterior_cov4}  depends on $\vp$ via~\eqref{eq:FintFbdry}, \eqref{eq:Fpr}, \eqref{eq:tFpr}, and \eqref{eq:posterior_cov4}. In particular, as $\vp$ only affects $\vB_{\rm int}(\vp)$ in~\eqref{eq:FintFbdry}, evaluating the gradient $\nabla_{\! \vp} \mathit{\Phi}_{\!\vA}(\vp)$ is computationally affordable. Hence, given an initial guess $\vp_0$ for the sensor positions, some steepest descent type method can be employed to find a (local) minimum for $\mathit{\Phi}_{\!\vA}$. We do not give the details of our minimization procedure here, since we expect any reasonable algorithm to produce results comparable to those documented in Section~\ref{sec:numerics} below (modulo the speed of convergence). However, it should be emphasized that $\mathit{\Phi}_{\!\vA}$ typically has several local minima, and so finding the global one is almost impossible in practice.

In the rest of this section, we first explain in more detail how the target function \eqref{eq:traceminim} and its gradient can be efficiently evaluated for the above described {\em sliding sensors method} for finding the A-optimal design.   { We consider this method and the proposed computational tools to be novel, although they are inspired by~\cite{Hyvonen2014}}. Afterwards, we briefly recall the {\em $l_0$-sparsification method} \cite{Alexanderian2014, Haber2008} that is employed as the reference approach in one of our numerical tests.

\subsection{Evaluating the target functional and its gradient}
\label{sec:eval_object}
The target functional $\mathit{\Phi}_{\!\vA}(\vp)$ defined in~\eqref{eq:traceminim} can be directly evaluated if the low-rank approximation \eqref{eq:Fapprox} for the boundary measurement has been precomputed. Another option would be to estimate the trace in~\eqref{eq:traceminim} by random matrix techniques~\cite{Alexanderian2014, Haber2008}, but such an approach is not considered in this work. By virtue of~\eqref{eq:posterior_cov4},
\begin{align*}
 \mathit{\Phi}_{\!\vA}(\vp)  = {\rm tr}(\vGamma_{\!\rm pr} \vA^T \!\vA)
    - {\rm tr} \Big( \vL^{-1}\widetilde\vF^{\rm pr}(\vp)^T \big( \widetilde\vF^{\rm pr}(\vp) \widetilde\vF^{\rm pr}(\vp)^T + \widetilde \vGamma_{\!\rm noise}  \big)^{-1} \widetilde\vF^{\rm pr}(\vp) \vL^{-T} \vA^T \! \vA \Big),
\end{align*}
where the first term can be neglected in the optimization since it is independent of $\vp$. As the number of rows in $\widetilde\vF^{\rm pr}(\vp)$ is low(ish), the term $\vC(\vp)^T := \vL^{-1}\widetilde\vF^{\rm pr}(\vp)^T$ and its transpose can be explicitly computed and the inverse of $\vH(\vp) := \widetilde\vF^{\rm pr}(\vp) \widetilde\vF^{\rm pr}(\vp)^T + \widetilde \vGamma_{\!\rm noise}$ can be directly applied.

If one wants to use the norm of $L^2(\Omega)$ to measure the distance between the posterior mean and the unknown, then according to the material in Section~\ref{sec:BayOpt}, the appropriate choice for the finite-dimensional metric is $\vA^T \! \vA = \vM$. In this case, the basic properties of the trace operator  yield
\begin{align}
\mathit{\Phi}_{\vA}(\vp)  = C - {\rm tr} \big( \vC(\vp)^T \vH(\vp)^{-1} \vC(\vp) \vM\big)
  	= C - \sum_{i, j=1}^n \big[ \vM \odot \big( \vC(\vp)^T \vH(\vp)^{-1} \vC(\vp) \big) \big]_{i, j},
  	\label{eq:tracecomp}
\end{align}
where $\odot$ denotes the entrywise matrix product and $C := {\rm tr}(\vGamma_{\!\rm pr} \vM)$ is independent of $\vp$. This expression is inexpensive to evaluate since the mass matrix $\vM$ is sparse and the additive constant $C$ can be ignored.

Although we exclusively resort to the $L^2(\Omega)$-motivated choice $\vA^T \! \vA = \vM$ in the numerical experiments of Section~\ref{sec:numerics}, it is worth mentioning that $\vA = \vL$, i.e.~$\vA^T\!\vA = \vGamma_{\!\rm pr}^{-1}$, leads to a particularly simple form for the A-optimality target functional:
\begin{equation}
\label{eq:tracecompL}
\mathit{\Phi}_{\vL}(\vp)  = n
    - {\rm tr} \Big(\widetilde\vF^{\rm pr}(\vp)^T \big( \widetilde\vF^{\rm pr}(\vp) \widetilde\vF^{\rm pr}(\vp)^T + \widetilde \vGamma_{\!\rm noise}  \big)^{-1} \widetilde\vF^{\rm pr}(\vp)  \Big).
 \end{equation}

\begin{remark}
\label{remark:trace}
{
The infinite-dimensional covariance operator corresponding to the prior~\eqref{eq:prior} is not in the trace class. 
Choosing such a prior and, e.g., $\vA^T \! \vA = \vM$ thus leads to a blowup in the term ${\rm tr}(\vGamma_{\!\rm pr} \vA^T \!\vA)$ when $n$ tends to infinity; see \eqref{eq:tracecompL} for a transparent example of this behavior.
In consequence, $\mathit{\Phi}_{\!\vA}$ is not discretization invariant. However, the $\vp$-dependent part of $\mathit{\Phi}_{\!\vA}(\vp)$ seems to behave well as $n \to \infty$, and thus the same is expected of the corresponding A-optimal sensor positions.}
\end{remark}

Let us then consider evaluating the gradient $\nabla_{\vp}\mathit{\Phi}_{\vA}(\vp)$. To this end, denote the partial derivative with respect to a component of the design parameter vector as $\partial_p$ and observe that by~\eqref{eq:FintFbdry}, \eqref{eq:Fpr} and~\eqref{eq:tFpr},
\begin{equation}
	\partial_p \widetilde{\vF}^{\rm pr}(\vp) = \begin{bmatrix} \partial_p \vB_{\rm int}(\vp) \vPsi \vM \vL^{-1} \\[1mm] 0 \end{bmatrix},
	\label{eq:dFdp}
\end{equation}
where $\partial_p$ only operates on the finite-dimensional observation map $\vB_{\rm int}$ that is a FE discretization of the continuum one defined in~\eqref{eq:pointmeas}. A full description of the corresponding Fr\'echet derivative in the continuum case is given by \eqref{eq:Fderivative} in Appendix~\ref{sec:derivative}. After FE discretization, the derivative $\partial_p \vB_{\rm int}(\vp)$, i.e.~\eqref{eq:Fderivative} with $q$ chosen as the appropriate Cartesian coordinate vector, has a particularly simple form when triangular/tetrahedral linear elements are used: The $i$th row in $\vB_{\rm int}(\vp)$ can be considered to be a weighted sum of triangular/tetrahedral barycentric coordinates at the $i$th measurement location, and the partial derivatives of the barycentric coordinates with respect to $p$ are trivial to evaluate.

After $\partial_p \vB_{\rm int}(\vp)$ has been computed, $\vM \! \vPsi^T  \! \partial_p \vB_{\rm int}(\vp)^T$ can be evaluated by mimicking \eqref{eq:Ftransp}, and subsequently $\partial_p \widetilde{\vF}^{\rm pr}(\vp)$ can be formed as in \eqref{eq:dFdp}. The linearity of the trace yields
\begin{align*}
    \partial_p \mathit{\Phi}_{\!\vA}(\vp) &= {\rm tr}\big(\partial_p \widetilde{\vGamma}_{\rm post}(\vp) \vA^T \! \vA \big) \\[1mm]
       &= - {\rm tr} \Big( \vL^{-1} \partial_p\big(\widetilde\vF^{\rm pr}(\vp)^T (\widetilde\vF^{\rm pr}(\vp) \widetilde\vF^{\rm pr}(\vp)^T + \widetilde \vGamma_{\!\rm noise} )^{-1} \widetilde\vF^{\rm pr}(\vp)\big) \vL^{-T} \vA^T \! \vA \Big),
\end{align*}
where the remaining derivative can be calculated using the product rule and the differentiation formula for an inverse matrix. Finally, for the choices $\vA^T \! \vA = \vM$ and $\vA^T \! \vA = \vGamma_{\!\rm pr}^{-1}$ the whole derivative can be evaluated following the ideas leading to~\eqref{eq:tracecomp} and \eqref{eq:tracecompL}, respectively.

\subsection{Sparsification method}
\label{sec:sparsification}
Another approach to finding an A-optimal design is the so-called $l_0$-sparsification method that is introduced  more thoroughly in~\cite{Alexanderian2014, Haber2008}. In short, there is initially a large number of sensor candidates whose locations do not change during the minimization process. In particular, $\vF_{\rm int}$ is constructed either explicitly or approximately during the initialization phase, and it remains unaltered for the rest of the algorithm. A weight $w_i \in (0, 1)$ is assigned to each sensor candidate, and an A-optimal design is sought by driving the weights towards binary $w_i \in \{ 0, 1 \}$ by a certain iterative procedure explained in \cite{Alexanderian2014, Haber2008}.

Because the $l_0$-sparsification method involves a sequence of convex minimization problems with box constraints $w_i \in [0, 1]$, as a numerical optimization problem it is arguably simpler than our sliding sensors method. On the negative side, the number of active sensor locations remaining after the procedure is unpredictable since it depends implicitly on the input parameters of the algorithm.


\section{Numerical Tests}
\label{sec:numerics}
\label{sec:numresults}
We demonstrate the introduced computational framework with numerical examples. The forward problem~\eqref{eq:fwdPDE} is spatially discretized by piecewise linear FEs using triangular or tetrahedral meshes depending on the spatial dimension. The implicit midpoint rule is used for time integration. We generate simulated measurements by solving the discretized PDE with a chosen {\em true} source $f_{\mathrm{true}}$ and adding independent realizations of zero-mean Gaussian random noise to the temperatures both at the internal and the boundary sensors. To be more precise, the noise covariance is of the form
\begin{equation}
	\vGamma_{\!\mathrm{noise}} = \begin{bmatrix} \gamma_{\mathrm{int}}^2 \vI & \bm{0} \\ \bm{0} & \gamma_{\mathrm{bdry}}^2 \vI \end{bmatrix} \in \R^{m \times m},
	\label{eq:noiseparam}
\end{equation}
where the diagonal blocks are of the sizes $m_{\rm int}m_t \times m_{\rm int}m_t$ and $m_{\rm bdry}m_t \times m_{\rm bdry}m_t$, respectively, and the  parameters $\gamma_{\mathrm{int}}, \gamma_{\mathrm{bdry}} > 0$ are the standard deviations of the measurement noise at the internal and boundary sensors, respectively. By stating that the amount of measurement noise is $p$\%, we mean that the size of $\gamma_{\mathrm{int}}$ (resp.~$\gamma_{\mathrm{bdry}}$) is $0.01p$ times the maximal simulated temperature at the internal (resp.~boundary) sensors over the measurement time interval. In all numerical tests we assume to know $\vGamma_{\!\mathrm{noise}}$ and use it in the formulas for the posterior mean and covariance, which is somewhat unrealistic but simplifies the considerations. The prior mean is always zero, whereas the prior covariance is of the form~\eqref{eq:prior}, where the free parameters $\alpha$ and $\beta$ are chosen separately for each numerical test.

The {\em reconstruction} is defined to be
$$
\widehat{f} = \sum_{i=1}^n \widehat{\vx}_i \phi_i,
$$
where $\widehat{\vx} = \widehat{\vx}(\vp)$ is the posterior mean of the nodal values, which depends on the design parameter $\vp \in \R^N$, and $\{\phi \}_{i=1}^n$ is the employed FE basis. The reconstruction quality is quantified by the relative $L^2(\Omega)$-error
\begin{equation}
	\textrm{err}_{\textrm{rel}} := \frac{ \big\| \widehat f - f_{\mathrm{true}} \big\|_{L^2(\Omega)}}{ \| f_{\mathrm{true}} \|_{L^2(\Omega)}},
	\label{eq:L2error}
\end{equation}
which naturally also depends on $\vp$ via $\widehat{f}$.

We consider three numerical experiments. The first one involves a simple rectangular geometry and demonstrates the positive effect of including a time transient measurement in addition to the steady state data. We also numerically demonstrate that the A-optimality is a reasonable quantifier for the quality of the  measurement setup and briefly compare the optimal internal sensor designs produced by the sliding sensors and sparsification methods.

The second experiment considers a more realistic two-dimensional geometry, with physically reasonable parameter values and a true source resembling an iron loss field obtained from an existing iron loss model. The performance of the sliding sensors method is tested by comparing the accuracy of reconstructions obtained using nonoptimized and A-optimized sensor locations.

Both aforementioned two-dimensional problems are computationally {\em small}, that is, all necessary matrices and associated linear operations can be explicitly stored and computed in a reasonable time without a need for the model order reduction methods described in Sections~\ref{sec:lowrank} and~\ref{sec:eval_object}. In practice, by ``small'' we mean fewer than $n=5000$ FE degrees of freedom and fewer than $m=5000$ measurement data. However, our third and final test is supposed to be {\em non-small}. It considers a semirealistic three-dimensional setup, with $n \approx 5\cdot 10^4$. The model order reduction method for the boundary measurements described in Section~\ref{sec:lowrank} is tested and used in connection with the sliding sensors algorithm.

A single desktop computer and the MATLAB software were used for all computations.  The A-optimization by the sliding sensors method was performed by a custom steepest descent type algorithm, whereas a black box solver (L-BFGS-B~\cite{Morales2011BFGS}) was utilized with the sparsification method. Although not tested in this paper, we believe the proposed methods to be easily parallelizable so that they could also be used for larger-scale problems; see Section~\ref{sec:discussion} for further discussion on this matter.

\subsection{Experiment I:~the unit square}
\label{sec:2dtestsquare}
\begin{figure}
	\includegraphics[width = 0.49\columnwidth]{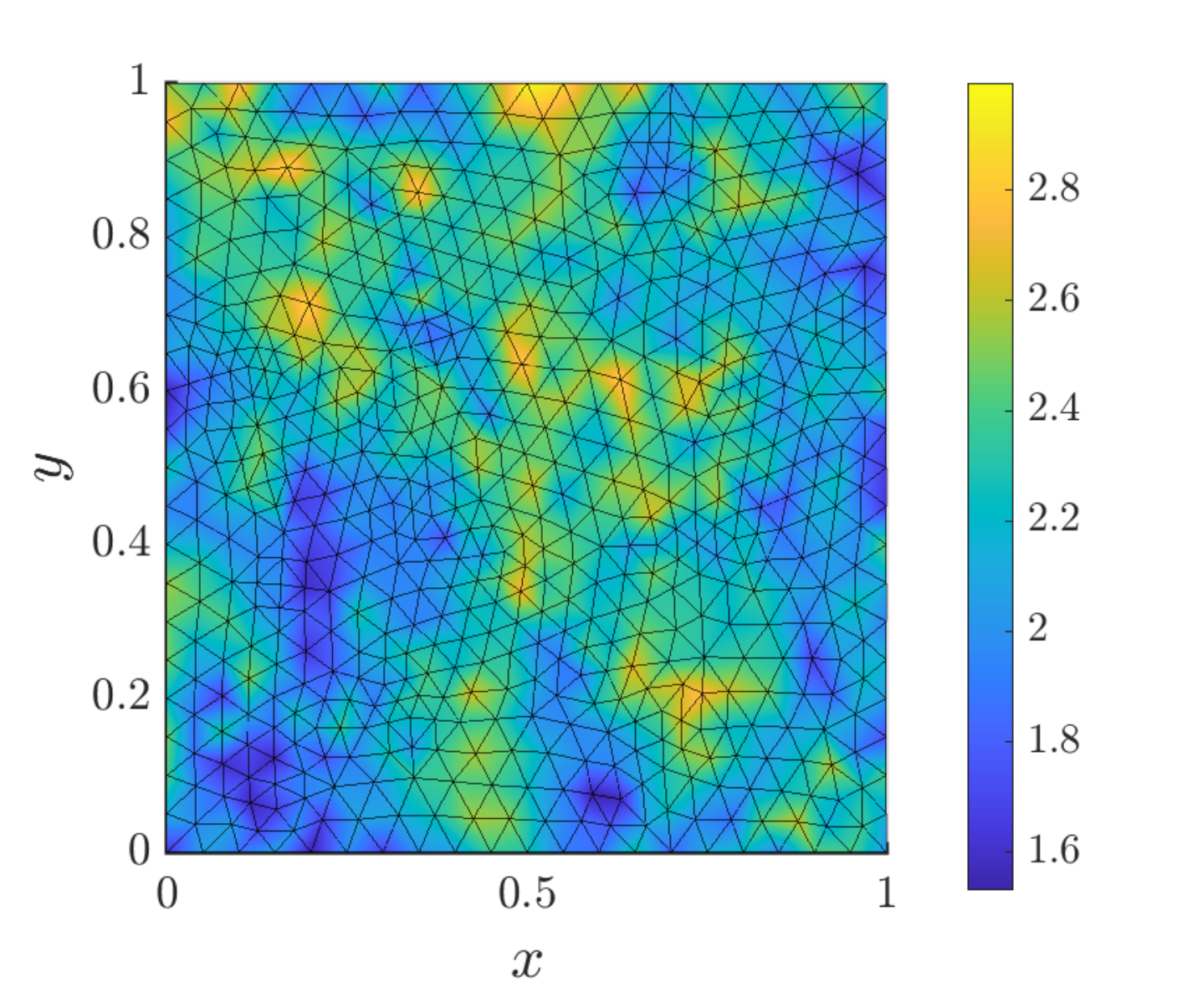}		
	\includegraphics[width = 0.49\columnwidth]{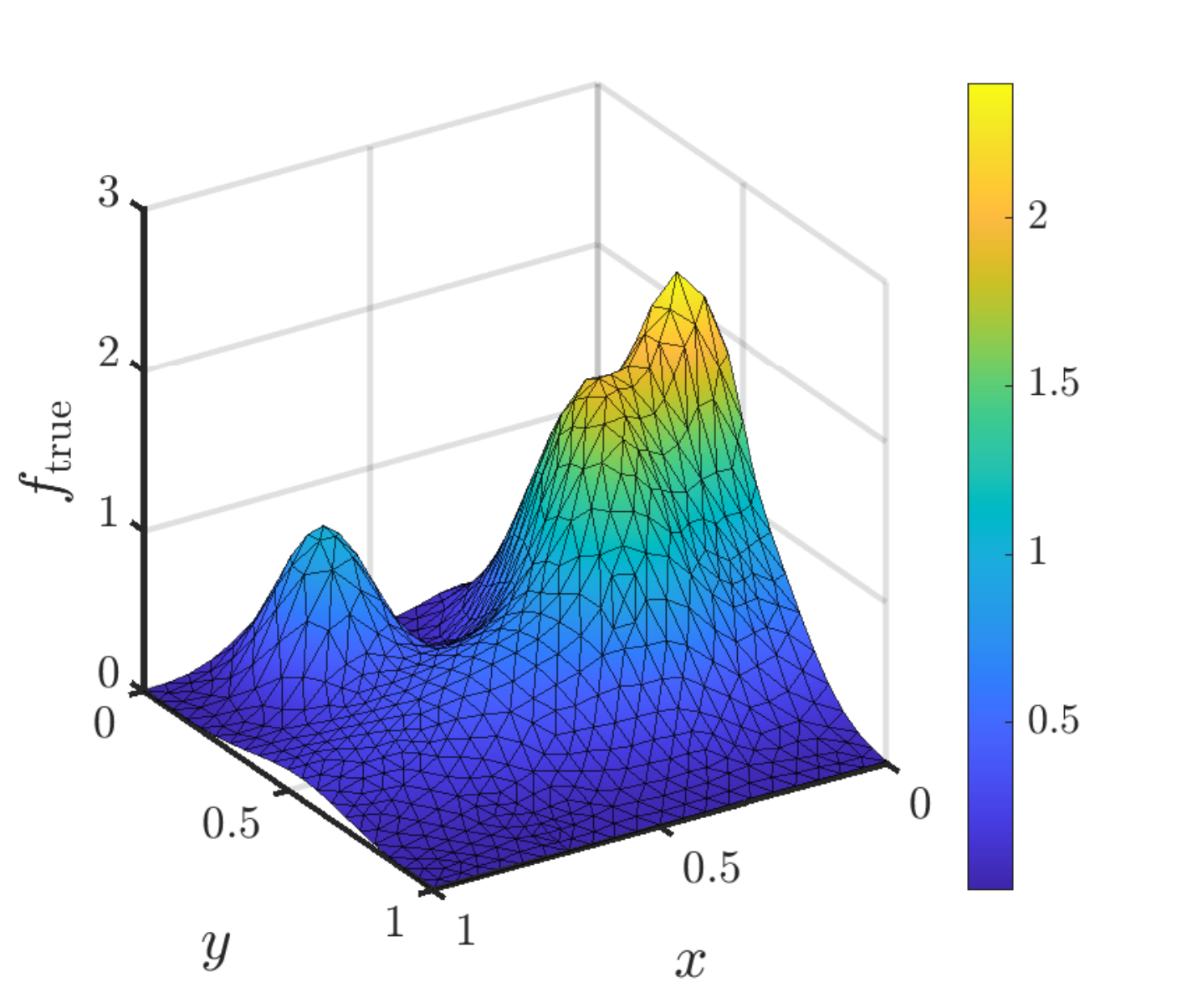}		
	\caption{\small Experiment I. {\sc Left}: Random draw from the assumed prior distribution for the source. {\sc Right}: Random realization of the true source \eqref{eq:frandom}.}
		\label{fig:random}
\end{figure}

Let $\Omega = (0, 1)^2 \subset \R^2$ be the unit square characterized by the  homogeneous parameters $\rho, \kappa, h \equiv 1$ and $\Gamma_{\rm R} = \partial \Omega$ in \eqref{eq:fwdPDE}. We consider $200$ `multi-modal Gaussian' realizations for the target source,
\begin{equation}
	f_{\mathrm{true}}(x,y) = \sum_{i=1}^{M} \exp\!\big(-a_i \big(x - \bar{x}^{(i)}\big)^2 - b_i \big(y - \bar{y}^{(i)}\big)^2\big), \qquad (x,y) \in \Omega,
	\label{eq:frandom}
\end{equation}
where the number of modes $M$ is drawn from the uniform distribution over $\{1, \dots 11\}$, the parameters $a_i, b_i\in \R_+$ are independently drawn from the uniform distribution over the interval $[0,100]$ and the center points $(\bar{x}^{(i)},\bar{y}^{(i)})$ are drawn from the uniform distribution over $\Omega$. The temperature is measured at the time instants $t_j = j T / m_t$, $j=1, \dots,m_t$, where $T=0.6$ and $m_t = 20$. The measurements are contaminated by 0.5\% of additive noise, and the free parameters in the prior covariance \eqref{eq:prior} are set to $\alpha = 10$ and $\beta = 0.1$. The left-hand image in Figure~\ref{fig:random} shows a random draw from the prior distribution of the source, and the right-hand image depicts a single realization of the random target source \eqref{eq:frandom}. In particular, it is obvious that \eqref{eq:frandom} is not completely inline with the assumed prior model for the unknown source.

\begin{figure}
	\includegraphics[width = 0.49\columnwidth]{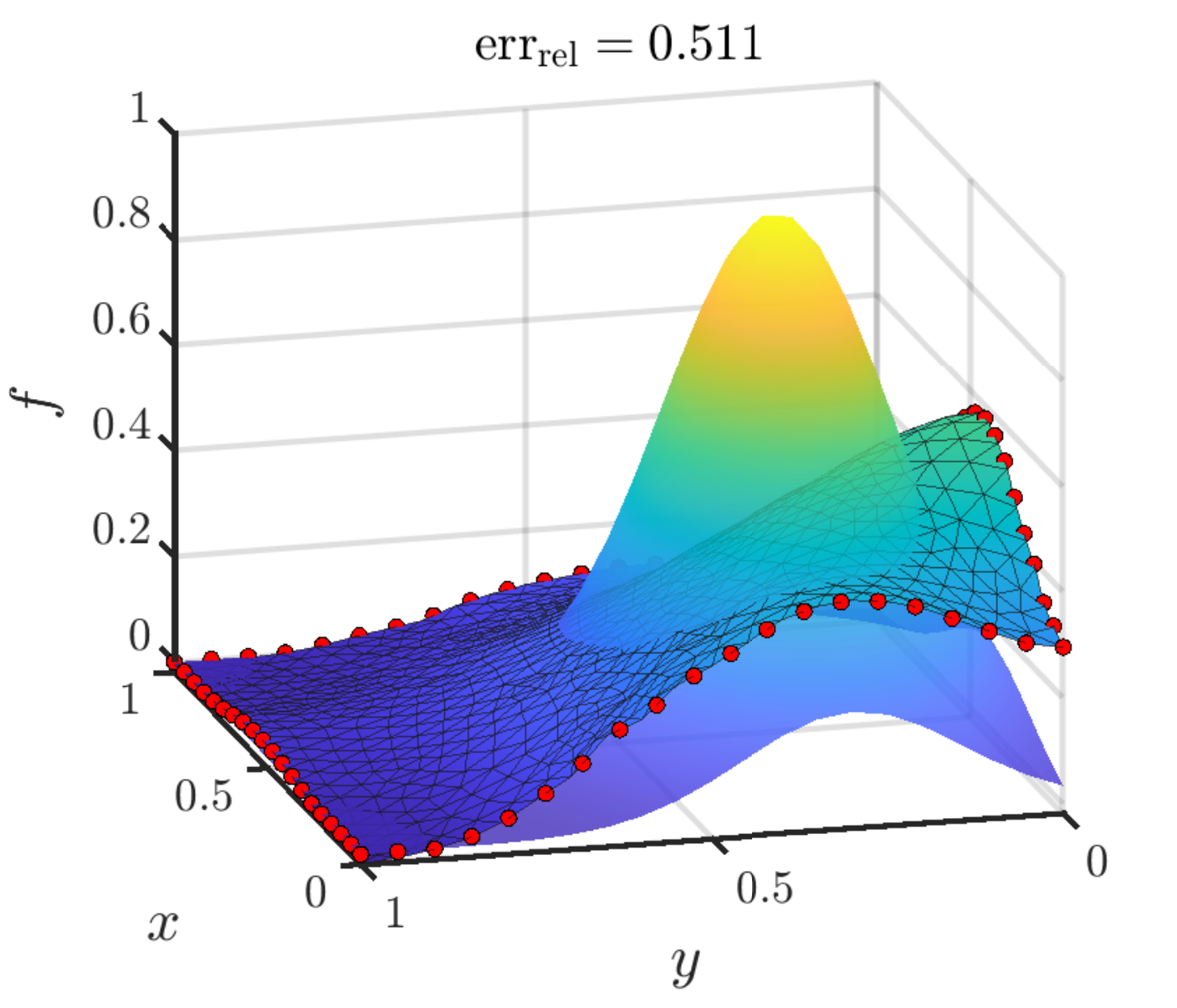}
	\includegraphics[width = 0.49\columnwidth]{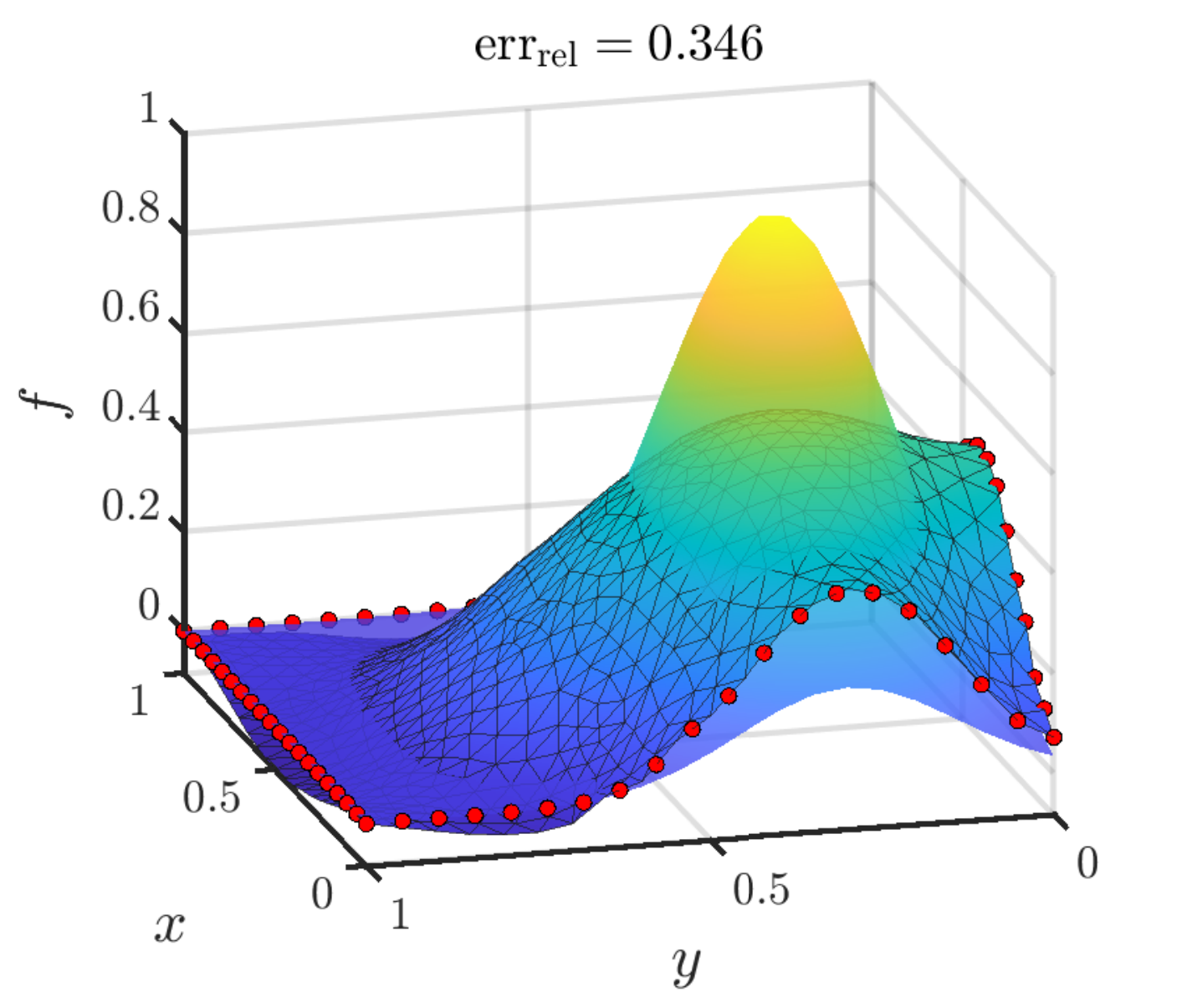}		
	\caption{\small Experiment I. Gaussian source with no internal measurements. The gridded surface represents the reconstruction $\widehat f$ and the gridless surface the true source $f_{\mathrm{true}}$. {\sc Left}: Steady state measurement. {\sc Right}: Time transient measurement.}
		\label{fig:steadyTransient}
\end{figure}

We start the actual tests by demonstrating how time transient data yields additional information about the unknown source. The temperature is measured at $m_{\rm bdry} = m = 76$ equidistant positions on the boundary $\partial \Omega$, while no internal sensors are utilized. Figure~\ref{fig:steadyTransient} shows two reconstructions and compares them with the corresponding true source that is of the form~\eqref{eq:frandom} with only a single mode: the left-hand reconstruction is based on steady state data, whereas the one on the right exploits transient data. It is obvious the steady state measurements carry much less information on the behavior of the target source in the interior of the domain, albeit mere boundary measurements do not lead to a good reconstruction in either case even though the considered source is rather simple.

\begin{figure}
	\includegraphics[width = 0.49\columnwidth]{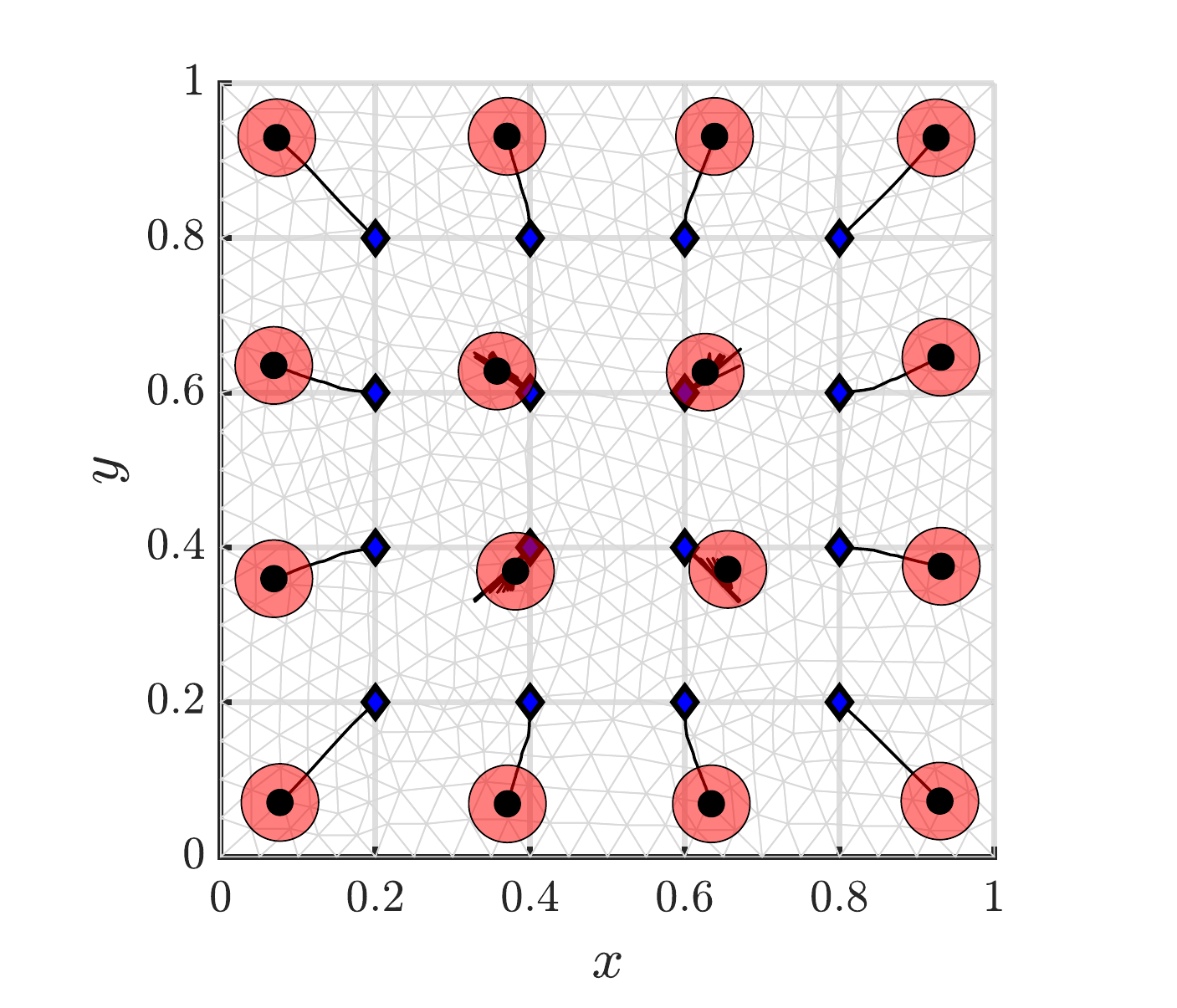}
	\includegraphics[width = 0.49\columnwidth]{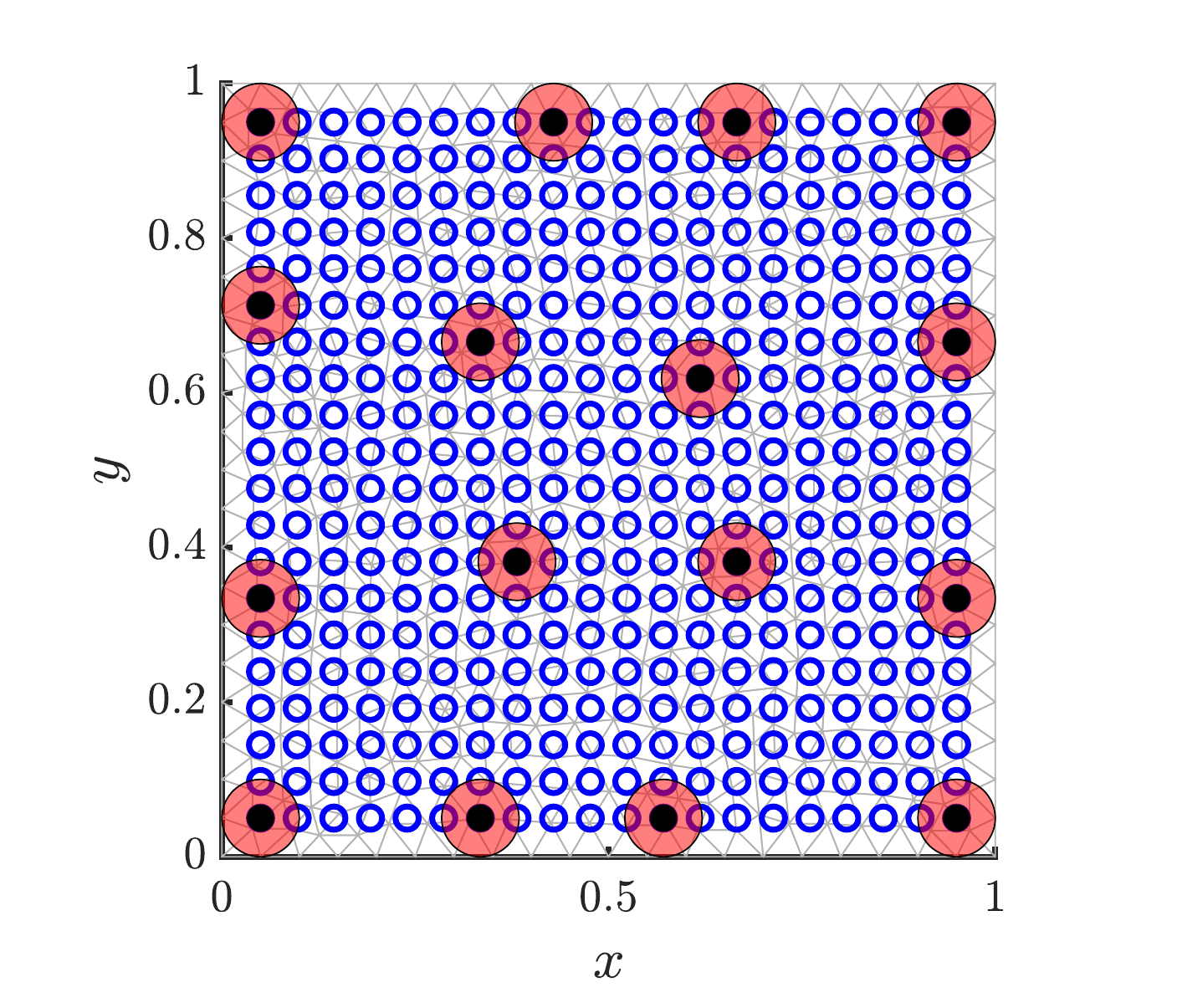}		
	\caption{\small Experiment I. Approximate A-optimal configurations for $m_{\rm int}= 16$ internal sensors. {\sc Left}: Sliding sensors. Initial locations (diamond), A-optimized locations (circle), sensor size (transparent). {\sc Right}: Sparsification method. Candidate locations (non-filled), accepted locations (filled), sensor size (transparent).}
		\label{fig:sliding_vs_sparsi}
\end{figure}

Let us next compare the A-optimal sensor positions produced by the sliding sensors method with those predicted by the sparsification algorithm. As always in our numerical tests, the A-optimality target functional is $\mathit{\Phi}_\vA$ from \eqref{eq:traceminim} with $\vA^T\!\vA = \vM$. We consider a setting with $m_{\rm int} = m = 16$ internal sensors and no boundary measurements. The sensors measure spatial averages over disks of radius 0.05 (cf.~\eqref{eq:pointmeas}). The left-hand image of Figure~\ref{fig:sliding_vs_sparsi} shows the approximate A-optimal positions for the internal sensors given by the sliding sensors method if the iteration is started from the regular grid $\{1/5,2/5,3/5,4/5\}^2$; the outermost sensors nearly touch the boundary of $\Omega$ at the end of the minimization process. For comparison, the right-hand image in Figure~\ref{fig:sliding_vs_sparsi} presents the A-optimal locations for $16$ sensors produced by the sparsification method if the process is started with a regular grid of $400$ possible measurement positions; observe that controlling the precise number of active sensors is not trivial in the sparsification method as it can only be achieved via trial and error by manually tuning the parameters of the algorithm (cf.~\cite{Alexanderian2014}). Up to the resolution of the candidate grid of the sparsification method, the two sets of A-optimal positions are in good agreement. However, the value of the target functional associated to the output of the sparsification method would still slightly decrease if the corresponding sensor positions were used as the initial guess for the sliding sensors method.

\begin{figure}
	\includegraphics[width = 0.49\columnwidth]{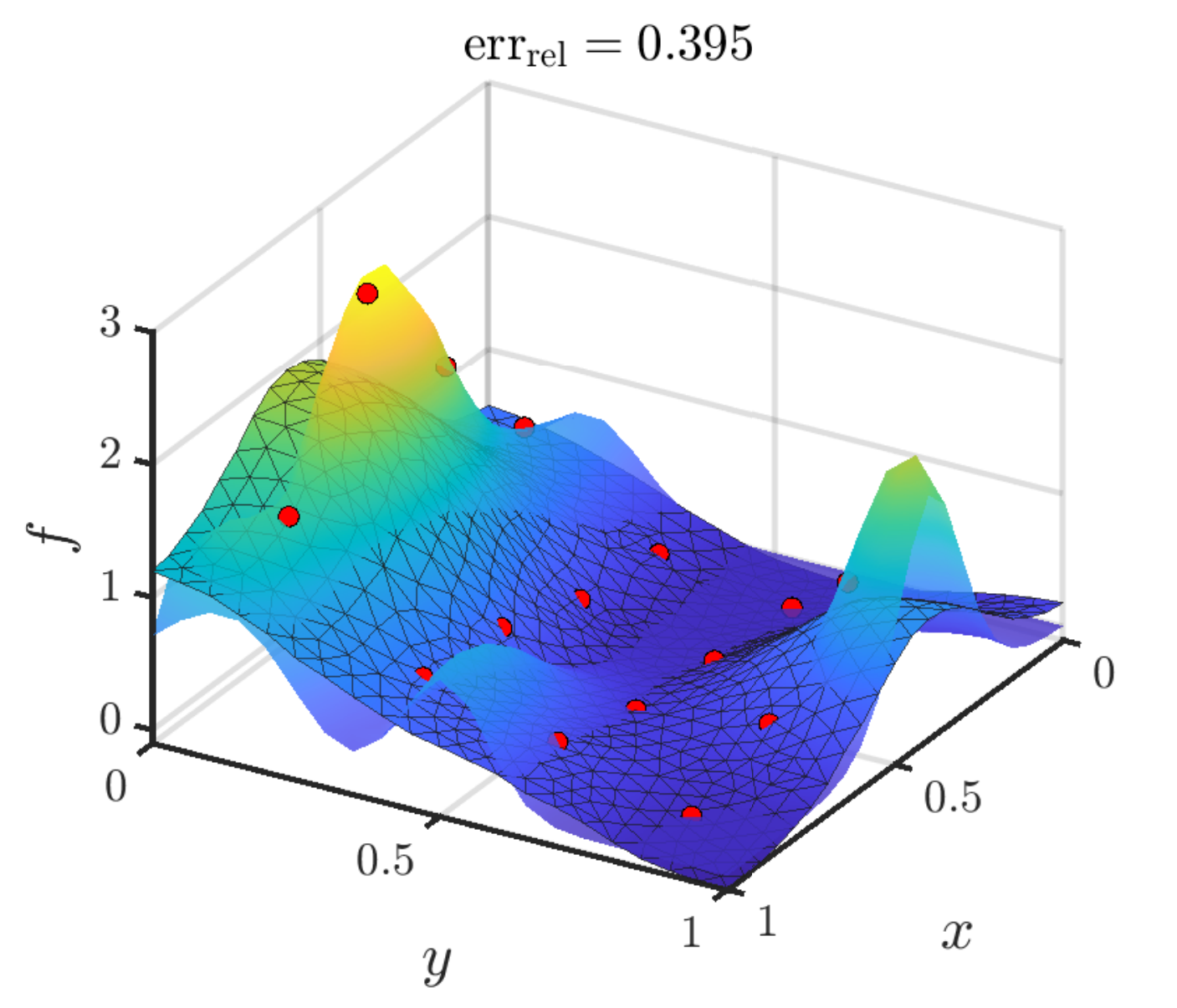}
	\includegraphics[width = 0.49\columnwidth]{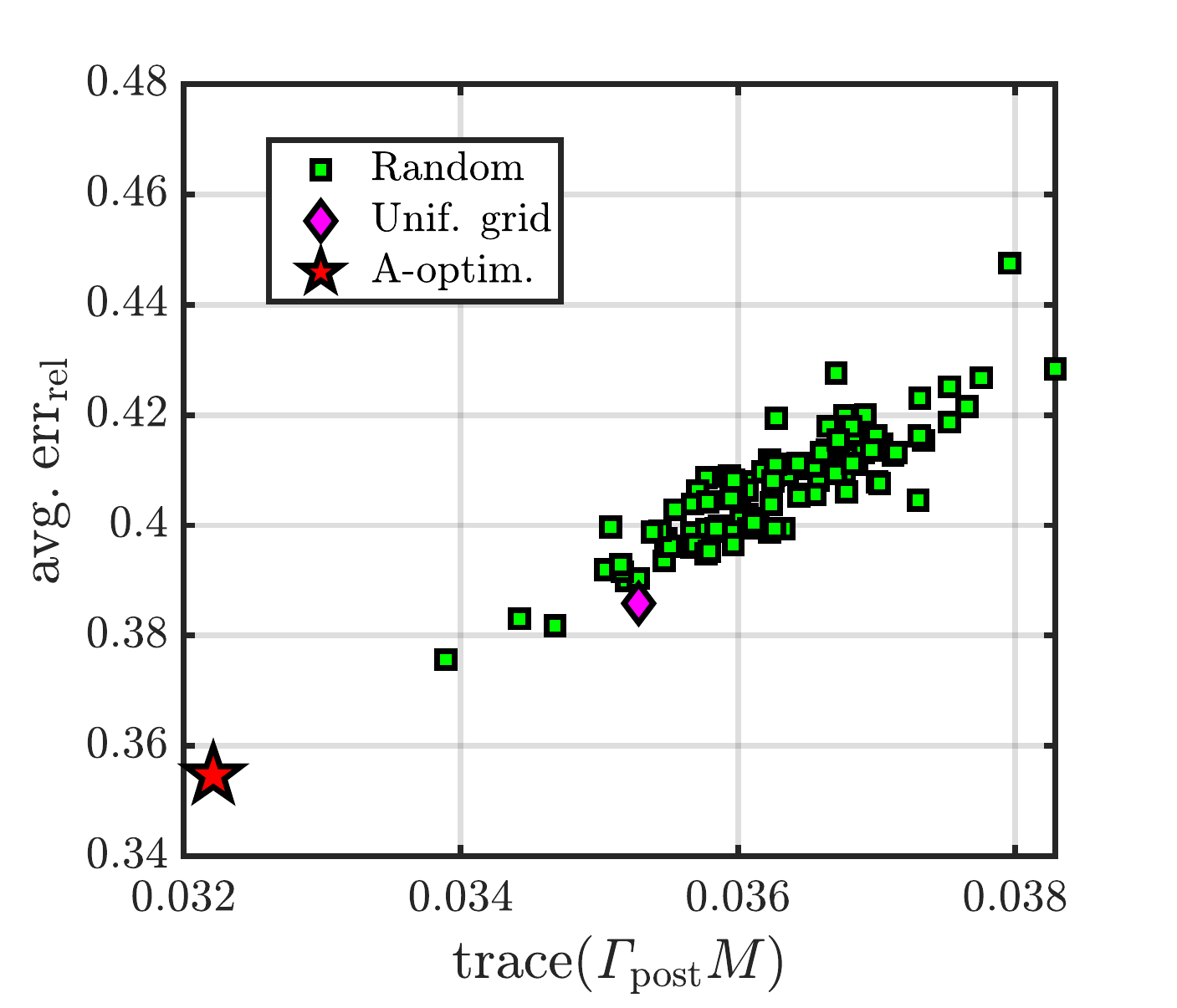}		
	\caption{\small Experiment I. {\sc Left}: An example of a randomly generated source of the form \eqref{eq:frandom} and random internal sensor positions. The gridded surface represents the reconstruction $\widehat f$ and the gridless surface the true source $f_{\mathrm{true}}$. {\sc Right}: A-optimality target $\mathit{\Phi}_{\vA} = \trace(\vGamma_{\!\rm post} \vM)$ versus the average relative $L^2(\Omega)$ reconstruction error $\textrm{err}_{\textrm{rel}}$ over 200 randomly generated sources of the form \eqref{eq:frandom} for 100 random sensor configurations (squares), the initial  uniform grid of sensors (diamond), and the A-optimized sensor positions (star).
	}
		\label{fig:trace_vs_error}
\end{figure}

To complete the first numerical experiment, we illustrate the connection between the value of the A-optimality target functional, i.e.~$\mathit{\Phi}_\vA$ from \eqref{eq:traceminim} with $\vA^T\!\vA = \vM$, and the reconstruction quality measured by \eqref{eq:L2error}. We introduce $100$ random configurations for the $16$ internal sensors by adding to the regular grid points $\{1/5,2/5,3/5,4/5\}^2$ independent realizations of a random vector distributed uniformly over $[-1/5,1/5]^2$. The reconstruction $\widehat f$ is then computed from noisy measurements corresponding to all $200$ randomly generated target sources of the form~\eqref{eq:frandom} for all 100 sensor configuration, as well as for the A-optimal one obtained by the sliding sensors method. The relative error~\eqref{eq:L2error} is computed for all $101\cdot 200 = 20\,200$ cases, and the mean of the relative errors for each sensor configuration is compared to the corresponding value of  $\mathit{\Phi}_{\!\vA} ={\rm tr}(\vGamma_{\!\rm post}\vM)$ in the right-hand image of Figure~\ref{fig:trace_vs_error}. Since the correlation between the A-optimality measure and the relative reconstruction error is apparent, $\mathit{\Phi}_{\!\vA}$ appears to be a plausible measure for the reconstruction quality --- even though the target sources are not drawn from the assumed prior distribution.

\subsection{Experiment II:~a semirealistic configuration in two dimensions}
\label{sec:2dtest}
\begin{figure}
   	\includegraphics[width = 0.350\columnwidth]{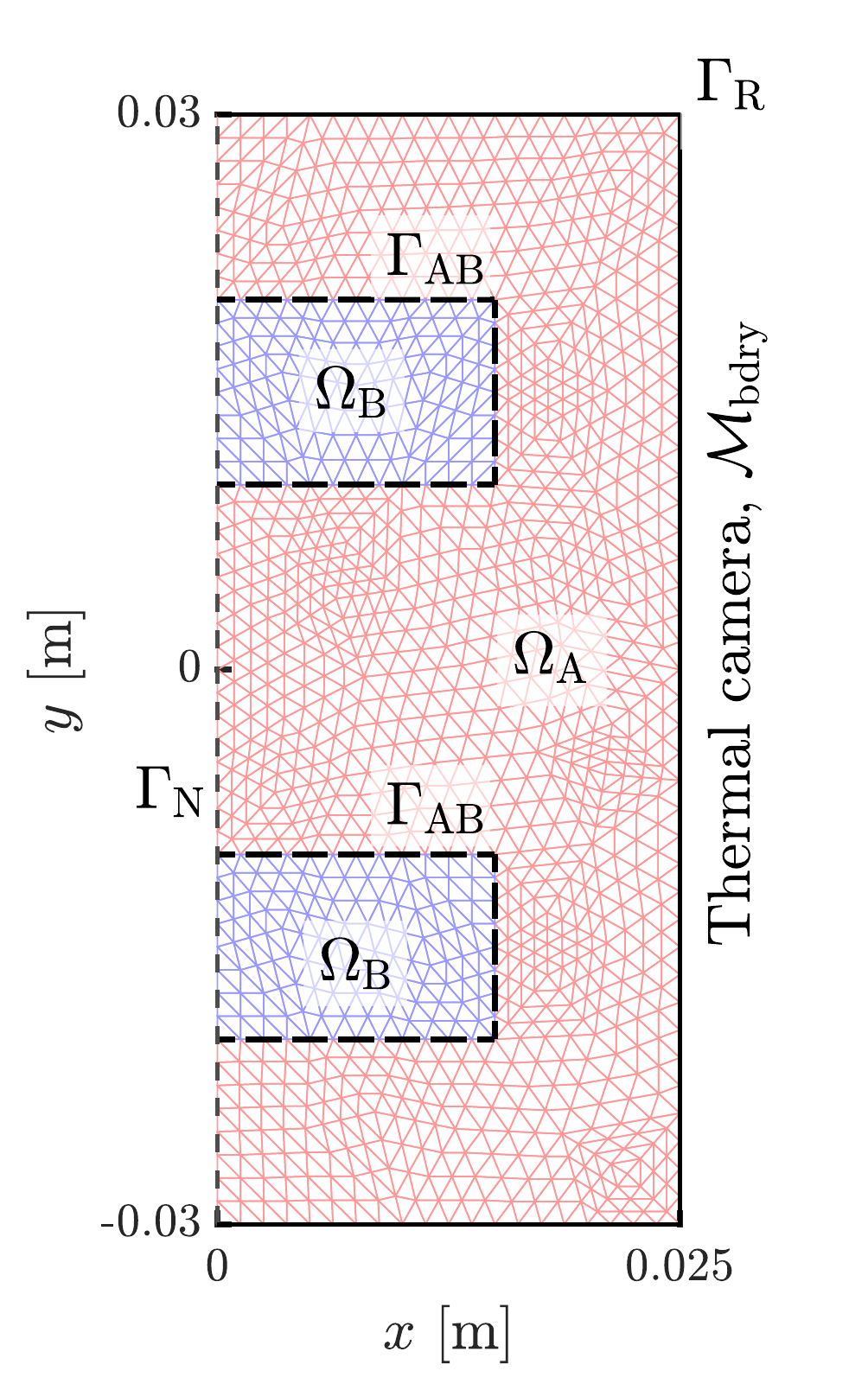}
	\includegraphics[width = 0.640\columnwidth]{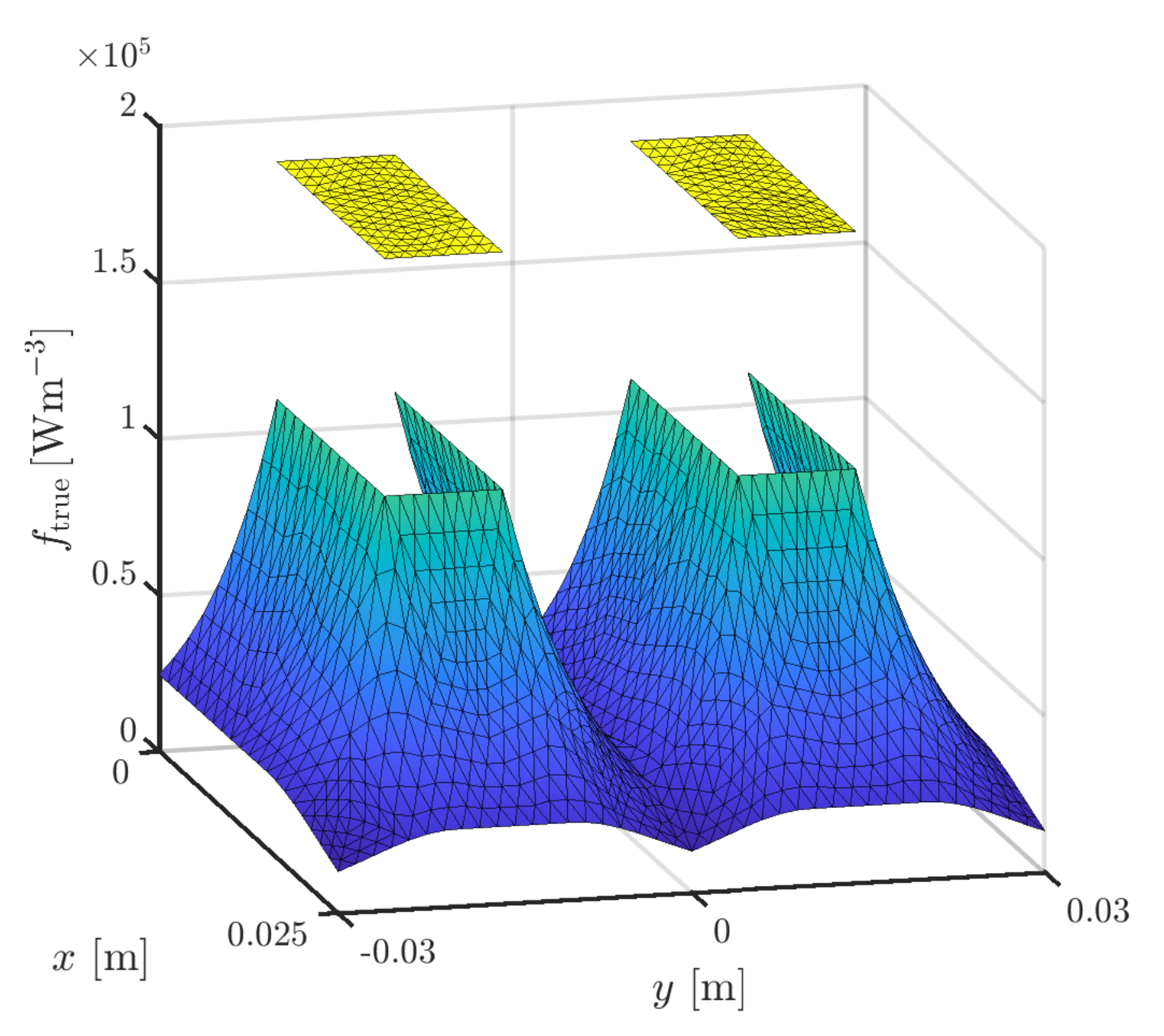}
    \caption{\small Experiment II. {\sc Left}: Half of a cross section of a transformer showing the  the iron core $\Omega_{\mathrm{A}}$ and the coils $\Omega_{\mathrm{B}}$, the boundaries $\Gamma_{\mathrm{N}}$, $\Gamma_{\mathrm{R}}$, $\Gamma_{\mathrm{AB}}$, and the measurement boundary $\mathcal M_{\mathrm{bdry}}$. {\sc Right}: The true heat source $f_{\mathrm{true}}$ representing the coil and iron losses. The component in $\Omega_{\mathrm{B}}$ is scaled down for a better visualization.}
    \label{fig:2dgeom}
\end{figure}

Our second example involves a more realistic geometry, parameters and measurement setup, with all quantities given in the appropriate SI units. The domain $\Omega = (0, 0.025) \times (-0.030, 0.030) \subset \R^2$, shown in the left-hand image of Figure~\ref{fig:2dgeom}, represents one half of a cross section of a small transformer. The domain is decomposed as $\Omega = \Omega_{\mathrm{A}} \cup \Gamma_{\!\mathrm{AB}} \cup \Omega_{\mathrm{B}}$, where $\Omega_{\mathrm A}$ is the iron core, $\Omega_{\mathrm B}$ is the coil and $\Gamma_{\mathrm{AB}} := (\partial \Omega_{\mathrm A} \cup \partial \Omega_{\mathrm B}) \setminus \partial \Omega$ corresponds to a thin layer filled with insulating material. The material parameters in the parabolic PDE~\eqref{eq:heateq} are assumed to be piecewise constant, that is, the pairs $(\kappa_{\mathrm A}, \rho_{\mathrm A}) \in \R_+^2$ and $(\kappa_{\mathrm B}, \rho_{\mathrm B}) \in \R_+^2$ characterize the properties of $\Omega_{\mathrm A}$ and $\Omega_{\mathrm B}$, respectively. Since $\Omega$ models one half of a transformer, the symmetry boundary is $\Gamma_{\mathrm{N}} = \{ (x, y) \in \partial \Omega \  | \  x=0 \}$, and the remainder $\Gamma_{\mathrm{R}} = \partial \Omega \setminus \overline{\Gamma}_{\mathrm{N}}$ is the Robin boundary that models the heat conduction into the surrounding air.

The heat conduction in the insulating layer $\Gamma_{\!\mathrm {AB}}$  is modeled by the (weak forms of the) jump conditions
\begin{subequations}
\label{eq:jump}
\begin{align}
    \kappa_{\rm A} \nu \cdot \nabla u_{\mathrm A} &= \frac{\kappa_{\mathrm{ins}}}{d_{\mathrm{ins}}} (u_{\mathrm B} - u_{\mathrm A}) \qquad {\rm on} \  \Gamma_{\mathrm{AB}},  \\[1mm]
    \kappa_{\rm B} \nu \cdot \nabla u_{\mathrm B} &= \frac{\kappa_{\mathrm{ins}}}{d_{\mathrm{ins}}} (u_{\mathrm B} - u_{\mathrm A}) \qquad {\rm on} \ \Gamma_{\mathrm{AB}},
\end{align}
\end{subequations}
where $\nu$ is the unit normal of $\Gamma_{\mathrm{AB}}$ pointing toward $\Omega_{\rm B}$, $u_{\mathrm A}$ and $u_{\mathrm B}$ are the temperatures in the subdomains $\Omega_{\mathrm A}$ and $\Omega_{\mathrm B}$, respectively, and $\kappa_{\mathrm{ins}}$ and $d_{\mathrm{ins}}$ are the constant heat conductivity and thickness of the insulating layer, respectively~\cite{Javili2012}. In practical terms, the heat conduction in $\Omega_{\mathrm A}$ and $\Omega_{\mathrm B}$ is modeled by two separate FE systems that are coupled via the above described boundary condition on $\Gamma_{\mathrm{AB}}$. By using this approximation, one avoids constructing an unnecessarily fine mesh in the thin layer, as it is well known that thin elements can be problematic for FE solvers. The parameters for~\eqref{eq:fwdPDE} and~\eqref{eq:jump} are $\kappa_{\mathrm A} = 10$, $\rho_{\mathrm A} = 3.43\cdot 10^6$, $\kappa_{\mathrm B} = 26$, $\rho_{\mathrm B} = 3.26\cdot 10^6$, $h = 14$, $d_{\mathrm{ins}} = 5 \cdot 10^{-4}$ and $\kappa_{\mathrm{ins}} = 0.028$.

The true heat source
\begin{equation}
	f_{\mathrm{true}}(x) =
	\begin{cases}
		2.557 \cdot 10^5, & \text{if } x \in \Omega_{\mathrm{B}}, \\[1mm]
                10^5 \exp \!\big( -150 \operatorname{dist}(x, \Omega_{\mathrm{B}}) \big), & \text{if } x \in \Omega_{\mathrm{A}},
	\end{cases}		
	\label{eq:ftrue2}
\end{equation}
is shown in Figure~\ref{fig:2dgeom}. Here $\operatorname{dist}(x, \Omega_{\mathrm{B}})$ denotes the Euclidean distance from $x$ to the coils $\Omega_{\rm B}$. The source in $\Omega_{\mathrm B}$ models the ohmic loss generated by the coil current, totaling 2.0W if uniform distribution in the $z$-direction over 0.025m is assumed. On the other hand, the source in $\Omega_{\mathrm A}$ represents the iron loss, totaling 1.5W, and this latter source component is the object of primary interest in the reconstruction process. In practice, the loss in the coil can be reliably estimated by measuring the current and resistance of the coil, and this information could also be included in our model by choosing a nonzero mean for the prior distribution of the heat source. However, such an approach is not considered in this work.

\begin{figure}
   \includegraphics[width = 0.49\columnwidth]{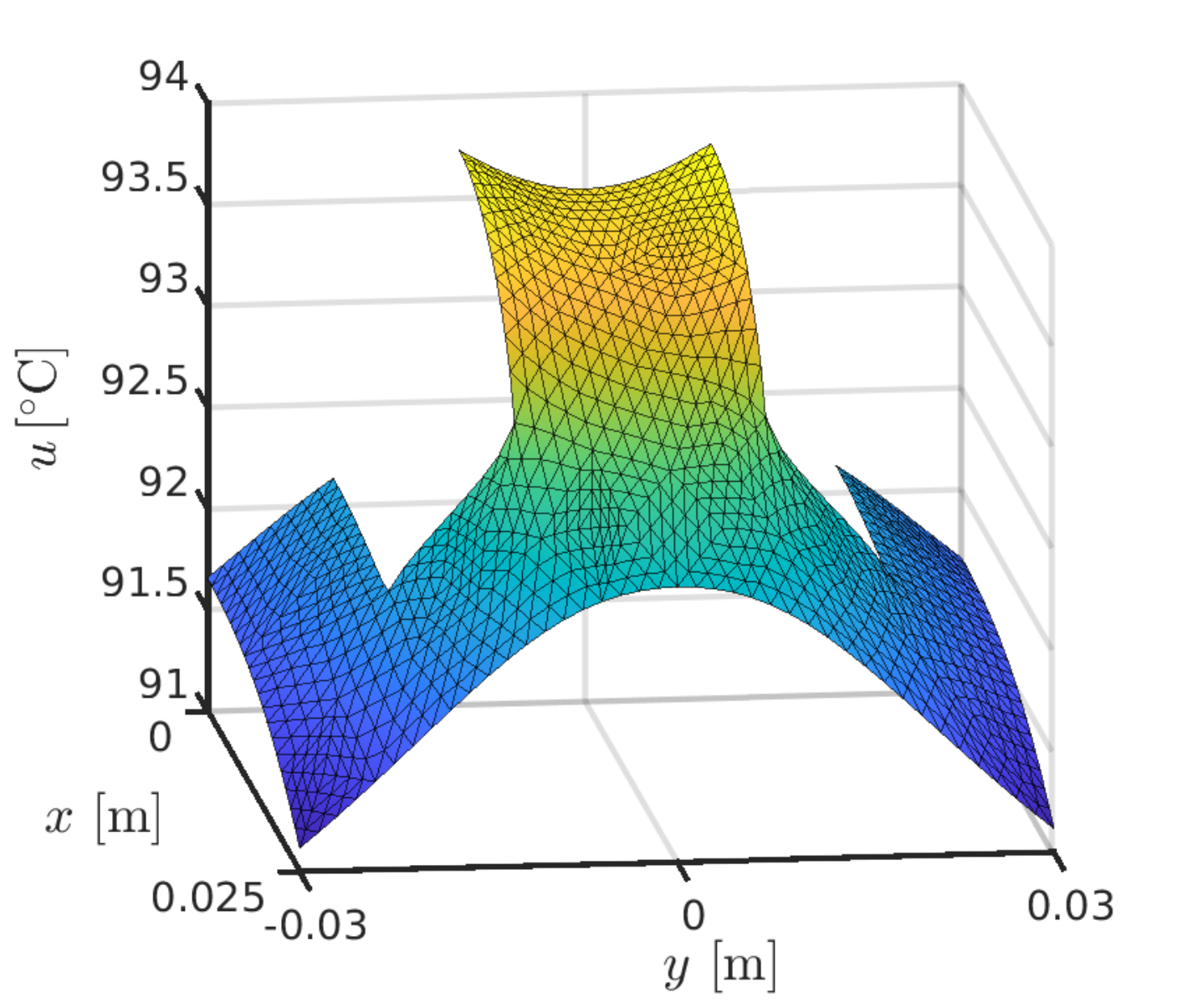}
   \includegraphics[width = 0.49\columnwidth]{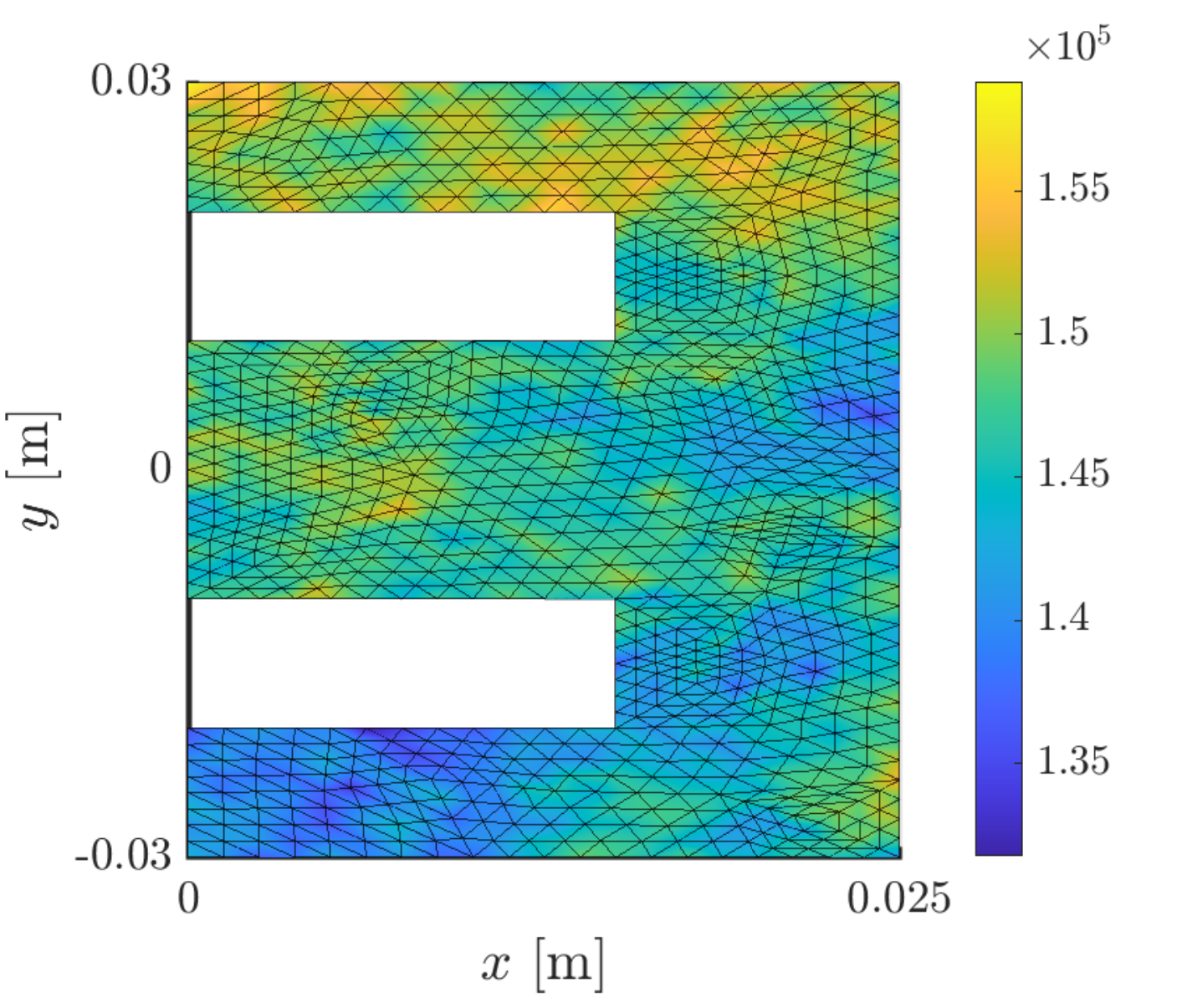}
   \caption{\small Experiment II. {\sc Left}: The iron core temperature at the steady state with the coil temperature of approximately 110\textdegree{}C omitted. The temperature is unrealistically high because the heat dissipation in the $z$-direction is neglected in the two-dimensional model. {\sc Right}: A sample from the prior distribution.}
    \label{fig:2dtemp}
\end{figure}

The temperature on $\mathcal{M}_{\rm bdry} = \{ (x, y) \in \partial \Omega \  | \  x=0.025 \}$ is measured by a thermal camera, whereas only a limited number of small circular sensors are inserted inside the object. In all tests, there are $m_{\rm bdry} = 60$ boundary sensors, i.e.~pixels in the thermal images, but the number of internal sensors $m_{\rm int}$ may vary. The measurement times are  $t_j = j T / m_t$, $j=1, \dots,m_t$, where $T=2 \cdot 10^4$, $m_t = 20$. The measurement noise level is assumed to be $0.1\%$. The measured temperatures inside the iron core are in the range $[0, 94]$\textdegree C, with the steady state temperature shown on the left in Figure~\ref{fig:2dtemp}. The coil temperature is omitted from the image as it is considerably higher, approximately 110\textdegree C. It should be noted that the simulated steady state temperature is actually unrealistically high because the heat dissipation in the $z$-direction is neglected in this two-dimensional model. We consider a zero-mean Gaussian prior distribution with constant parameters $\alpha = 10^{-8}$ and $\beta = 10^{-7}$ in~\eqref{eq:prior}. The right-hand image in Figure~\ref{fig:2dtemp} shows a random draw from the prior.

\begin{figure}
   \includegraphics[width = 0.49\columnwidth]{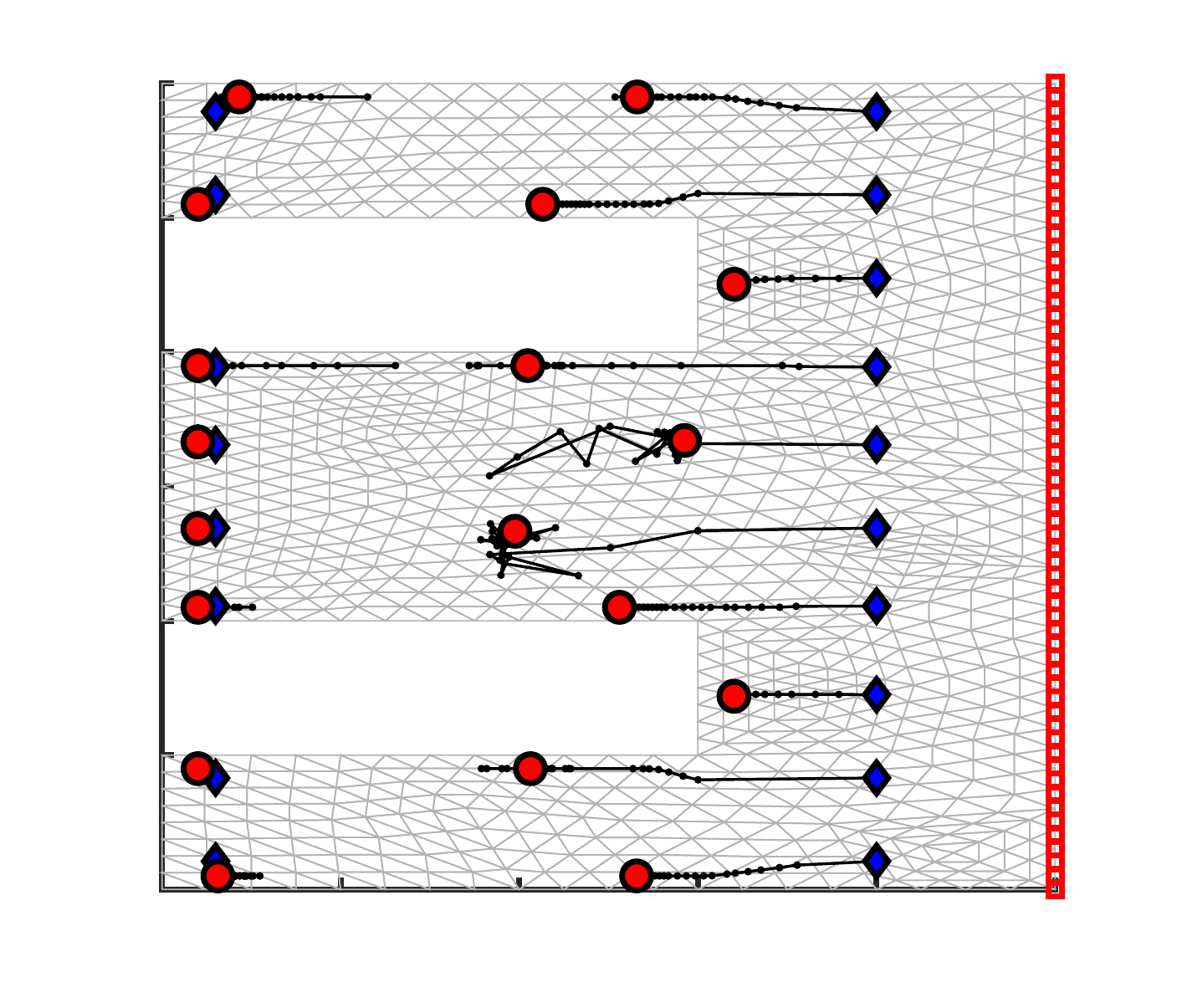}
   \includegraphics[width = 0.49\columnwidth]{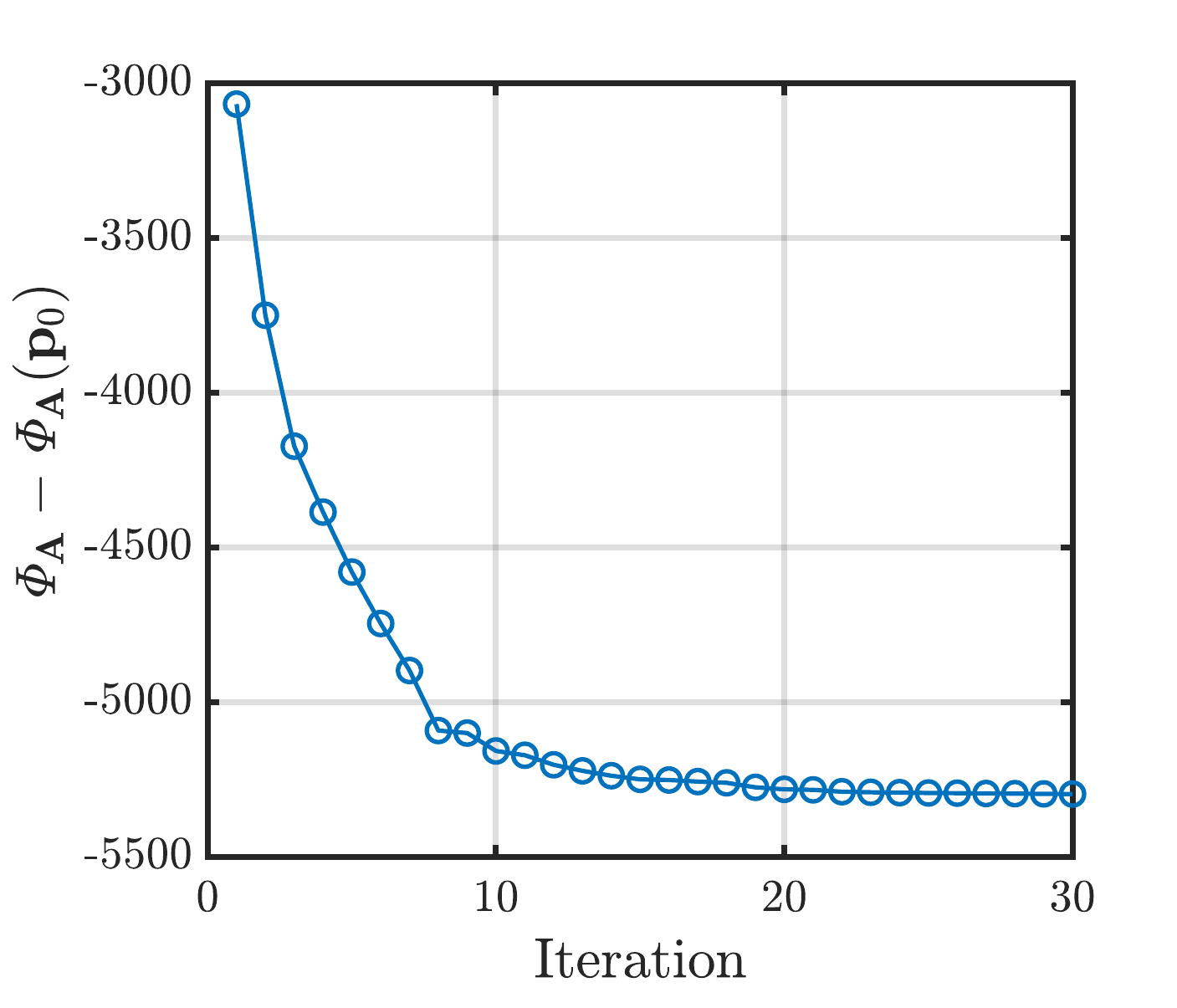}
    \caption{\small Experiment II. {\sc Left}: Sliding sensors iteration for finding the A-optimal configuration, showing initial and final positions of the sensors as diamonds and circles, respectively. The boundary measurement pixels are emphasized along the right-hand edge. {\sc Right}: Value of the shifted A-optimality target $\mathit{\Phi}_{\vA} - \mathit{\Phi}_{\vA}(\vp_0)$, where $\vp_0$ defines the initial sensor locations, at each iteration of the algorithm.}
    \label{fig:minimiter}
\end{figure}

Figure~\ref{fig:minimiter} illustrates the progress of the sliding sensors method for finding the A-optimal positions with $m_{\rm int} = 18$ internal sensors. The initial sensor positions are as indicated in Figure~\ref{fig:minimiter}, and the sliding sensors algorithm gradually pushes them away from $\mathcal{M}_{\rm bdry}$ towards other sections of $\partial \Omega_{\rm A}$. The evolution of the A-optimality target $\mathit{\Phi}_{\vA}$, with $\vA^T\!\vA = \vM$, is also visualized in Figure~\ref{fig:minimiter}.

\begin{figure}
    \includegraphics[width = 0.49\columnwidth]{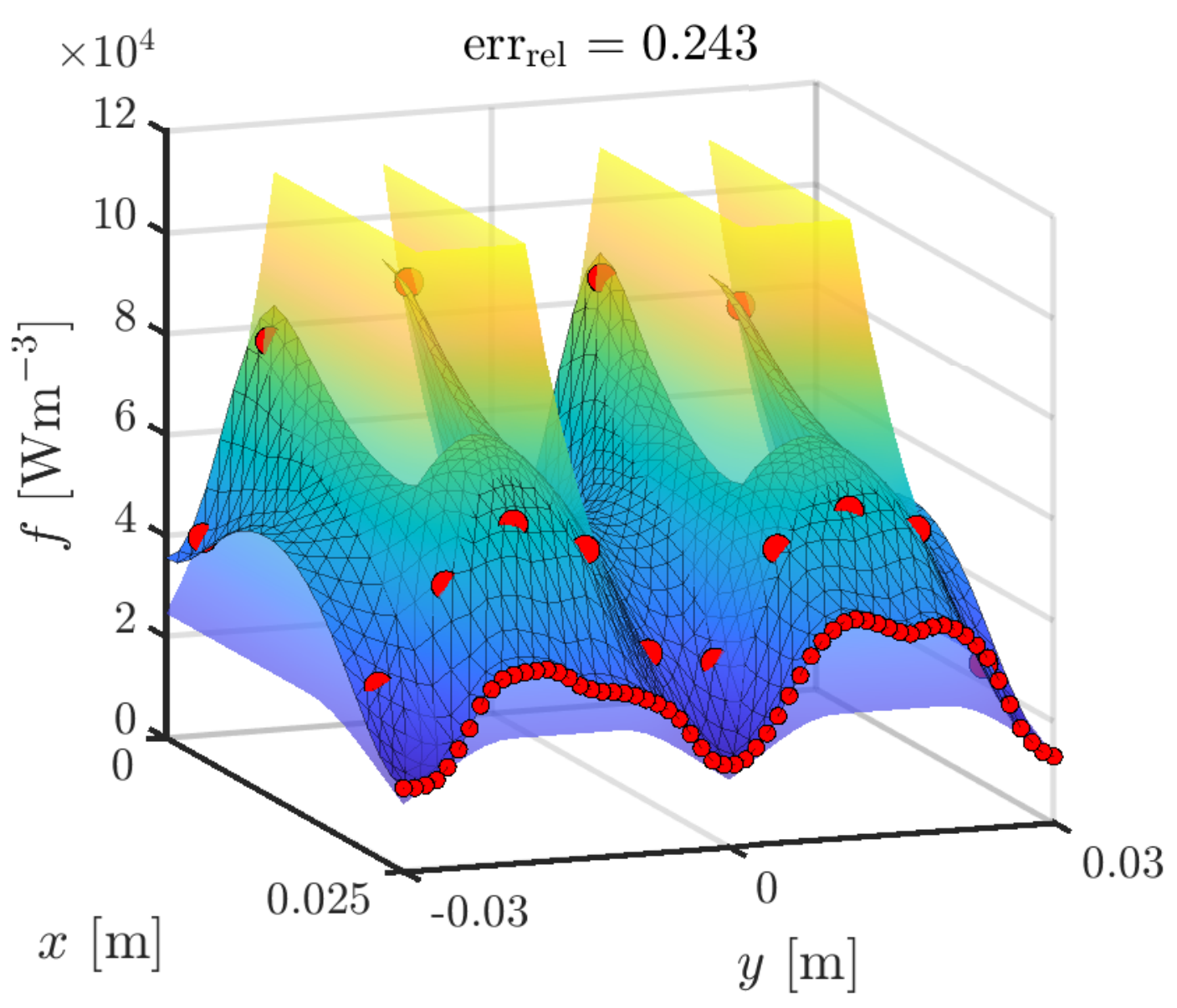}
    \includegraphics[width = 0.49\columnwidth]{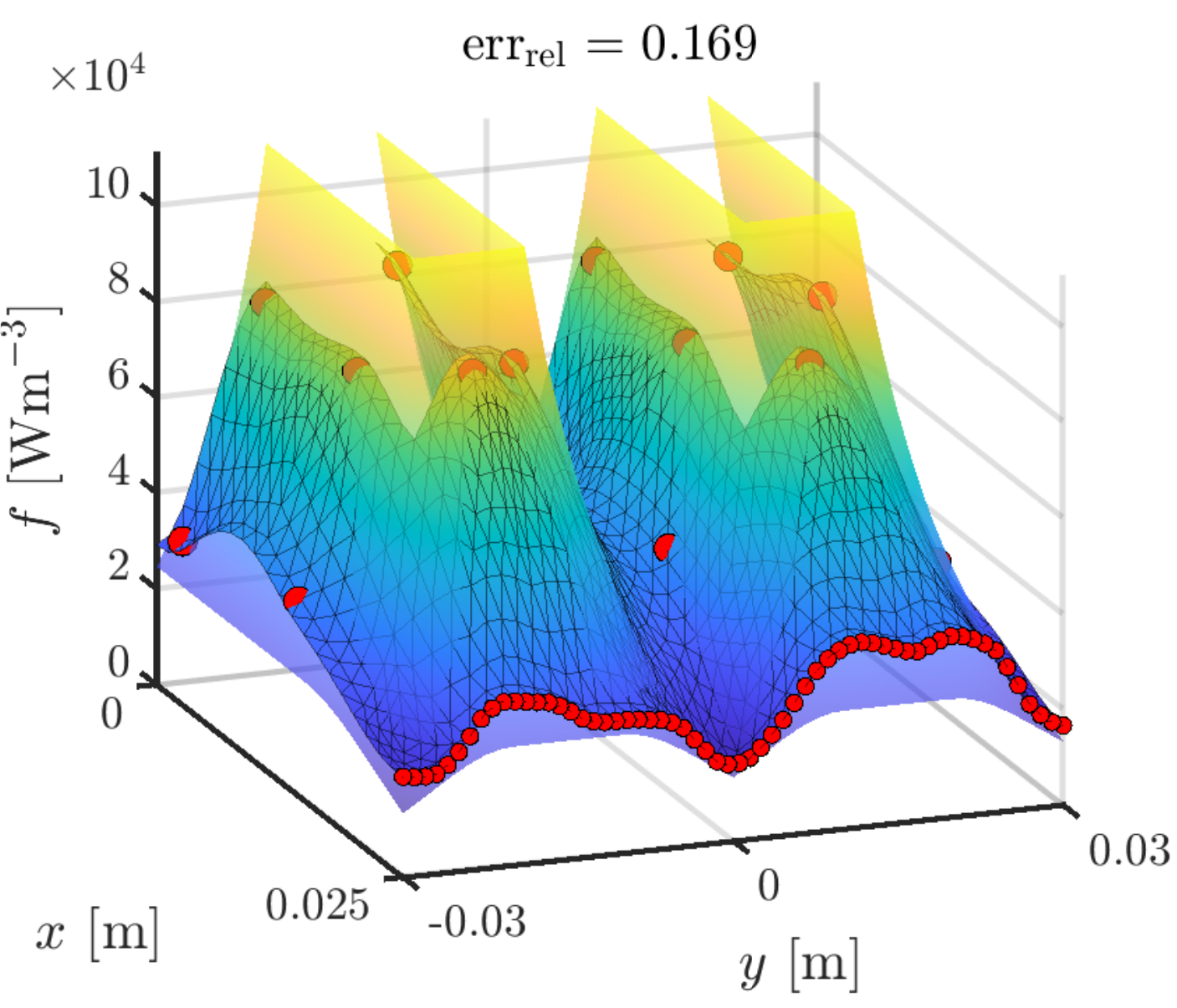}
    \caption{\small Experiment II. Reconstructions with 18 inner sensors. The gridded surface represents the reconstruction $\widehat f$ and the gridless surface the true source $f_{\mathrm{true}}$. {\sc Left}: Regular grid of sensors. {\sc Right}: A-optimized sensor positions.}
    \label{fig:2dlocations_nsens14}
\end{figure}
\begin{figure}
    \includegraphics[width = 0.49\columnwidth]{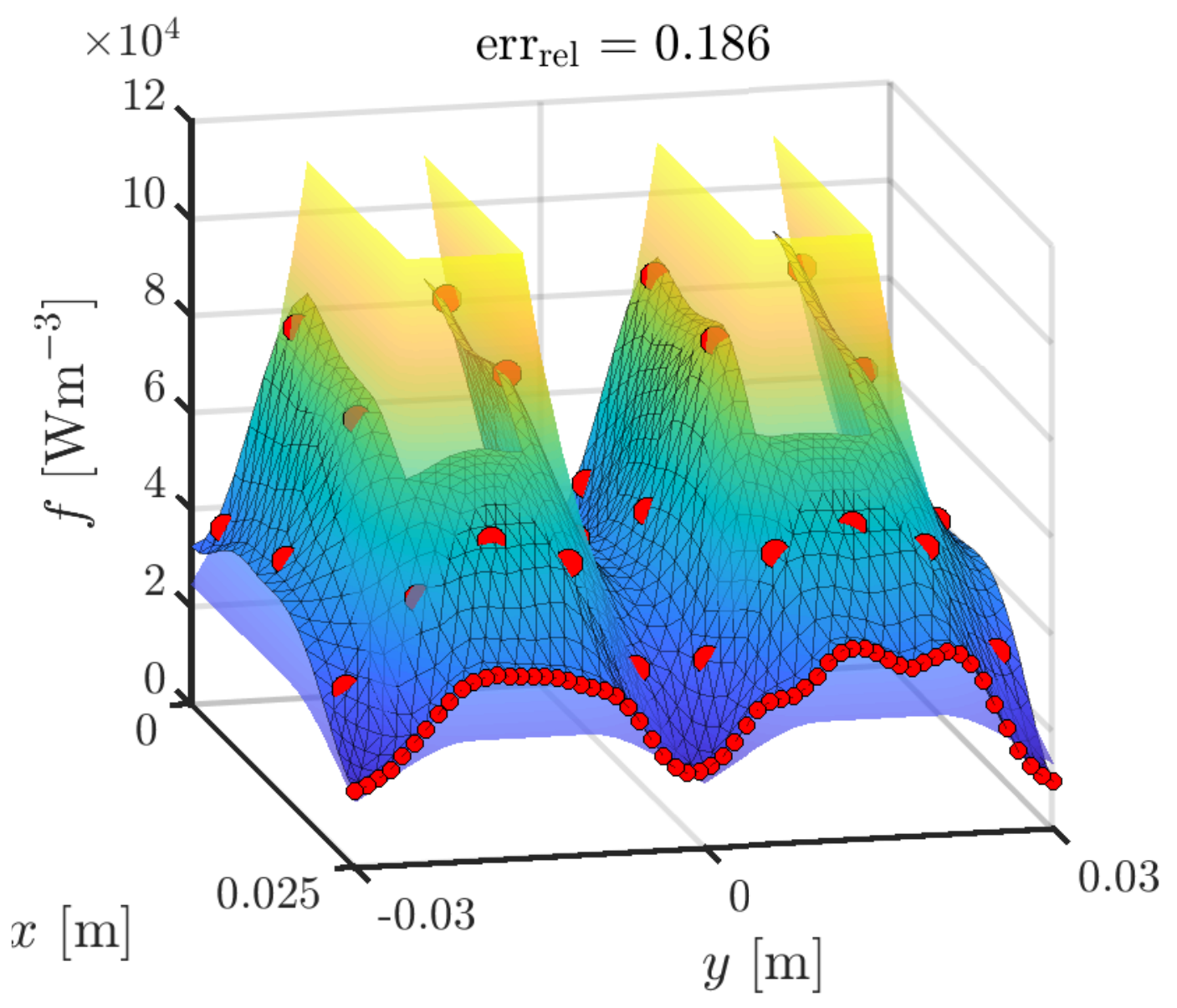}
    \includegraphics[width = 0.49\columnwidth]{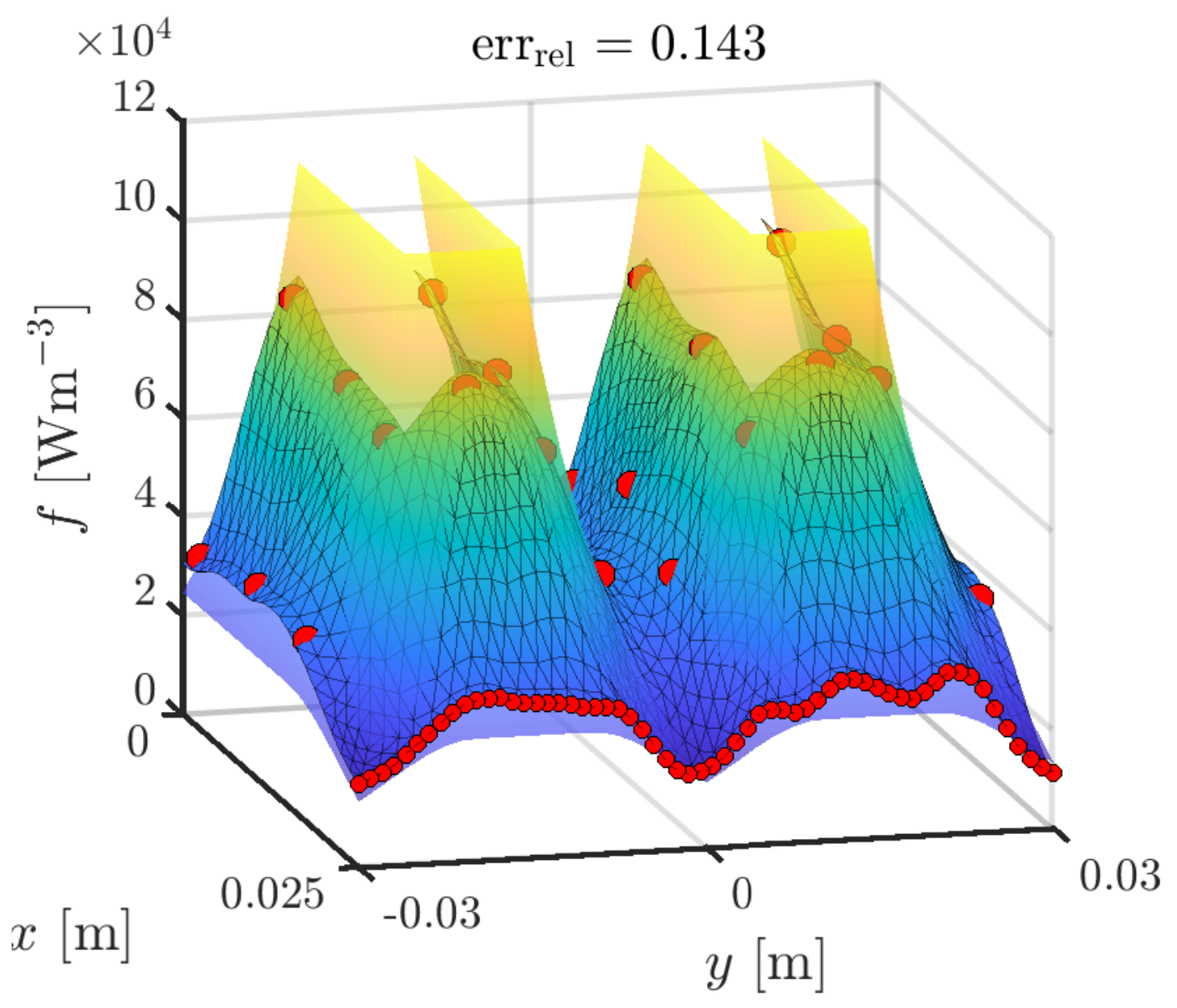}
    \caption{\small Experiment II. Reconstructions with 26 inner sensors. The gridded surface represents the reconstruction $\widehat f$ and the gridless surface the true source $f_{\mathrm{true}}$. {\sc Left}: Regular grid of sensors. {\sc Right}: A-optimized sensor positions.}
    \label{fig:2dlocations_nsens26}
\end{figure}

Reconstructions of the true source \eqref{eq:ftrue2} corresponding to a regular grid of $m_{\rm int}=18$ internal sensors and to the A-optimized set of $m_{\rm int}=18$ sensors from Figure~\ref{fig:minimiter} are visualized in Figure~\ref{fig:2dlocations_nsens14}. Figure~\ref{fig:2dlocations_nsens26} shows the corresponding results for $m_{\rm int}=26$. In both cases, the A-optimized sensor configuration produces significantly better reconstructions of the target source than the regular grid of internal sensors, as is quantified by the relative $L^2(\Omega)$ errors listed in Figures~\ref{fig:2dlocations_nsens14} and~\ref{fig:2dlocations_nsens26}. We observe that a lower number of sensors ($m_{\rm int}=18$) with A-optimized locations yields a better reconstruction than a considerably higher number  of sensors ($m_{\rm int}=26$) in a regular grid. Although not shown in Figures~\ref{fig:2dlocations_nsens14} and~\ref{fig:2dlocations_nsens26}, altogether omitting the internal sensors would result in intolerably bad reconstructions of the target source, yielding essentially a constant distribution deeper inside the domain.

\begin{figure}
	\includegraphics[width = 0.49\columnwidth]{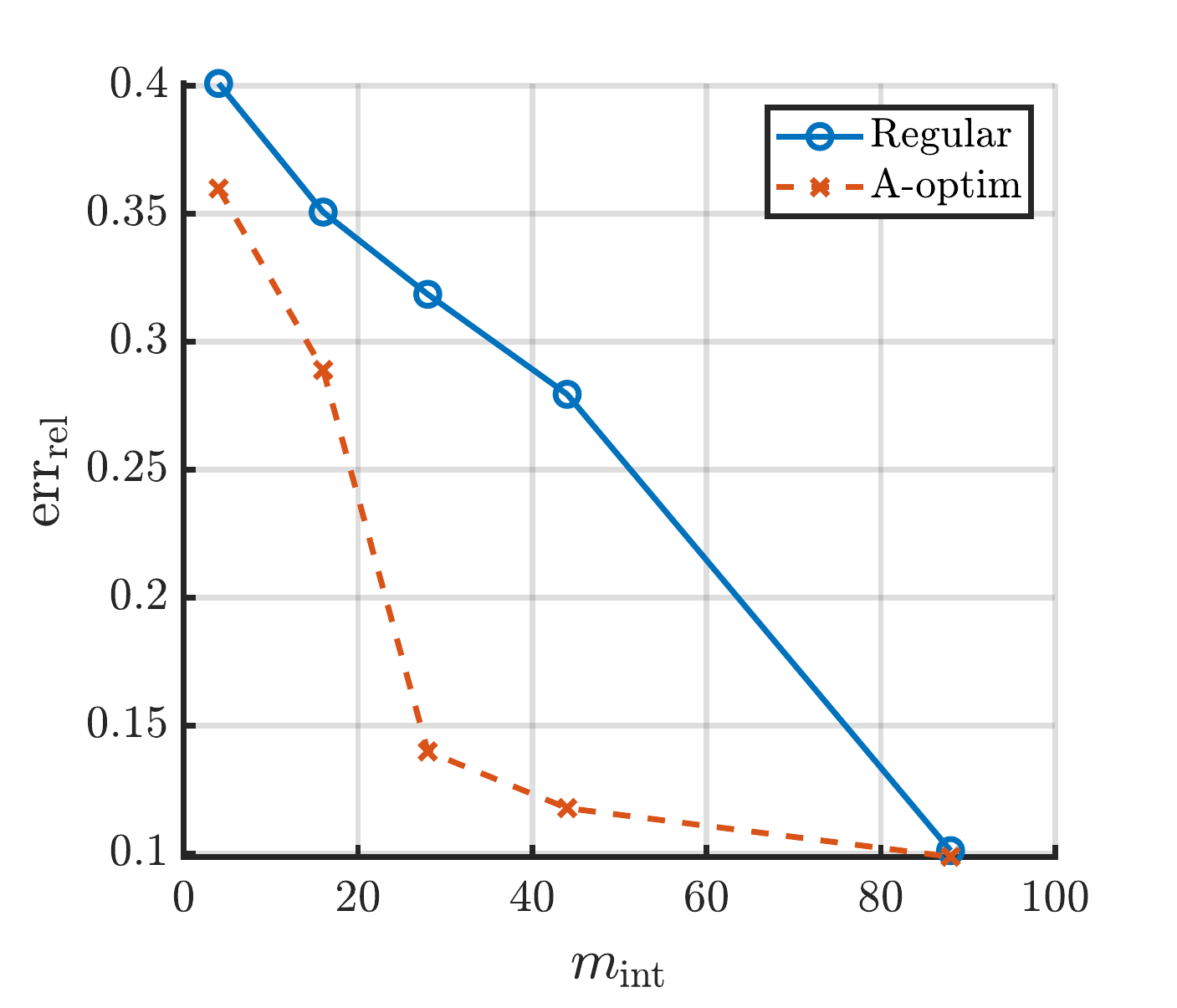}
	\includegraphics[width = 0.49\columnwidth]{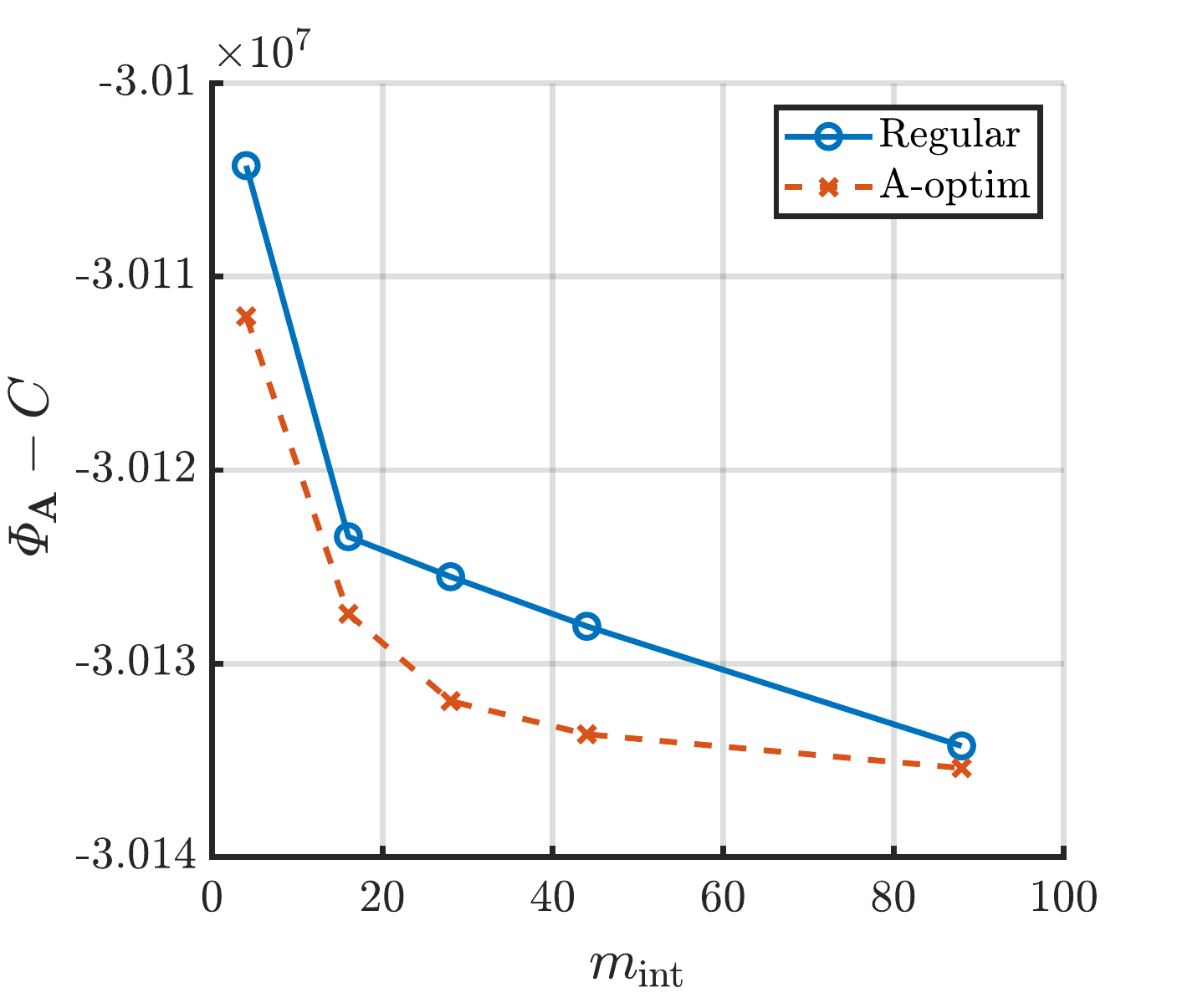}	
    \caption{\small Experiment II. {\sc Left}: The relative $L^2(\Omega)$ reconstruction error ${\rm err}_{\rm rel}$ as a function of the number of inner sensors $m_{\mathrm{int}}$ for regular grid positions and A-optimized positions. {\sc Right}: The corresponding values of the (shifted) A-optimality target functional $\mathit{\Phi}_{\!\vA} - C$ for regular grid positions and A-optimized positions.}
    \label{fig:2derrors}
\end{figure}

The relative reconstruction error~\eqref{eq:L2error} is plotted as a function of the number of internal sensors for certain regular grids and A-optimized configurations in Figure~\ref{fig:2derrors}. The figure also shows the corresponding values for the shifted A-optimality target functional $\mathit{\Phi}_{\vA} - C$, where the optimization-wise irrelevant constant $C$ is as in~\eqref{eq:tracecomp}. The A-optimized positions perform consistently better than the regular locations. At about $m_{\rm int} = 30$, the introduction of new sensors to the A-optimized configuration starts to yield less and less improvement in the reconstruction quality; for the considered regular positions, such a stagnation only happens at about $m_{\rm int} = 90$. Naturally, these numbers and also the overall performance of the A-optimized sensor positions compared to other sensor configurations depends,~e.g.,~on the (accuracy of the) prior, the form of the target source, the measurement geometry and the noise model. Be that as it may, at least in this simple test case, optimizing the sensor positions definitely seems worthwhile.

\subsection{Experiment III:~a semirealistic configuration in three dimensions}
\begin{figure}
    \includegraphics[width = 0.44\columnwidth]{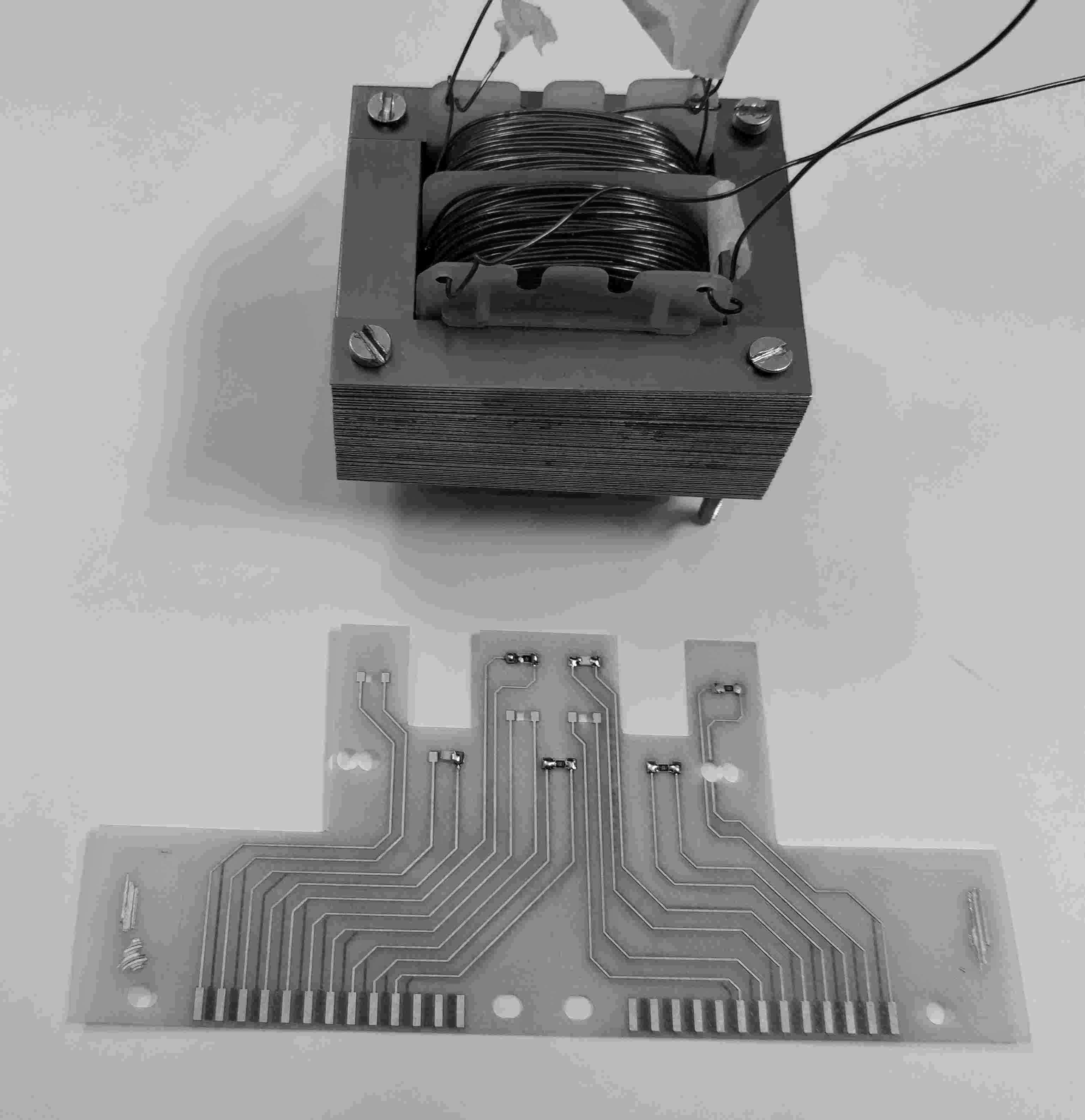}
    \includegraphics[width = 0.54\columnwidth, clip=true, trim = 0.3cm 0.6cm 0.0cm 0.9cm]{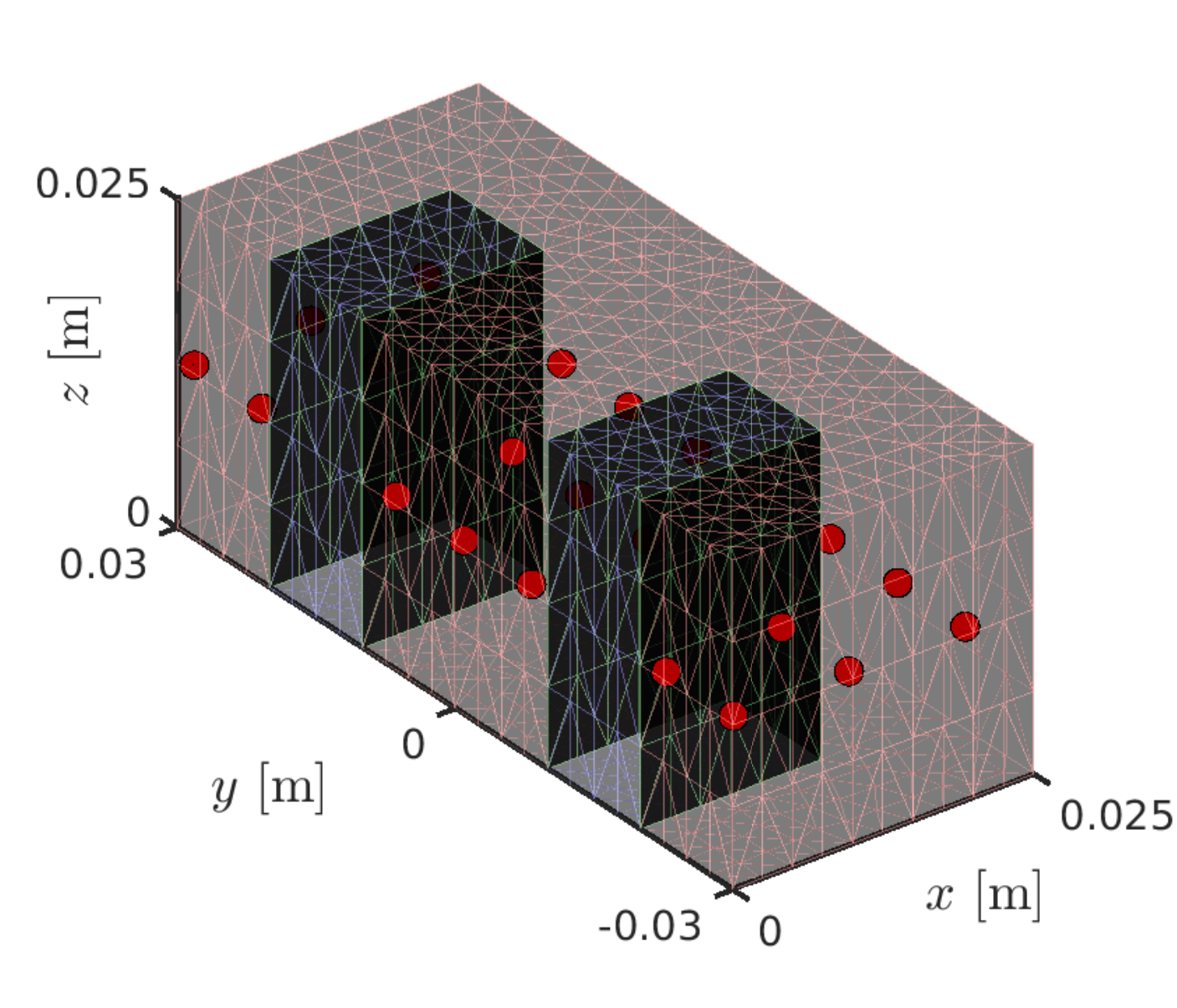}
    \caption{\small Experiment III. {\sc Left}: A small transformer and a prototype of a circuit board with internal temperature sensors. {\sc Right}: A simplified FE model for one half of the transformer, showing possible positions for the internal sensors at the mid cross section.}
    \label{fig:3dgeom}
\end{figure}

The considered three-dimensional geometry is obtained by extruding the geometry of Section~\ref{sec:2dtest} in the $z$-direction, resulting in the domain $\Omega = (0, 0.025) \times (-0.030, 0.030) \times (0, 0.025)$; see Figure~\ref{fig:3dgeom}. We still denote the iron core by $\Omega_{\rm A}$, the coils by $\Omega_{\rm B}$ and their interface by $\Gamma_{\rm AB}$ (cf.~Figure~\ref{fig:2dgeom}). The coil ends are omitted from the model for simplicity as we are dealing with simulated, not measured, data in this initial study. The Neumann, Robin and measurement boundaries are the natural extensions of the corresponding two-dimensional definitions in Section~\ref{sec:2dtest}, that is, $\Gamma_{\mathrm{N}} = \{ (x, y, z) \in \partial \Omega \ | \ x=0 \}$, $\Gamma_{\mathrm{R}} = \partial \Omega \setminus \Gamma_{\mathrm{N}}$ and $\mathcal{M}_{\rm bdry} = \{ (x, y, z) \in \partial \Omega \  | \  x=0.025 \}$.

\begin{figure}
    \includegraphics[width = 0.49\columnwidth]{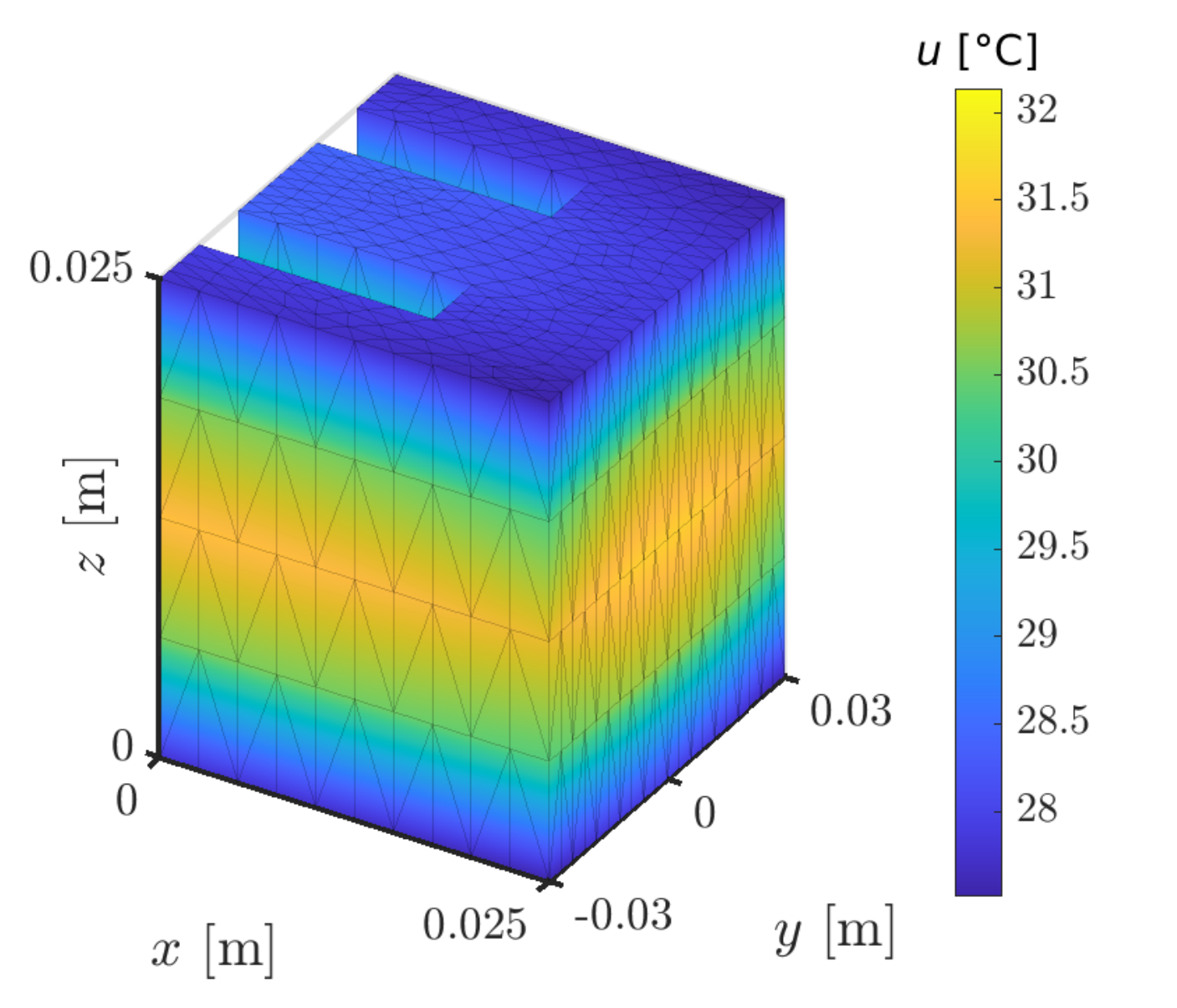}
    \includegraphics[width = 0.49\columnwidth]{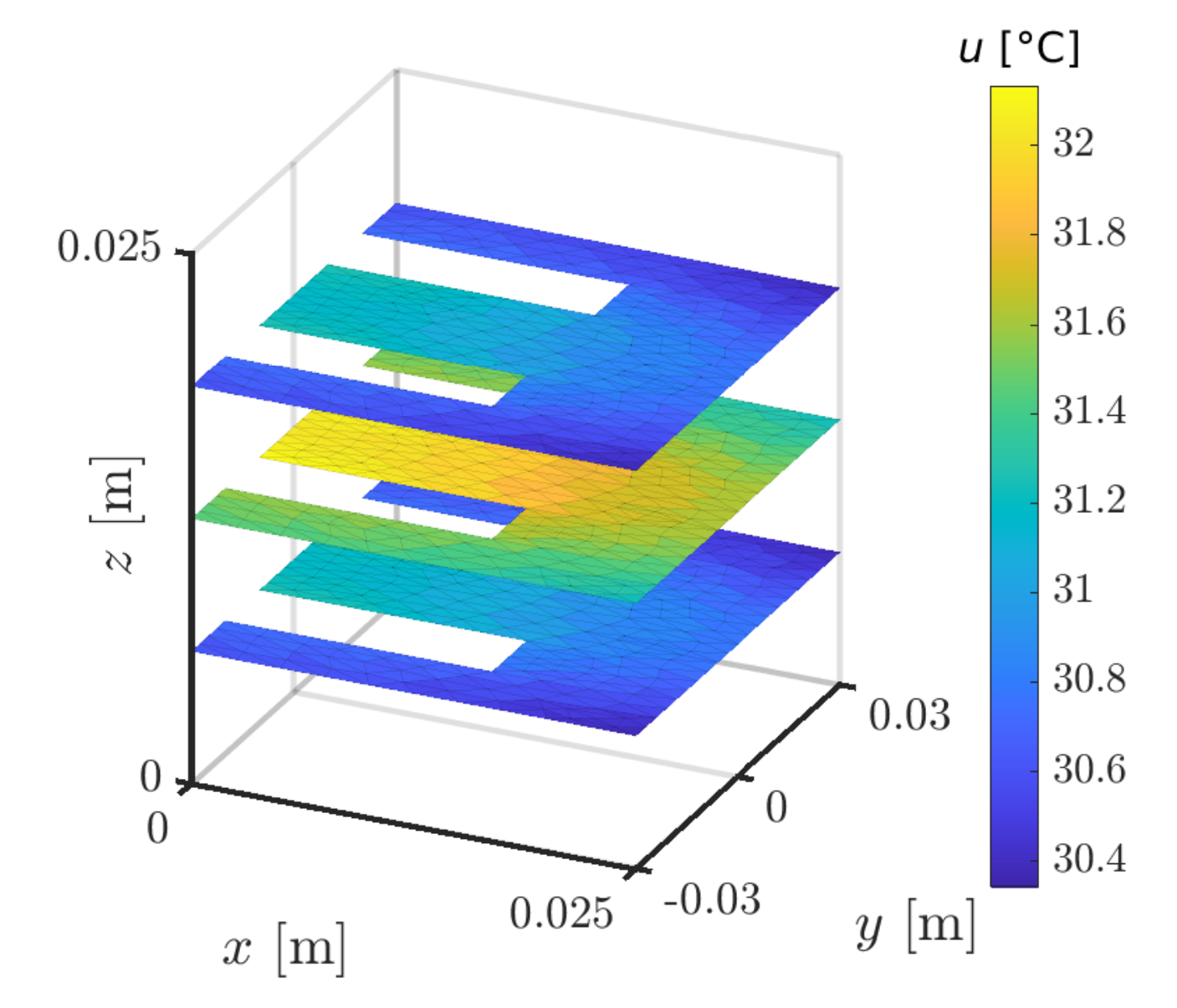}
    \caption{\small Experiment III. {\sc Left}: Surface temperature of the iron core in the steady state. {\sc Right}: Internal temperature of the iron core in the steady state.}
    \label{fig:3dtemp}
\end{figure}

The heat conductivity is modeled as anisotropic: in the windings, the conductivity is higher along the wires,~i.e.~in the $z$-direction, whereas in the stacked core steel sheets the conductivity is higher along the sheets,~i.e.~along the $xy$-plane. These effects are modeled by a piecewise constant diagonal heat conductivity matrix that takes the value $\kappa = \operatorname{diag}(\kappa_{{\rm A},x}, \kappa_{{\rm A},y}, \kappa_{{\rm A},z})$ in $\Omega_{\rm A}$ and $\kappa = \operatorname{diag}(\kappa_{{\rm B},x}, \kappa_{{\rm B},y}, \kappa_{{\rm B},z})$ in $\Omega_{\rm B}$, with $\kappa_{\mathrm{A}, x} = \kappa_{\mathrm{A}, y} = 26$, $\kappa_{\mathrm{A}, z} = 0.6$ and $\kappa_{\mathrm{B}, x} = \kappa_{\mathrm{B}, y} = 10$, $\kappa_{\mathrm{B}, z} = 400$.  The parameter values for the forward problem \eqref{eq:fwdPDE} and the model for the insulating layer \eqref{eq:jump} are otherwise as in Section~\ref{sec:2dtest}. The true source $f_{\mathrm{true}}$ is chosen to be as in~\eqref{eq:ftrue2}, extended uniformly in the $z$-direction. The forward solution for $f_{\mathrm{true}}$ is illustrated in Figure~\ref{fig:3dtemp}, which shows the steady state temperature in the iron core and on its boundary. The simulated temperature rise of 32\textdegree{}C is in a realistic range.

The temperature on the face $\mathcal{M}_{\rm bdry}$ is measured by a thermal camera at $m_{\mathrm{bdry}} = 361$ pixels, and $m_{\rm int}= 23$ sensors are inserted onto the mid cross section of the object; see Figure~\ref{fig:3dgeom}. The temperature is recorded at $m_t = 20$ observation times distributed uniformly over $[T/80,T]$, where $T=10^4$ (approximately three hours), so the total number of boundary and internal data points are $m_\bdry m_t = 7220$ and $m_{\rm int} m_t = 460$, respectively, leading altogether to $m=7680$. The measurements are contaminated by $0.1$\% of zero-mean additive Gaussian noise. It is expected that the heat source varies more quickly in the $xy$-plane than in the $z$-direction, and this extra information is encoded as an anisotropic parameter $\alpha = \operatorname{diag}(\alpha_x, \alpha_y, \alpha_z)$ in the prior covariance~\eqref{eq:prior}, with $\alpha_x = \alpha_y = 3 \cdot 10^{-8}$, $\alpha_z = 3 \cdot 10^{-6}$, $\beta = 8 \cdot 10^{-6}$.

\begin{figure}
    \includegraphics[width = 0.49\columnwidth]{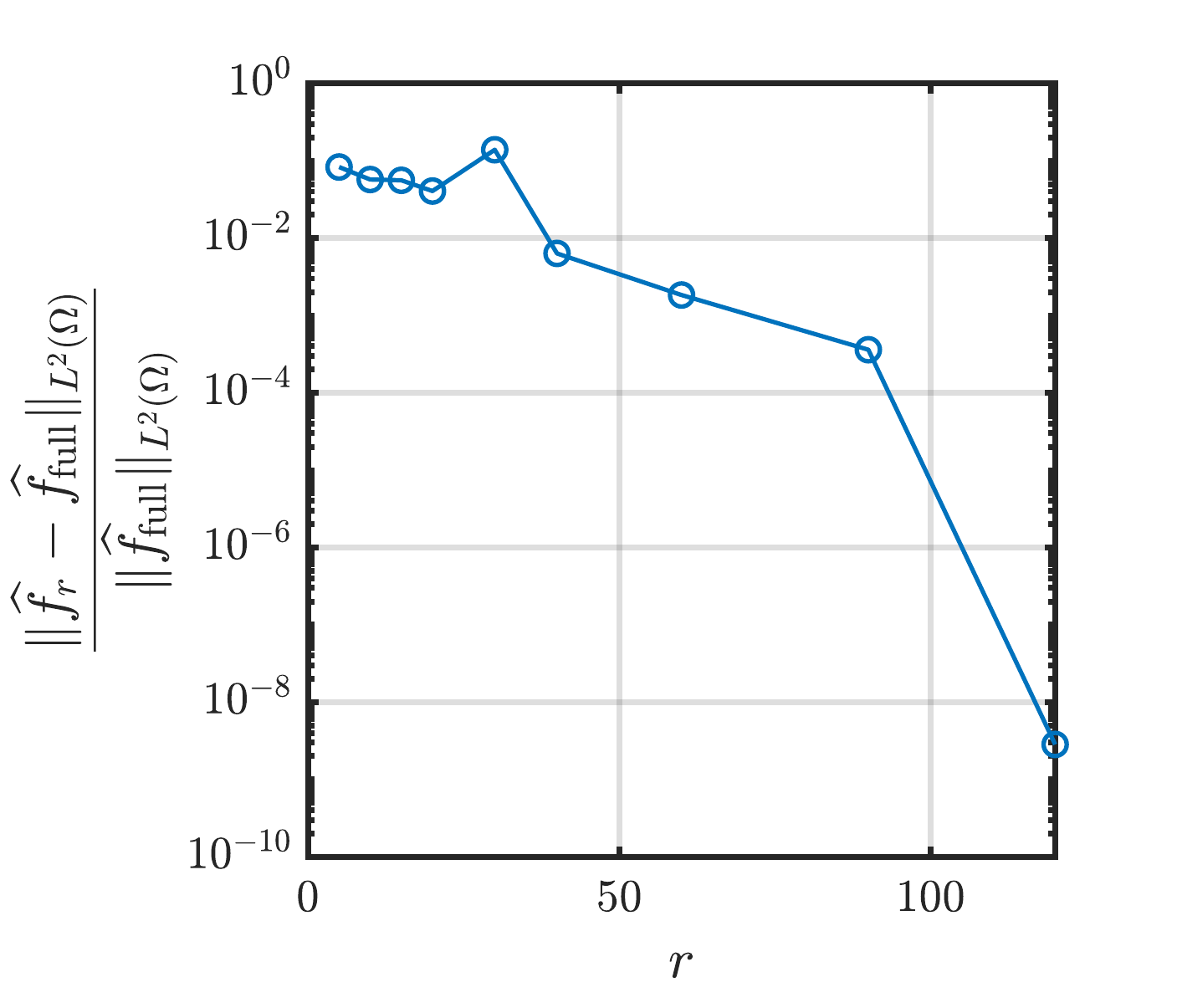}
    \includegraphics[width = 0.49\columnwidth]{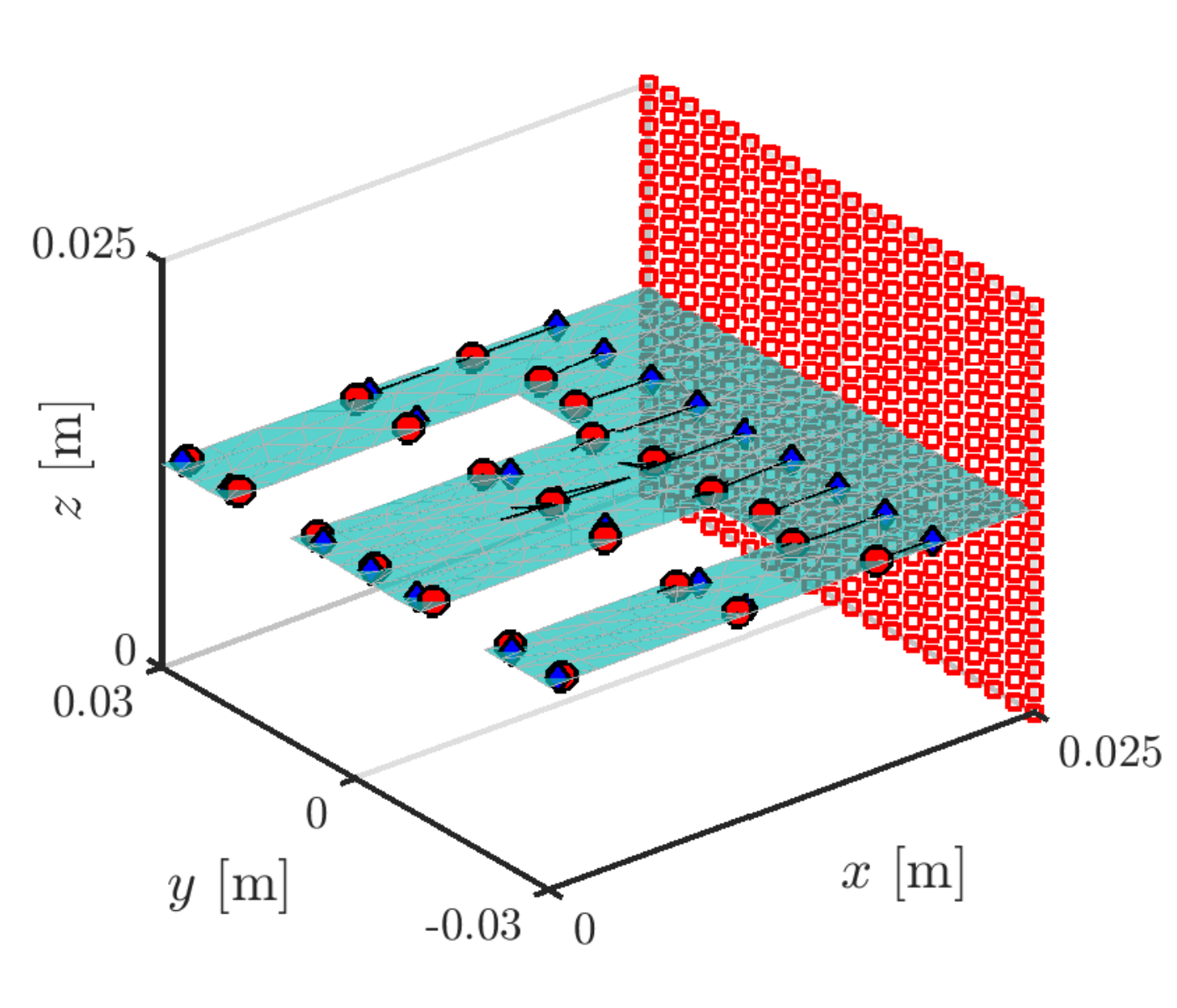}
    \caption{\small Experiment III. {\sc Left}: Relative $L^2(\Omega)$ discrepancy in the reconstruction caused by the low-rank approximation as a function of the reduced dimension $r$ in the boundary measurement model that has the full dimension $m_{\mathrm{bdry}} m_t = 7220$. {\sc Right}: The progress of the sliding sensors algorithm with the initial and final internal sensor positions marked by diamonds and balls, respectively. The boundary sensors are marked by pixels.}
    \label{fig:r_vs_error}
\end{figure}

To begin with, let us investigate how the error introduced by the low-rank approximation of the boundary measurements described in Section~\ref{sec:lowrank} propagates into the reconstruction. We assume the internal sensors form a regular grid and consider a small problem that corresponds to a discretization of the studied geometry with only $n = 3590$ FE degrees of freedom. As a consequence, all matrices involved in forming the reconstruction can be explicitly computed and applied without any approximations; on the negative side, the sparse FE discretization certainly causes considerable numerical errors that are ignored here for simplicity. The left-hand image of Figure~\ref{fig:r_vs_error} shows the relative discrepancy $\| \widehat{f}_r - \widehat{f}_{\mathrm{full}} \|_{L^2(\Omega)} / \| \widehat{f}_{\mathrm{full}} \|_{L^2(\Omega)}$ as a function of $r$. Here the benchmark reconstruction $\widehat{f}_{\mathrm{full}}$ is computed using the full (discretized) boundary measurement operator $\vF_\bdry$, whereas $\widehat{f}_{r}$ corresponds to the associated low-rank approximation $\widetilde \vF_\bdry \in \R^{r \times n}$. According to this simple low-dimensional test, it seems to be possible to achieve a reduction of, say, $m_{\mathrm{bdry}} m_t = 7220 \hookrightarrow r=120$ in the dimension of the boundary measurement operator without any relevant loss of information from the standpoint of the considered inverse source problem.

\begin{figure}
    \includegraphics[width = 0.49\columnwidth]{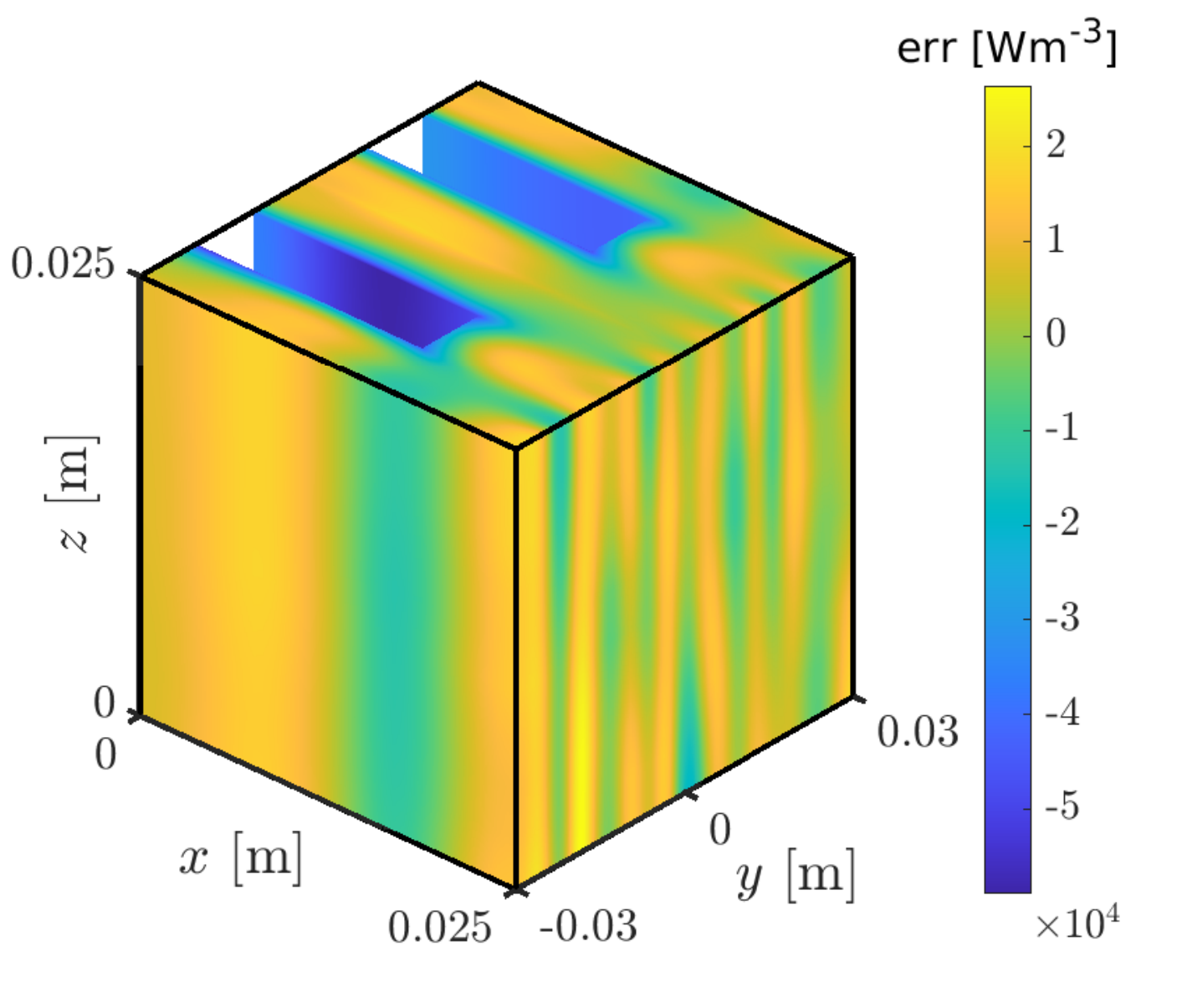}
    \includegraphics[width = 0.49\columnwidth]{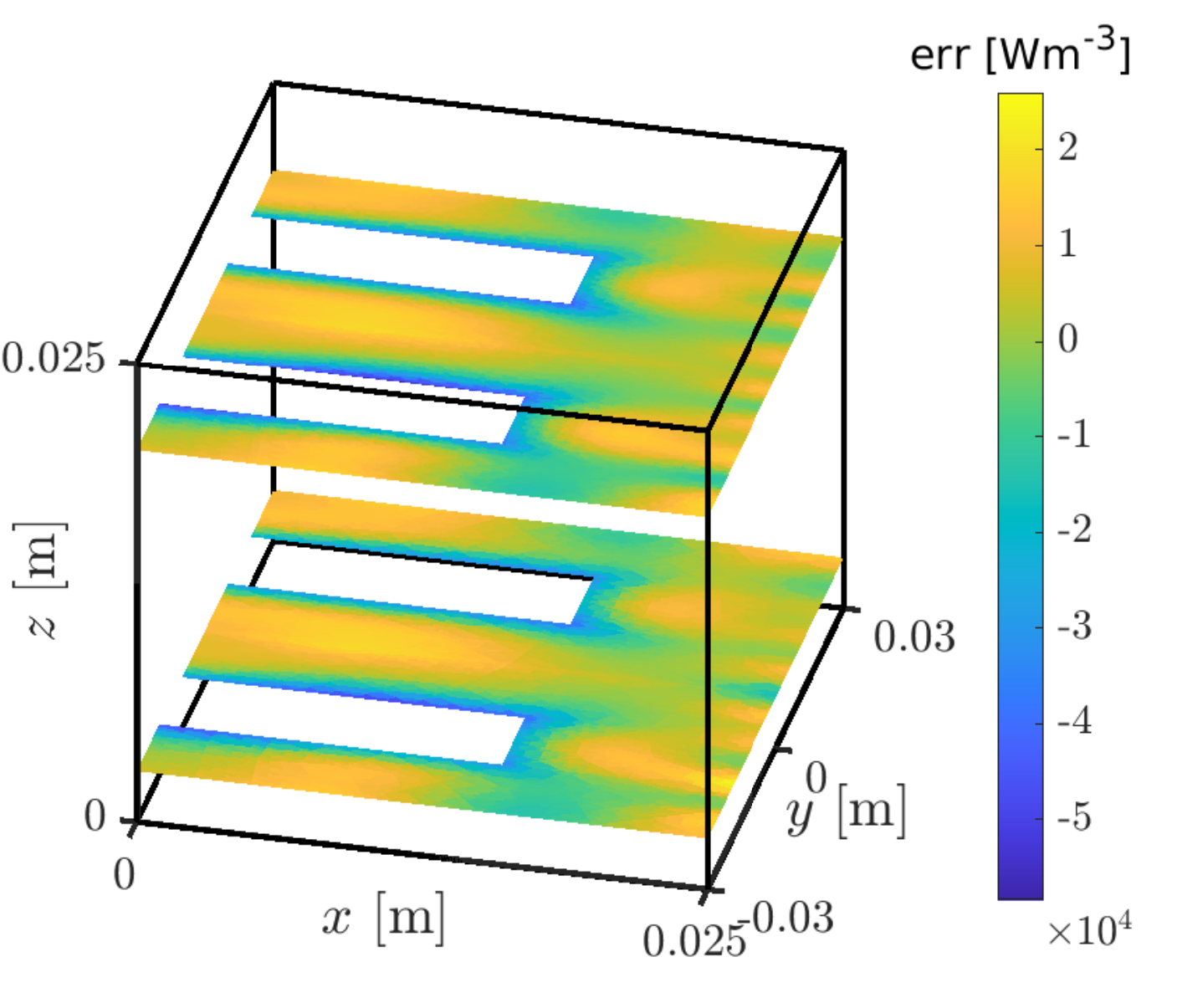}
    \caption{\small Experiment III.  Reconstruction error $\widehat f - f_{\mathrm{true}}$ with the A-optimized sensor locations. {\sc Left}: Boundary of the iron core. {\sc Right}: Horizontal cross sections of the iron core.}
    \label{fig:test3_errors}
\end{figure}

Next, the A-optimality seeking sliding sensors algorithm is combined with the dimension reduction method. As we switch to a FE model with a realistic number of degrees of freedom $n = 49\,310$, the full matrix $\vF \in \R^{m \times n}$ can no longer be constructed and used in practical computations. Hence, a low-rank approximation corresponding to $m_{\mathrm{bdry}} m_t = 7220 \hookrightarrow r=120$ is computed for $\vF_\bdry$. To mimic the real application, the inner sensors are restricted to lie on a single core sheet, i.e., their positions can only be optimized in the $x$ and $y$-directions. The initial and A-optimized sensor locations are shown in Figure~\ref{fig:r_vs_error}. A comparison with Figure~\ref{fig:minimiter} indicates that the optimal positions are qualitatively similar to those found in the two-dimensional case. Figure~\ref{fig:test3_errors} shows the reconstruction error $\widehat f - f_{\mathrm{true}}$ on certain parts of the boundary $\partial \Omega_{\mathrm{A}}$ and at two horizontal cross sections of the domain. It can be seen that the reconstructed source is severely underestimated close to the coils, where the true source and its gradient are the largest. The relative $L^2(\Omega)$ reconstruction error with the A-optimized locations is ${\rm err}_{\rm rel} = 0.262$, whereas with the initial locations the error is ${\rm err}_{\rm rel} = 0.317$.

\begin{table}
\caption{Experiment III. Computational cost.}
\label{tab:comptime}
\begin{tabular}{c||c|c|c}
$n$ & \makecell{Low-rank \\ computation time} & \makecell{A-optimization \\ computation time} & \makecell{Total number of \\ forward solutions} \\
\hline
1795  & 65\,s & 110\,s  & 9089 \\
9093  & 965\,s & 1316\,s & 8905 \\
49310 & 12\thinspace{}024\,s & 15\thinspace{}629\,s & 9020
\end{tabular}
\end{table}

Finally, we briefly discuss the computational cost for problems of different sizes. The  above described three-dimensional setting is discretized using three different refinement levels for the employed FE mesh. The number of degrees of freedom in the discretization $n$, the computational time required to construct a low rank approximation with $r=120$, the computational time to perform 20 iterations of the sliding sensors algorithm based on the reduced model, and the total number of forward solutions are listed in Table~\ref{tab:comptime}. For the two finest meshes, the A-optimized sensor locations are almost identical (although not shown here), suggesting that the corresponding discretizations are sufficiently accurate. Although not listed explicitly in Table~\ref{tab:comptime}, the number of forward solutions required to construct an accurate enough low-rank approximation seems to be independent of $n$ (cf.~\cite{Flath2011}). Using the techniques described in Section~\ref{sec:eval_object}, the computational cost related to the needed evaluations of matrix traces and their derivatives is negligible compared to the cost of the forward solutions, meaning that the computation times associated to the sliding sensors algorithm are mainly related to the latter. Because the computation of forward solutions for different sources is trivial to parallelize, the sliding sensors method should be readily scalable to be used for even larger $n$.

\section{Discussion and conclusion}
\label{sec:discussion}
This work introduced a computational framework for a simplified linear inverse problem that models the determination of the iron loss field inside an electric machine. Although many issues that are important in practice, such as nonlinearities in the forward model, the proper choice of a prior distribution and a realistic noise model, were ignored, our precursory numerical studies anyway demonstrate that temperature measurements on the boundary of the imaged machine do not suffice for accurate enough reconstruction of an internal heat source. Moreover, both theoretical and numerical considerations clearly indicate that measuring time transient temperature data facilitates the reconstruction process.

In our numerical experiments, the use of (approximately) A-optimal positions for the temperature sensors inserted between the stacked core steel sheets had a significant positive effect on the reconstruction quality. In particular, the introduced sliding sensors algorithm for predicting A-optimal configurations  was able to considerably improve the performance of measurement setups. From a computational standpoint, the repetitive solution of the parabolic forward problem constitutes the most expensive part of the sliding sensors method. However, since such forward solutions are easily parallelizable, we anticipate the proposed method can, in fact, be used for much larger problems than the ones tackled in this paper.

The prior model used in our numerical tests corresponds to combined $L^2$--$H^1$-regularization. Such a traditional approach was chosen due to its simplicity and computational efficiency. Constructing a prior distribution that captures the generic properties of (all reasonable) iron loss models is a research topic of its own. Another potential direction for future research would be to design statistical tests based on temperature measurements for confirming or rejecting proposed iron loss models.

In the real world, constructing an accurate heat conduction model for an electric machine is nontrivial as it involves a nonlinear parabolic PDE with temperature-dependent coefficients. In addition, the parameter $h$ appearing in~\eqref{eq:robin} is not precisely known in the real application. As a consequence, the modeling errors would most likely dominate the measurement accuracy if our simplified linear model were employed in practice. The introduced mismodeling could to a certain extent be incorporated in the noise model, or alternatively, the nonlinearities could be properly modeled and $h$ could be included as an additional unknown in the reconstruction process. The former corresponds to the so-called {\em approximation error approach} \cite{Kaipio2004statistical} and the latter would lead to a more complicated nonlinear inverse problem; both of these alternatives are interesting topics for future studies.

\appendix

\section{Adjoints of the forward operators}
\label{sec:adjoints}
Let us first deduce the adjoint for the idealized forward operator $\mathcal{F}: L^2(\Omega) \to L^2(\mathcal{M}_{\rm bdry} \times (0,T))$ defined in \eqref{eq:forwardop}. To this end, consider the parabolic `backwards' initial/boundary value problem
\begin{subequations} \label{eq:adjPDE}
\begin{align}
	\rho \partial_t w + \nabla \cdot (\kappa \nabla w) &= 0  \qquad	 \text{ in } \Omega \times (0, T),  \label{eq:heateq2} \\
	\nu \cdot \kappa \nabla w + h  w &= g \qquad	 \text{ on } \Gamma_{\mathrm R} \times (0, T),  \label{eq:robin2}\\
	 \nu \cdot \kappa \nabla w &= g \qquad	 \text{ on } \Gamma_{\mathrm N} \times (0, T),  \label{eq:neumann2}\\
	w &= 0 	 \qquad \text{ on } \Omega \times \{ t = T \}, \label{eq:initcond2}
\end{align}
\end{subequations}
where the coefficients are the same as in \eqref{eq:fwdPDE} and $g \in L^2(\mathcal{M}_{\rm bdry} \times (0,T))$ is interpreted as an element of $L^2(\partial \Omega \times (0,T))$ via zero continuation. The variational formulation of \eqref{eq:adjPDE} is
\begin{equation}
  \label{eq:weak_dual}
- \langle \partial_t w, v \rangle_\rho +  a(w, v) =   \int_{\mathcal{M}_{\rm bdry}} \!\!\!  g v \, {\rm d} S \qquad \text{ for all } v \in H^1(\Omega),
\end{equation}
with a vanishing `initial condition' at $t=T$. Via the change of variables $\tau = T - t$, it is easy to see that \eqref{eq:weak_dual} has a unique solution in $\mathcal{H}^1((0, T); \Omega)$ by virtue of the standard theory on parabolic partial differential equations \cite[Chapter~10]{Renardy93}.

Recalling \eqref{eq:forwardop} and comparing \eqref{eq:weakform} with \eqref{eq:weak_dual}, we get
\begin{align}
  \label{eq:dual_deduc}
  \int_{0}^T\int_{\mathcal{M}_{\rm bdry}} g \, \mathcal{F}f  \, {\rm d} S \, {\rm d} t &= -  \int_{0}^T \big(\langle \partial_t  w, u \rangle_\rho -  a(w, u) \big) {\rm d} t   \nonumber \\
  &=  - \int_{0}^T \big( \langle \partial_t  w, u \rangle_\rho + \langle \partial_t  u, w \rangle_\rho  \big) {\rm d} t + \int_{0}^T (f, w)_{L^2(\Omega)} \, {\rm d} t.
\end{align}
Since $w \in \mathcal{H}^1((0, T); \Omega)$ and  $u \in \mathcal{H}^2((0, T); \Omega)$ as noted before \eqref{eq:Psit}, it is straightforward to deduce that $(w, u )_\rho \in W^{1,1}(0,T)$ with the weak derivative
$$
\partial_t (w, u )_\rho = \langle \partial_t w, u \rangle_\rho + \langle \partial_t u, w \rangle_\rho.
$$
Hence, the first term on the right-hand side of \eqref{eq:dual_deduc} vanishes by virtue of the fundamental theorem of calculus accompanied by \eqref{eq:initcond} with $u_{\rm init} = 0$ and \eqref{eq:initcond2}; recall that we systematically identify the elements of $\mathcal{H}^1((0, T); \Omega)$ with their time-continuous representatives in $\mathcal{C}([0,T]; L^2(\Omega))$ \cite[Lemma~10.4]{Renardy93}. As \eqref{eq:dual_deduc} holds for all $f\in L^2(\Omega)$ and $g \in L^2(\mathcal{M}_{\rm bdry} \times (0,T))$, it follows that
\begin{equation}
  \label{eq:adjoint_forwardop}
\mathcal{F}^*: \left\{
  \begin{array}{l}
   {\displaystyle g \mapsto \int_0^{T} w \, {\rm d} t }, \\[4mm]
    L^2\big( \mathcal{M}_\bdry \times (0,T)\big) \to L^2(\Omega),
  \end{array}
  \right.
\end{equation}
defines the adjoint of the idealized forward operator.

By formally repeating the above calculations with the combined boundary/internal source
\begin{equation}
  \label{eq:delta_source}
g = \sum_{i=1}^{m_s} \sum_{j=1}^{m_t} y_{ij} \frac{1}{|S_i|} \chi_{i} \delta_{j},
\end{equation}
where $y = [y_{ij}] \in \R^{m}$, $\chi_{i}$ is the characteristic function of the $i$th sensor $S_i$, and $\delta_{j}$ is a Dirac delta functional in time supported at $t_j$, one arrives at the conclusion
$$
y^{T} \! Ff =  \Big(f, \int_0^{T} w \, {\rm d}t \Big)_{L^2(\Omega)} .
$$
That is, the adjoint of the realistic forward map $F: L^2(\Omega) \to \R^m$ from \eqref{eq:Fdisc} is formally given by
\begin{equation}
  \label{eq:realistic_adjoint}
F^*: y \mapsto \int_0^{T} w \, {\rm d} t,
\end{equation}
where $w$ is the solution to \eqref{eq:adjPDE} with the source \eqref{eq:delta_source}. To properly prove this claim, one should consider the solvability and regularity of \eqref{eq:adjPDE} for `time-irregular' sources of the form \eqref{eq:delta_source}. Be that as it may, we content ourselves here with commenting that all ambiguity in the above deduction of $F^*$ completely disappears if the Dirac deltas in \eqref{eq:delta_source} are replaced by (localized) weight functions in $L^2(0,T)$ with unit masses, and analogously the time evaluations of $\mathit{\Psi}$ in \eqref{eq:Fdisc} are replaced by the weighted time-averages of $\mathit{\Psi}$ defined by those same functions; cf.~Remark~\ref{remark:weight}.

In any case, \eqref{eq:realistic_adjoint} gives an intuitive explanation for the computational advantage in forming the discretized forward map $\vF$ via transposition as in \eqref{eq:Ftransp}: To evaluate $F^*$ for all sources of the form
$$
g_j = \frac{1}{|S_i|} \chi_{i} \delta_{j}, \qquad j = 1, \dots, m_t,
$$
with a fixed sensor index $1 \leq i \leq m_s$, one only needs to solve \eqref{eq:adjPDE} once over the whole interval $(0,T)$. Indeed, because the solution to \eqref{eq:adjPDE} with $g = g_j$ is clearly identically zero on the interval $(t_j, T]$, the time integration in \eqref{eq:realistic_adjoint} makes it equivalent to place a delta source at time $t=T$ and solve \eqref{eq:adjPDE} backwards in time up to $t = T-t_j$ as it is to solve \eqref{eq:adjPDE} backwards over the whole interval $(0,T)$ with the delta-like source only activating at $t_j$. Combining this logic with the change of variables $\tau = T - t$ in \eqref{eq:adjPDE} explains the advantageous structure of \eqref{eq:Ftransp}.

  \section{Sensor location derivative}
  \label{sec:derivative}
Let us assume that a sensor is only allowed to move within an open subdomain $D_0 \subset \Omega$ where the heat conductivity is Lipschitz continuous, i.e., $\kappa|_{D_0} \in [\mathcal{C}^{0,1}(D_0)]^{d \times d} \cap [L^{\infty}_+(D_0)]^{d \times d}$. Since the weak solution to \eqref{eq:fwdPDE} belongs to $\mathcal{H}^2((0,T); \Omega) \subset \mathcal{C}^1([0,T]; L^2(\Omega))$ as noted before \eqref{eq:Psit}, for any $t \in [0,T]$ it holds that
$$
\nabla \cdot (\kappa \nabla u(\, \cdot \,, t)) = \rho \partial_t u(\, \cdot \,, t) - f \in L^2(\Omega)
$$
without any extra assumptions on $\rho \in L^\infty_+(\Omega)$ or $f \in L^2(\Omega)$. In consequence, due to interior regularity of solutions to elliptic partial differential equations~\cite{Necas2012}, $u(\, \cdot \,, t)|_{D} \in H^2(D)$ for any domain $D \subset \subset D_0$. In the following, we assume $D$ has a Lipschitz boundary and otherwise satisfies the above listed properties.

As the realistic measurements modeled by \eqref{eq:Fdisc} only depend on the locations of the internal sensors via operators of the form \eqref{eq:pointmeas}, it suffices to consider a time-independent model measurement map $B(p) \in H^2(D)^*$ defined by
$$
B(p)\!:  v \mapsto \dfrac{1}{|S|} \int_{S(p)} v \, {\rm d} x, \qquad p \in \R^d.
$$
For a given location $p \in \R^d$, the sensor $S(p)$ is defined as
$$
S(p) = \{ x + p \, | \, x \in S \},
$$
where the bounded Lipschitz domain  $S = S(0) \subset \R^d$ models the shape (and the orientation) of the sensor. We only consider such $p$ that $S(p) \subset \subset D$, and denote the corresponding open subset of $\R^d$ by $U$. Our aim is to prove that the Fr\'echet derivative of the mapping $U \ni p \mapsto B(p) \in  H^2(D)^*$ in the direction $q \in \R^d$ is given by the functional
\begin{equation}
  \label{eq:Fderivative}
D B(p; q):  v  \mapsto \dfrac{1}{|S|} \int_{S(p)} q \cdot \nabla v \, {\rm d} x = \dfrac{1}{|S|} \int_{\partial S(p)} \nu \cdot q \, v \, {\rm d} S,
\end{equation}
where $\nu \in L^\infty(\partial S(p), \R^d)$ is the exterior unit normal of $\partial S(p)$ and the equality follows from the Gauss divergence theorem. It is obvious that $D B(p; q)$ defines an element of $H^2(D)^*$ for all $p \in U$ and $q\in \R^d$, with the dependence of $DB(p; q)$ on its the latter variable being linear.

For any $v \in \mathcal{C}^\infty(\overline{D})$, $p \in U$ and small enough $q \in \R^d$, we have
\begin{align*}
  \label{eq:mean_value}
  \big(B(p + q) - B(p)\big) v  &=
  \dfrac{1}{|S|} \int_{S(p)} \int_0^1 q \cdot \nabla v(x + s q) \, {\rm d} s \, {\rm d} x =  \dfrac{1}{|S|} \int_0^1 \int_{\partial S(p)} \nu \cdot q \, v(x + s q) \, {\rm d} S_x \, {\rm d} s.
\end{align*}
By resorting to the mean value theorem for integrals, it thus follows that
$$
\big| \big( B(p + q) - B(p) - D B(p; q) \big) v\big |
\leq \dfrac{|q|}{|S|} \int_{\partial S(p)} \big |v(x + \eta_{v,q} q) - v(x)\big | {\rm d} S_x
$$
for some $\eta_{v,q} \in [0,1]$. In particular,
$$
\big| \big(B(p + q) - B(p) - D B(p; q) \big) v \big | \leq |q|^{1+\alpha} \dfrac{|\partial S|}{|S|} \| v \|_{\mathcal{C}^{0,\alpha}(\overline{D})} \leq C |q|^{1+\alpha} \| v \|_{H^2(D)}
$$
with $\alpha = 1/2$, for both $d=2$ and $d=3$, due to a Sobolev embedding theorem \cite{Adams75}. Since $\mathcal{C}^\infty(\overline{D})$ is dense in $H^2(D)$, the above inequality actually holds for all $v \in H^2(D)$. Taking the supremum over $v \in  H^2(D)$ with $\| v \|_{H^2(\Omega)} = 1$, we finally arrive at
$$
\big \| \big(B(p + q) - B(p)\big) - D B(p; q) \big \|_{H^2(D)^*} \leq C |q |^{1+\alpha},
$$
which proves the claim.

\bibliographystyle{plain}
\bibliography{biblio}

\end{document}